\documentclass[11pt]{article}
\usepackage{amsmath}
\usepackage{amsmath,amsthm,amsfonts,amssymb,bm,euscript}
\usepackage{fancyhdr,amscd}
\usepackage{anysize}
\usepackage{graphicx,epsfig}
\usepackage{amsthm,latexsym,amsmath,amssymb}
\usepackage{amssymb,amsmath}
\usepackage{mathrsfs}
\usepackage{graphicx}
\usepackage{color}
\usepackage{mathrsfs}
\usepackage{cite}
\usepackage{indentfirst}
 \usepackage[titletoc]{appendix}
 \usepackage [latin1]{inputenc}
\let\f=\frac
\let\p=\psi

\let\D=\Delta

\let\wt=\widetilde

\def\cA{{\cal A}}

\def\cC{{\cal C}}

\def\cP{{\cal P}}

\def\cS{{\cal S}}

\def\cZ{{\cal Z}}

\def\na{\nabla}

\def\p{\partial}

\def\dv{\mbox{div}}

\def\eqdefa{\buildrel\hbox{\footnotesize def}\over =}

\def\C{\mathop{\bf C\kern 0pt}\nolimits}
\def\DD{\mathop{\bf D\kern 0pt}\nolimits}
\def\K{\mathop{\bf K\kern 0pt}\nolimits}
\def\N{\mathop{\bf N\kern 0pt}\nolimits}
\def\Q{\mathop{\bf Q\kern 0pt}\nolimits}
\def\R{\mathop{\bf R\kern 0pt}\nolimits}

\def\ddq{\dot \Delta_q}

\renewcommand{\div}{\mbox{\rm div}\;\!}

\newcommand{\Dv}{{\rm div}}

\newcommand{\beq}{\begin{equation}}
\newcommand{\eeq}{\end{equation}}
\newcommand{\ben}{\begin{eqnarray}}
\newcommand{\een}{\end{eqnarray}}
\newcommand{\beno}{\begin{eqnarray*}}
\newcommand{\eeno}{\end{eqnarray*}}

\newtheorem{Theorem}{Theorem}[section]
\newtheorem{Definition}[Theorem]{Definition}
\newtheorem{Proposition}[Theorem]{Proposition}
\newtheorem{Lemma}[Theorem]{Lemma}

\newtheorem{Remark}[Theorem]{Remark}

\numberwithin{equation}{section}

\allowdisplaybreaks \numberwithin{equation} {section}

\begin{document}
\title{On the  well-posedness and decay rates of
strong solutions to a  multi-dimensional  non-conservative viscous
compressible two-fluid system
     \thanks {Research supported by the
National Natural Science Foundation of China (11501332,11771043,11371221,11871302), the  Natural Science Foundation of Shandong Province (ZR2015AL007),
 and Young Scholars Research Fund of Shandong University of Technology.}
}
\author{ Fuyi  Xu$^{a\dag}$ \ \
         Meiling  Chi$^{a}$ \ \ Lishan   Liu$^{b,c}$ \ \ Yonghong  Wu$^{c}$  \\[2mm]
 { \small $ ^a$ School of Mathematics and  Statistics, Shandong University of Technology,}\\
  { \small Zibo,    255049,  Shandong,    China}\\
{\small $ ^b$ School  of Mathematics Science, Qufu normal University,} \\ { \small Qufu,  263516, Shandong, China}\\
{\small $ ^c$ Department of Mathematics and Statistics, Curtin University,}\\
{\small Perth, 6845, WA,  Australia}}
         \date{}
         \maketitle
\noindent{\bf Abstract}\ \ \ The present paper deals with
the Cauchy problem of a  multi-dimensional  non-conservative viscous
compressible two-fluid system. We first study the well-posedness  of the model  in spaces with critical regularity indices with respect to the scaling of the associated equations.
In the functional setting as close as possible to the physical energy spaces, we  prove the unique global solvability  of  strong
solutions  close to a stable
equilibrium state. Furthermore,  under a mild additional decay assumption involving only the low frequencies of the  data,  we establish  the  time decay rates  for  the  constructed global
solutions. The proof relies on an application of  Fourier analysis to a complicated  parabolic-hyperbolic system, and on a refined time-weighted inequality.
 \vskip   0.2cm \footnotetext[1]{$^\dag$Corresponding author.}
\vskip   0.2cm \footnotetext[2]{E-mail addresses: zbxufuyi@163.com(F.Xu),\ chimeiling0@163.com(M. Chi),\ mathlls@163.com(L. Liu),\  y.wu@curtin.edu.au(Y. Wu).} \setlength{\baselineskip}{20pt}

\section{Introduction and main results}
\setcounter{section}{1}\setcounter{equation}{0} \ \ \ \ \
It is well known that models of two-phase or multiphase flows are widely applied to  study the  hydrodynamics in industry, for example, in manufacturing, engineering,
and biomedicine, where the fluids under investigation contain more than one component.
In fact, it has been estimated that over half of everything produced in a modern
industrial society depends, to some degree, on a multiphase flow process for its optimum
design and safe operation. In nature, there is a variety of different multiphase
flow phenomena, such as sediment transport, geysers, volcanic eruptions, clouds, and
rain \cite{BDG,GN}. In addition, models of multiphase flows also naturally appear in many contexts
within biology, ranging from tumor biology and anticancer therapies to developmental
biology and plant physiology\cite{MT}. The principles of single-phase flow fluid dynamics
and heat transfer are relatively well understood; however, the thermofluid dynamics
of two-phase flows is an order of magnitude more complicated than that of the
single-phase flow due to the existence of a moving and deformable interface and its
interactions with two phases \cite{NP,NK1,NK2}.

 In this paper, we are concerned with the following  mathematical model of multiphase flow, namely a multi-dimensional non-conservative viscous
compressible two-fluid system $\mathbb{R}^{N}(N\geq2)$:
\begin{align}\label{equ:CTFS}
\left\{
\begin{aligned}
&\alpha^{+}+\alpha^{-}=1,\\
&\p_t(\alpha^{\pm}\rho^{\pm})+\textrm{div}(\alpha^{\pm}\rho^{\pm}u^{\pm})=0, \\
&\p_t(\alpha^{\pm}\rho^{\pm}u^{\pm})+\textrm{div}(\alpha^{\pm}\rho^{\pm}
u^{\pm}\otimes u^{\pm})+\alpha^{\pm}\nabla P^{\pm}(\rho^{\pm})
=\textrm{div}(\alpha^{\pm}\tau^{\pm}), \\
&P^{+}(\rho^{+})=A^{+}(\rho^{+})^{\bar{\gamma}^{+}}=P^{-}(\rho^{-})=A^{-}(\rho^{-})^{\bar{\gamma}^{-}},
\end{aligned}
\right.
\end{align}
where the variable $0\leq\alpha^{+}(x,t)\leq1$ is the volume fraction of fluid $+$ in one of the two gases, and $0\leq\alpha^{-}(x,t)\leq1$  is the  volume fraction of the other fluid$-$. Moreover, $\rho^{\pm}(x,t)\geq0$, $u^{\pm}(x,t)$ and $P^{\pm}(\rho^{\pm})=A^{\pm}(\rho^{\pm})^{\overline{\gamma}^{\pm}}$ are, respectively, the
densities, the velocities, and the two pressure functions of the fluids. It is assumed
that $\overline{\gamma}^{\pm}\geq1,A^{\pm}>0$ are constants. In what follows, we set $A^{+}=A^{-}=1$ without loss of any generality. Also, $\tau^{\pm}$ are the viscous stress tensors
\begin{equation}\label{1.2}
\tau^{\pm}:=2\mu^{\pm}D(u^{\pm})+\lambda^{\pm}\textrm{div}u^{\pm}\textrm{Id},
\end{equation}
where $D(u^{\pm})\eqdefa\frac{\nabla u^{\pm}+\nabla^{t}u^{\pm}}{2}$ stand for the deformation tensor,  the constants  $\mu^{\pm}$ and $\lambda^{\pm}$ are the (given) shear and bulk viscosity coefficients satisfying $\mu^{\pm}>0$ and $\lambda^{\pm}+2\mu^{\pm}>0$.  This model is known as a two-fluid flow system with algebraic closure, and  we refer readers to Refs \cite{BDG,BHL,EWW} for more discussions on the system.

From the viewpoint of partial differential equations,  system
\eqref{equ:CTFS} is a highly nonlinear system coupling between hyperbolic
equations and parabolic equations. As a
matter of fact, there is no diffusion on the mass conservation system satisfying hyperbolic
equations, whereas velocity evolves according to the parabolic equations due to the  viscosity phenomena.
We should point out that the system \eqref{equ:CTFS} includes
 important single phase flow models such as  the compressible Navier-Stokes equations  when one of the two phases volume fraction tends to zero (i.e.,$\alpha^{+}=0$ or $\alpha^{-}=0$ ).  As  an
extremely important system to describe compressible fluids (e.g.,
gas dynamics), the compressible Navier-Stokes equations  have attracted a lot of attention among many analysts and many important results have been developed.  Here we briefly review some  of the most relevant papers
about global well-posedness and  large time behaviors of the  solutions to  the system. Lions \cite{Lions} proved the global existence of weak solutions
for large initial data. However, the question of uniqueness of weak solutions remains open,
even in the two dimensional case.
Matsumura and Nishida \cite{MN1,MN2} first studied the global existence of classical
solutions to the compressible Navier-Stokes equations   for data $(\rho_{0}, u_{0})$ with high regularity order and close to a stable
equilibrium  in the 3D
whole space and obtained the time decay rates based on the $L^{2}$-framework. Later, Ponce \cite{Ponce} established the optimal $L^{p}$ $(2 \leq p\leq \infty)$
decay rates. Applying Fourier analysis  to the linearized homogeneous system and capturing
the dissipation of the hyperbolic component in the solution, Kawashima \cite{Ka1,Ka2} and Shizuta and Kawashima \cite{SK}
developed a general approach to obtain the time-decay of solutions. It is worth mentioning that Li
and Zhang \cite{LZ} obtained the optimal $L^{p}$ time-decay rates for the compressible Navier-Stokes
equations in three dimensions when initial data belong to some space $H^{s} \cap \dot{B}_{1,\infty}^{-s}$ (see  Definition \ref{def2.2} for details)and $s\in[0,1]$.
Guo and Wang \cite{GW} obtained the optimal decay rates for
the compressible Navier-Stokes equations when  the initial data are  close to a stable
equilibrium state in negative Sobolev spaces by using a pure energy method.
In \cite{Dan1}, Danchin first proved   the existence  and uniqueness of the global strong solution for  initial data  close to a stable
equilibrium state  in critical Besov spaces. Later, Danchin \cite{Dan2} further established  the time decay rates of  the  global
solutions  constructed  in \cite{Dan1}.

Since a single phase flow model may be considered as a special case of  a two-phase
flow model in the limit when one of the two phases volume fractions tends to zero, the mathematical structure of the
 two-phase system is much more complex than that in the case of single phase
flow model. So, extending  the currently available results for single phase flow
models  to two-phase models is not an easy task.  Nowadays,  more and more researchers  pay more attention to
 the mathematical problems of the generic two-phase model.  In \cite{BDG},  Bresch et al. first established the existence of global weak solutions to the 3D  generic two-fluid flow model with  capillary pressure
effects in terms of a third order derivative of $\alpha^{\pm}\rho^{\pm}$. Based on
detailed analysis of the Green function to the linearized system and on elaborate energy estimates
to the nonlinear system, the authors in \cite{CWYZ} obtained  global existence of  smooth solutions and the time decay rates to the 3D model where  the initial data are close to an
equilibrium state  in $H^s(\mathbb{R}^{3}_{x})(s\geq3)$ with high Sobolev regularity and belong to $
L^1(\mathbb{R}^{3}_{x})$. More recently, Lai-Wen-Yao \cite{LWY}  studied the vanishing capillarity limit of the smooth solutions  to the 3D model  with  unequal pressure functions.
When a generic two-fluid flow model does not include  capillary pressure
effects, the model reduces to the system \eqref{equ:CTFS}. Bresch-Huang-Li \cite{BHL} extended the result in \cite{BDG} and proved the existence of global weak solutions to \eqref{equ:CTFS} in one space dimension.  In 2016, Evje-Wang-Wen \cite{EWW}   proved the global existence
of strong solutions to the model \eqref{equ:CTFS} with constant viscosity
coefficients and unequal pressure functions by the standard energy method under the condition that the initial data are
close to the constant equilibrium state in $H^{2}(\mathbb{R}^{3}_{x})$ and obtained the optimal decay rates for the constructed
global strong solutions in $L^{2}$-norm  if the initial data belong to $L^{1}$ additionally. However, to the best of our
knowledge, very few results have been established on the  global  well-posedness and the decay rates of
strong solutions to a  multi-dimensional  non-conservative viscous
compressible two-fluid systems in critical regularity framework.  The purpose of this work is to investigate the mathematical properties of  system \eqref{equ:CTFS} in critical regularity framework. More specifically, we address the question of whether available mathematical
results such as the global well-posedness and time decay rate in critical Besov spaces to a single fluid governed by the compressible barotropic
Navier-Stokes equations may be extended to   multi-dimensional  non-conservative viscous
compressible two-fluid system.

First, we will derive
another expression of the pressure gradient in terms of the gradients of $\alpha^{+}\rho^{+}$ and $\alpha^{-}\rho^{-}$
 by using the pressure equilibrium assumption. The method comes from \cite{BDG}. For the convenience of the reader,  we also show some  derivations  in this part.
 The relation between the
pressures of \eqref{equ:CTFS} implies the following differential identities
\begin{equation}\label{1.3}\textrm{d}P^{+}=s_{+}^{2}\textrm{d}\rho^{+},\quad
\textrm{d}P^{-}=s_{-}^{2}\textrm{d}\rho^{-},\quad \hbox{where}\quad
s_{\pm}:=\sqrt{\frac{\textrm{d}P^{\pm}}{\textrm{d}\rho^{\pm}}(\rho^{\pm})}
=\sqrt{\overline{\gamma}^{\pm}\frac{P^{\pm}(\rho^{\pm})}{\rho^{\pm}}},
\end{equation}
where $s_{\pm}$ denote the sound speed of each phase respectively.

Let
\begin{equation}\label{1.4}
R^{\pm}=\alpha^{\pm}\rho^{\pm}.
\end{equation}
Resorting to $\eqref{equ:CTFS}_{1}$,  we have
\begin{equation}\label{1.5}
\textrm{d}\rho^{+}=\frac{1}{\alpha_{+}}(\textrm{d}R^{+}-
\rho^{+}\textrm{d}\alpha^{+}),
~~\textrm{d}\rho^{-}=\frac{1}{\alpha_{-}}(\textrm{d}R^{-}+
\rho^{-}\textrm{d}\alpha^{+}).
\end{equation}
Combining with \eqref{1.4} and \eqref{1.5}, we  conclude that
\begin{equation*}\label{1.6}
\textrm{d}\alpha^{+}=\frac{\alpha^{-}s_{+}^{2}}
{\alpha^{-}\rho^{+}s_{+}^{2}+\alpha^{+}\rho^{-}s_{-}^{2}}\textrm{d}R^{+}
-\frac{\alpha^{+}s_{-}^{2}}
{\alpha^{-}\rho^{+}s_{+}^{2}+\alpha^{+}\rho^{-}s_{-}^{2}}
\textrm{d}R^{-}.
\end{equation*}
Substituting  the above equality into \eqref{1.5}, we obtain
\begin{equation*}
\textrm{d}\rho^{+}=\frac{\rho^{+}\rho^{-}s_{-}^{2}}
{R^{-}(\rho^{+})^{2}s_{+}^{2}+R^{+}(\rho^{-})^{2}s_{-}^{2}}
\Big(\rho^{-}\textrm{d}R^{+}
+\rho^{+}\textrm{d}R^{-}\Big),
\end{equation*}
and
\begin{equation*}
\textrm{d}\rho^{-}=\frac{\rho^{+}\rho^{-}s_{+}^{2}}
{R^{-}(\rho^{+})^{2}s_{+}^{2}+R^{+}(\rho^{-})^{2}s_{-}^{2}}
\Big(\rho^{-}\textrm{d}R^{+}
+\rho^{+}\textrm{d}R^{-}\Big),
\end{equation*}
which give, for the pressure differential $\textrm{d}P^{\pm}$,
\begin{equation*}
\textrm{d}P^{+}=\mathcal{C}^{2}\big(\rho^{-}\textrm{d}R^{+}
+\rho^{+}\textrm{d}R^{-}\big),
\end{equation*}
and
\begin{equation*}
\textrm{d}P^{-}=\mathcal{C}^{2}\big(\rho^{-}\textrm{d}R^{+}
+\rho^{+}\textrm{d}R^{-}\big),
\end{equation*}
where
\begin{equation*}
\mathcal{C}^{2}\eqdefa\frac{s_{-}^{2}s_{+}^{2}}{\alpha^{-}\rho^{+}s_{+}^{2}
+\alpha^{+}\rho^{-}s_{-}^{2}}.
\end{equation*}
Recalling $\alpha^{+}+\alpha^{-}=1$, we get the following identity:
\begin{equation}\label{1.61}
\frac{R^{+}}{\rho^{+}}+\frac{R^{-}}{\rho^{-}}=1,\quad \hbox{and ~therefore~}
\rho^{-}=\frac{R^{-}\rho^{+}}{\rho^{+}-R^{+}}.
\end{equation}
Then it follows from the pressure relation $\eqref{equ:CTFS}_{4}$ that
\begin{equation}\label{1.7}
\varphi(\rho^{+}):=P^{+}(\rho^{+})-P^{-}(\frac{R^{-}\rho^{+}}
{\rho^{+}-R^{+}})=0.
\end{equation}
 Differentiating $\varphi$ with respect to $\rho^{+}$,   we have
$$\varphi^{'}(\rho^{+})=s^{2}_{+}+s^{2}_{-}\frac{R^{-}R^{+}}{(\rho^{+}-R^{+})^{2}}.$$
By the definition of $R^{+}$, it is natural to look for $\rho^{+}$ which belongs to $(R^{+},+\infty).$  Since $\varphi'>0$ in $(R^{+},+\infty)$  for any given $R^{\pm}>0,$ and $\varphi:(R^{+},+\infty)\longmapsto(-\infty,+\infty),$  this determines that $\rho^{+}=\rho^{+}(R^{+},R^{-})\in(R^{+},+\infty)$ is the unique solution of the equation \eqref{1.7}. Due to \eqref{1.5}, \eqref{1.61} and $\eqref{equ:CTFS}_{1}$, $\rho^{-}$ and $\alpha^{\pm}$  are defined as follows:
$$\rho^{-}(R^{+},R^{-})=\frac{R^{-}\rho^{+}(R^{+},R^{-})}{\rho^{+}(R^{+},R^{-})-R^{+}},$$
$$\alpha^{+}(R^{+},R^{-})=\frac{R^{+}}{\rho^{+}(R^{+},R^{-})},$$
$$\alpha^{-}(R^{+},R^{-})=1-\frac{R^{+}}{\rho^{+}(R^{+},R^{-})}=\frac{R^{-}}{\rho^{-}(R^{+},R^{-})}.$$
Based on  the above analysis,  the system \eqref{equ:CTFS} is equivalent to the following form
\begin{align}\label{equ:CTFS1}
\left\{
\begin{aligned}
&\p_tR^{\pm}+\textrm{div}(R^{\pm}u^{\pm})=0, \\
&\p_t(R^{+}u^{+})+\textrm{div}(R^{+}u^{+}\otimes u^{+})+\alpha^{+}\mathcal{C}^{2}[\rho^{-}\nabla R^{+}+\rho^{+}\nabla R^{-}]
\\&\qquad=\textrm{div}\big(\alpha^{+}[\mu^{+}(\nabla u^{+}+\nabla^{t}u^{+})+\lambda^{+}\textrm{div}u^{+}\textrm{Id}]\big), \\
&\p_t(R^{-}u^{-})+\textrm{div}(R^{-}u^{-}\otimes u^{-})+\alpha^{-}\mathcal{C}^{2}[\rho^{-}\nabla R^{+}+\rho^{+}\nabla R^{-}]
\\&\qquad=\textrm{div}\big(\alpha^{-}[\mu^{-}(\nabla u^{-}+\nabla^{t}u^{-})+\lambda^{-}\textrm{div}u^{-}\textrm{Id}]\big).
\end{aligned}
\right.
\end{align}
In this
paper, we are concerned with the Cauchy problem of the system
 \eqref{equ:CTFS1} in $\mathbb{R}_{+}\times \mathbb{R}^N$ subject to the
initial data
\begin{equation}\label{eq:initial data}
(R^{+},\,u^{+},\,R^{-},\,u^{-})(x,t)|_{t=0}=(R^{+}_{0},\,u^{+}_{0},\,R^{-}_{0},\,u^{-}_{0})(x), \quad x\in\mathbb{R}^{N},
\end{equation}
and
$$u^{+}(x,t)\rightarrow0,\quad u^{-}(x,t)\rightarrow0,\quad R^{+}\rightarrow R^{+}_{\infty}>0,\quad R^{-}\rightarrow R^{-}_{\infty}>0, \hbox{~as~} |x|\rightarrow\infty,$$
where $R^{\pm}_{\infty}$ denote the background doping profile, and in the present paper $R^{\pm}_{\infty}$ are taken as $1$ for simplicity.

At this stage, we are going to use scaling considerations for  \eqref{equ:CTFS} to
guess which spaces may be critical. One can check that if $(\alpha^{+}\rho^{+},\,u^{+},\,\alpha^{-}\rho^{-},\,u^{-})$  solves  \eqref{equ:CTFS}, so does $\big((\alpha^{+}\rho^{+})_{\lambda},\,u^{+}_{\lambda}, (\alpha^{-}\rho^{-})_{\lambda},\,u^{-}_{\lambda}\big)$ where:
\begin{equation}
\begin{split}\label{udivc}&(\alpha^{+}\rho^{+})_{\lambda}(t, x)=(\alpha^{+}\rho^{+})(\lambda^{2}t, \lambda x),\quad u^{+}_{\lambda}(t, x)=\lambda u^{+}(\lambda^{2}t, \lambda x),\\ &(\alpha^{-}\rho^{-})_{\lambda}(t, x)=(\alpha^{-}\rho^{-})(\lambda^{2}t, \lambda x),\quad u^{-}_{\lambda}(t, x)=\lambda u^{-}(\lambda^{2}t, \lambda x)
\end{split}
\end{equation}
provided that the pressure laws $P$ have been changed into $\lambda^{2}P$.  This suggests us to choose
initial data $\big((\alpha^{+}\rho^{+})_{0},\,u^{+}_{0},\,(\alpha^{-}\rho^{-})_{0},\,u^{-}_{0}\big)$ in critical  spaces whose norm is invariant for all $\lambda>0$ by the
transformation $\big((\alpha^{+}\rho^{+})_{0},\,u^{+}_{0},\,(\alpha^{-}\rho^{-})_{0},\,u^{-}_{0}\big)\rightarrow \big((\alpha^{+}\rho^{+})_{0}(\lambda\,\cdot),\,\lambda u^{+}_{0}(\lambda\,\cdot),\,(\alpha^{-}\rho^{-})_{0}(\lambda\,\cdot),\,\lambda u^{-}_{0}(\lambda\,\cdot)\big)$. Due to the mixed hyperbolic-parabolic property of the partial differential system \eqref{equ:CTFS}, motivated by Danchin's excellent work in \cite{Dan1}, the different dissipative
mechanisms of low frequencies and high frequencies inspire us to
deal with $(\alpha^{+}\rho^{+}-1,\,u^{+},\,\alpha^{-}\rho^{-}-1,\,u^{-})$ in $\tilde{B}^{\f N2-1,\f N2}_{2,1}\times \dot{B}^{\f N2-1}_{2,1}\times\tilde{B}^{\f N2-1,\f N2}_{2,1}\times \dot{B}^{\f N2-1}_{2,1}$(see  Definitions \ref{def2.2} and \ref{def2.6} for details). However,  we  can not  obtain the  desired bounds directly  in critical regularity framework if   the convection terms are treated as perturbations.
More precisely,  there exists a difficulty coming from the convection terms $u^{\pm}\cdot\na R^{\pm}$ in the transport equations without any diffusion in high frequencies, as  one derivative loss about the function $R^{\pm}$ will appear no matter how smooth is $u^{\pm}$ if they are viewed as  perturbation terms.  To overcome
the difficulty, employing the Littlewood-Paley theory and some commutator estimates,  we shall, as in \cite{Dan1} for  the standard barotropic Navier-Stokes equations, study a complicated hyperbolic-parabolic linear system
including  convection terms and then deduce the smooth effect for $(R^{+}-1,\,u^{+},\,R^{-}-1,\,u^{-})$  in the low frequencies regime and the $L^{1}$ decay on the density $R^{\pm}$ in the high frequencies regime. In particular, based on  the damping effect of $R^{\pm}$, we further exploit the smooth effect for $(u^{+},\,u^{-})$  in the high frequencies regime with
$\nabla R^{\pm}$ being viewed as  perturbation terms  and finally establish a uniform  \emph{priori}  estimate for the complicated system (see the following Lemma \ref{prop3.1} for details). Here,  it should be pointed out that, different from the standard barotropic compressible  Navier-Stokes equations, we need to make more careful analyses to cancel some mixed terms from the  two-phase
flows. Next,  one may wonder  how global strong solutions constructed above look like for large time. Under a suitable additional condition involving only the low frequencies of the  data and in the  $L^{2}$-critical regularity framework,  we  exhibit the  time decay rates  for  the  constructed global strong
solutions. In this part, our main ideas are based on  the low-high frequency decomposition   and  a refined time-weighted energy functional.  In  low frequencies, making good use of
Fourier   localization  analysis to a linearized parabolic-hyperbolic system in order to obtain  smoothing effect of  the Green function  in the low-frequency part  and  avoid some complicate analysis of the Green function (see Lemma \ref{lem6.1}), which
is a $8\times 8 $ matrix. In  high frequencies,  we can deal with  the estimates of the nonlinear terms in the
system employing  the Fourier localization method
and commutator estimates. Finally, in order to close the energy estimates, we exploit some  decay estimates with gain of
regularity  for the high frequencies of $\nabla u^{\pm}$.

Now we state our main results as follows:
\begin{Theorem}\label{th:main1}  Assume that $(R^{+}_{0}-1,\,u^{+}_{0},\,R^{-}_{0}-1,\,u^{-}_{0})\in \tilde{B}^{\f{N}{2}-1,\f{N}{2}}_{2,1}\times
\dot{B}^{\f{N}{2}-1}_{2,1}\times \tilde{B}^{\f{N}{2}-1,\f{N}{2}}_{2,1}\times
\dot{B}^{\f{N}{2}-1}_{2,1}$. Then there exists a constant $\eta>0$ such that if
\begin{equation}
\|(R^{+}_{0}-1,\,R^{+}_{0}-1)\|_{\tilde{B}^{\f{N}{2}-1,\f{N}{2}}_{2,1}}+\|(u^{+}_{0},\,u^{-}_{0})\|_{
\dot{B}^{\f{N}{2}-1}_{2,1}}\leq \eta,  \label{1Hn/2}
\end{equation}%
then the Cauchy  problem \eqref{equ:CTFS1}-\eqref{eq:initial data} admits a unique global solution $(R^{+}-1,\,u^{+},\,R^{-}-1,\,u^{-})$ satisfying that for all $t\geq 0$,
\begin{equation}\label{1.6A}
\begin{split}
\,X(t)\lesssim \|(R^{+}_{0}-1,\,R^{+}_{0}-1)\|_{\tilde{B}^{\f{N}{2}-1,\f{N}{2}}_{2,1}}+\|(u^{+}_{0},\,u^{-}_{0})\|_{
\dot{B}^{\f{N}{2}-1}_{2,1}},
\end{split}
\end{equation}
where
\begin{equation}
\begin{split}
\label{1.6}
\,X(t)&\eqdefa\|(R^{+}-1,\,u^{+},\,R^{-}-1,\,u^{-})\|_{\wt L^\infty_t(\dot B^{\frac{N}{2}-1}_{2,1})}^\ell
+\|(R^{+}-1,\,u^{+},\,R^{-}-1,\,u^{-})\|_{L^1_t(\dot B^{\frac{N}{2}+1}_{2,1})}^\ell\\
&\quad+\|(u^{+},\,u^{-})\|_{\wt L^\infty_t(\dot B^{\frac{N}{2}-1}_{2,1})}^h
+\|(R^{+}-1,\,R^{-}-1)\|_{\wt L^\infty_t(\dot B^{\frac{N}{2}}_{2,1})}^h
\\
&\quad+\|(u^{+},\,u^{-})\|_{L^1_t(\dot B^{\frac{N}{2}+1}_{2,1})}^h
+\|(R^{+}-1,\,R^{-}-1)\|_{L^1_t(\dot B^{\frac{N}{2}}_{2,1})}^h.
\end{split}
\end{equation}
\end{Theorem}
\begin{Theorem}\label{th:decay}
Let the data $(R^{+}_{0}-1,\,u^{+}_{0},\,R^{-}_{0}-1,\,u^{-}_{0})$ satisfy the assumptions of Theorem \ref{th:main1}. Denote $\langle \tau\rangle\eqdefa\sqrt{1+\tau^2}$
and $\alpha\eqdefa\min\{2+\frac{N}{4},\frac{N}{2}+\frac{1}{2}-\varepsilon\}$ with $\varepsilon>0$ arbitrarily small.
 There exists a positive constant $c$ such that if in addition
\begin{equation}\label{eq:D0}
D_0\eqdefa\|(R^{+}_{0}-1,\,u^{+}_{0},\,R^{-}_{0}-1,\,u^{-}_{0})\|^{\ell}_{\dot B^{-\frac N2}_{2,\infty}}\leq c,
\end{equation}
then  the global solution $(R^{+}-1,\,u^{+},\,R^{-}-1,\,u^{-})$
given by Theorem \ref{th:main1} satisfies for all $t\geq0,$
\begin{equation}
\label{1.8}
D(t)\leq C\bigl(D_0+\|(\nabla R^{+}_{0},\,u^{+}_{0},\,\nabla R^{-}_{0},\,u^{-}_{0})\|^h_{\dot B^{\frac N2-1}_{2,1}}\bigr)
\end{equation}
with
\begin{equation}
\begin{split}
\label{1.9}
D(t)&\eqdefa \sup_{s\in(\varepsilon-\frac N2,2]}\|\langle\tau\rangle^{\frac N4+\frac s2}(R^{+}-1,\,u^{+},\,R^{-}-1,\,u^{-})\|_{L^\infty_t(\dot B^s_{2,1})}^\ell
\\ &\quad+\|\langle\tau\rangle^{\alpha}(\nabla R^{+},\,u^{+},\,\nabla R^{-},\,u^{-})\|_{\wt L^\infty_t(\dot B^{\frac N2-1}_{2,1})}^h\\
&\quad+\|\tau(\nabla  u^{+},\,\nabla  u^{-})\|_{\wt L^\infty_t(\dot B^{\frac{N}{2}}_{2,1})}^h.
\end{split}
\end{equation}
\end{Theorem}
\begin{Remark}\label{1.2} In Theorem \ref{th:decay}, we obtain the time decay rates for multi-dimensional  non-conservative viscous
compressible two-fluid system
 \eqref{equ:CTFS}  in critical regularity framework. Additionally,
the regularity  index $s$ can take both negative and nonnegative values, rather than only nonnegative integers, which   improves the classical decay results  in high Sobolev regularity, such as \cite{EWW} when $f=P^{+}(\rho^{+})=P^{-}(\rho^{-})$. In fact, for the solution $(R^{+}-1,\,u^{+},\,R^{-}-1,\,u^{-})$  constructed in  Theorem \ref{th:main1}, applying  to homogeneous Littlewood-Paley decomposition
 for $R^{+}-1$, we have
$$\|\Lambda^s (R^{+}-1)\|_{L^2}\lesssim  \sum_{q\in \mathbb{Z}} \|\Delta_q\Lambda^s (R^{+}-1)\|_{L^2}=\|\Lambda^s (R^{+}-1)\|_{\dot B^0_{2,1}}.$$
Based on  Bernstein's inequalities and the low-high frequencies decomposition, we may write
$$
\sup_{t\in[0,T]} \langle t\rangle^{\frac N4+\frac s2}\|\Lambda^s (R^{+}-1)\|_{\dot B^0_{2,1}}\lesssim
 \|\langle t\rangle^{\frac N4+\frac s2}(R^{+}-1)\|_{L^\infty_T(\dot B^s_{2,1})}^\ell
 +  \|\langle t\rangle^{\frac N4+\frac s2}(R^{+}-1)\|_{L^\infty_T(\dot B^s_{2,1})}^h.
$$
If follows from Inequality \eqref{1.8} and  definitions of $D(t)$ and $\alpha$ that
$$
  \|\langle t\rangle^{\frac N4+\frac s2}(R^{+}-1)\|_{L^\infty_T(\dot B^s_{2,1})}^\ell\lesssim D_0+\|(\nabla R^{+}_{0},\,u^{+}_{0},\,\nabla R^{-}_{0},\,u^{-}_{0})\|^h_{\dot B^{\frac N2-1}_{2,1}}
\quad\hbox{if }\ \  -N/2<s\leq 2
$$
and that, because we have $\alpha\geq\frac N4+\frac s2$ for all $s\leq \min\{2,N/2 \},$
$$
  \|\langle t\rangle^{\frac N4+\frac s2}(R^{+}-1)\|_{L^\infty_T(\dot B^s_{2,1})}^h
  \lesssim D_0+\|(\nabla R^{+}_{0},\,u^{+}_{0},\,\nabla R^{-}_{0},\,u^{-}_{0})\|^h_{\dot B^{\frac N2-1}_{2,1}}\quad\hbox{if }\ \  s\leq \min\{2,N/2 \}.$$
This yields the following desired result for $R^{+}-1$
$$\|\Lambda^{s}(R^{+}-1)\|_{L^2}
\leq C\bigl(D_0+\|(\nabla R^{+}_{0},\,u^{+}_{0},\,\nabla R^{-}_{0},\,u^{-}_{0})\|^h_{\dot B^{\frac N2-1}_{2,1}}\bigr)\langle t\rangle^{-\frac N4-\frac s2}
\quad  \hbox{ if } \ \ -N/2<s\leq \min\{2,N/2 \},$$  where the fractional derivative
 operator $\Lambda^{\ell}$ is defined by $\Lambda^{\ell}f\triangleq\mathcal{F}^{-1}(|\cdot|^{\ell}\mathcal{F}f)$. Similarly, we have
$$\displaylines{
\|\Lambda^{s}u^{+}\|_{L^2}
\leq C\bigl(D_0+\|(\nabla R^{+}_{0},\,u^{+}_{0},\,\nabla R^{-}_{0},\,u^{-}_{0})\|^h_{\dot B^{\frac N2-1}_{2,1}}\bigr)\langle t\rangle^{-\frac N4-\frac s2}
\quad \hbox{ if } \ \ -N/2<s\leq \min\{2,N/2-1 \},\cr
\|\Lambda^{s}(R^{-}-1)\|_{L^2}
\leq C\bigl(D_0+\|(\nabla R^{+}_{0},\,u^{+}_{0},\,\nabla R^{-}_{0},\,u^{-}_{0})\|^h_{\dot B^{\frac N2-1}_{2,1}}\bigr)\langle t\rangle^{-\frac N4-\frac s2}
\quad  \hbox{ if } \ \ -N/2<s\leq \min\{2,N/2 \},\cr
\|\Lambda^{s}u^{-}\|_{L^2}
\leq C\bigl(D_0+\|(\nabla R^{+}_{0},\,u^{+}_{0},\,\nabla R^{-}_{0},\,u^{-}_{0})\|^h_{\dot B^{\frac N2-1}_{2,1}}\bigr)\langle t\rangle^{-\frac N4-\frac s2}
\quad \hbox{ if } \ \ -N/2<s\leq \min\{2,N/2-1 \},}
$$
In particular, taking $s=0$ leads back to  the standard optimal   $L^{1}$-$L^{2}$ decay rate of $(R^{+}-1,\,u^{+},\,R^{-}-1,\,u^{-})$ as in \cite{EWW} when $N=3$.
 \end{Remark}
\begin{Remark}\label{1.2}Due to the embedding $L^{1}(\mathbb{R}^3)\hookrightarrow \dot{B}^{-\frac 32}_{2,\infty}(\mathbb{R}^3)$,  our results in Theorem \ref{th:decay} extend the known conclusions in \cite{EWW}. In particular, our condition involves only the low frequencies of the data and is based on the  $L^{2}(\mathbb{R}^{3}_{x})$-norm framework.  In particular, the decay rates of strong solutions is
in the so-called critical Besov spaces in any dimension $N \geq 2$ and the dimension of space is more extensive and is not limited to $N=3$.
\end{Remark}
\begin{Remark} In this paper,  we can not deal with the case of the model \eqref{equ:CTFS} with unequal pressure functions as in \cite{EWW} since we take advantages of the   symmetrizers  methods in our process of the proof.
\end{Remark}

\noindent{\bf Notations.}  We assume $C$ be
a positive generic constant throughout this paper that may vary at
different places and
denote  $A\le CB$ by  $A\lesssim B$.
We shall also use the  following notations
 $$z^\ell\eqdefa\sum_{j\leq k_0}\dot{\Delta}_{j}z\quad\hbox{and}\quad z^h\eqdefa z-z^\ell, \quad\hbox{for some } j_0.$$
$$\|z\|^\ell_{\dot B^s_{2,1}}\eqdefa \sum_{j\leq k_0}2^{js}\|\dot{\Delta}_{j}z\|_{L^{2}}\quad\hbox{and}\quad \|z\|^h_{\dot B^s_{2,1}}\eqdefa \sum_{j\geq k_0}2^{js}\|\dot{\Delta}_{j}z\|_{L^{2}}, \quad\hbox{for some } j_0.$$
Noting the small overlap between low and high frequencies,  we have
$$\|z^\ell\|_{\dot B^s_{2,1}}\lesssim \|z\|^\ell_{\dot B^s_{2,1}}\quad\hbox{and}\quad \|z^h\|_{\dot B^s_{2,1}}\lesssim \|z\|^h_{\dot B^s_{2,1}}.$$
\par

\section{ Littlewood-Paley theory and some useful lemmas }
\ \ \ \ \ Let us introduce the Littlewood-Paley decomposition.
Choose a radial function  $\varphi \in {\cS}(\mathbb{R}^N)$
supported in ${\cC}=\{\xi\in\mathbb{R}^N,\,
\frac{3}{4}\le|\xi|\le\frac{8}{3}\}$ such that \beno \sum_{q\in
\mathbb{Z}}\varphi(2^{-q}\xi)=1 \quad \textrm{for all}\,\,\xi\neq 0.
\eeno The homogeneous frequency localization operators
$\dot{\Delta}_q$ and $\dot{S}_q$ are defined by
\begin{align}
\dot{\Delta}_qf=\varphi(2^{-q}D)f,\quad \dot{S}_qf=\sum_{k\le
q-1}\dot{\Delta}_kf\quad\mbox{for}\quad q\in \mathbb{Z}. \nonumber
\end{align}
With our choice of $\varphi$, one can easily verify that
\begin{equation*}\begin{split}
&\dot{\Delta}_q\dot{\Delta}_kf=0\quad \textrm{if}\quad|q-k|\ge
2\quad \textrm{and}
\quad \\
&\dot{\Delta}_q(\dot{S}_{k-1}f\dot{\Delta}_k f)=0\quad
\textrm{if}\quad|q-k|\ge 5.
\end{split}\end{equation*}
 We denote the space ${\cZ'}(\mathbb{R}^N)$ by the dual space of
${\cZ}(\mathbb{R}^N)=\{f\in {\cS}(\mathbb{R}^N);\,D^\alpha
\hat{f}(0)=0; \forall\alpha\in\mathbb{ N}^N \,\mbox{multi-index}\}.$
It also can be identified by the quotient space of
${\cS'}(\mathbb{R}^N)/{\cP}$ with the polynomials space ${\cP}$. The
formal equality \beno f=\sum_{q\in\mathbb{Z}}\dot{\Delta}_qf \eeno
holds true for $f\in {\cZ'}(\mathbb{R}^N)$ and is called the
homogeneous Littlewood-Paley decomposition.

The following Bernstein's inequalities will be frequently used.
\begin{Lemma}\cite{Che-book}\label{Lem:Bernstein}
Let $1\le p_{1}\le p_{2}\le+\infty$. Assume that $f\in L^{p_{1}}(\mathbb{R}^N)$,
then for any $\gamma\in(\mathbb{N}\cup\{0\})^N$, there exist
constants $C_1$, $C_2$ independent of $f$, $q$ such that \beno
&&{\rm supp}\hat f\subseteq \{|\xi|\le A_02^{q}\}\Rightarrow
\|\partial^\gamma f\|_{p_2}\le C_12^{q{|\gamma|}+q
N(\frac{1}{p_1}-\frac{1}{p_2})}\|f\|_{p_1},
\\
&&{\rm supp}\hat f\subseteq \{A_12^{q}\le|\xi|\le
A_22^{q}\}\Rightarrow \|f\|_{p_1}\le
C_22^{-q|\gamma|}\sup_{|\beta|=|\gamma|}\|\partial^\beta f\|_{p_1}.
\eeno
\end{Lemma}
\begin{Definition}\label{def2.2} Let $s\in \mathbb{R}$, $1\le p,
r\le+\infty$. The homogeneous Besov space $\dot{B}^{s}_{p,r}$ is
defined by
$$\dot{B}^{s}_{p,r}=\Big\{f\in {\cZ'}(\mathbb{R}^N):\,\|f\|_{\dot{B}^{s}_{p,r}}<+\infty\Big\},$$
where \beno \|f\|_{\dot{B}^{s}_{p,r}}\eqdefa \Big\|2^{qs}
\|\dot{\Delta}_qf(t)\|_{p}\Big\|_{\ell^r}.\eeno
\end{Definition}
\begin{Remark}\label{2.3}
Some properties about the  Besov spaces are as follows
\begin{itemize}
\item\,\, Derivation: $$\|f\|_{\dot{B}^{s}_{2,1}}\approx\|\nabla f\|_{\dot{B}^{s-1}_{2,1}};$$
\item\,\, Algebraic properties: for $s > 0$, $\dot{B}^{s}_{2,1} \cap L^{\infty}
$ is an algebra;
\item\,\,Interpolation: for
$s_1, s_2\in\mathbb{R}$ and $\theta\in[0,1]$,
we have $$\|f\|_{\dot{B}^{\theta s_1+(1-\theta)s_2}_{2,1}}\le
\|f\|^\theta_{\dot{B}^{s_1}_{2,1}}\|f\|^{(1-\theta)}_{\dot{B}^{s_2}_{2,1}}.$$
\end{itemize}
\end{Remark}
\begin{Definition} Let $s\in \mathbb{R}$, $1\le p,\rho,
r\le+\infty$. The homogeneous space-time  Besov space
$L^\rho_{T}(\dot{B}_{p,  r}^s)$ is defined by
$$L^\rho_{T}(\dot{B}_{p,  r}^s)=\Big\{f\in \mathbb{R}_{+}\times{\cZ'}(\mathbb{R}^N):\,\|f\|_{L^\rho_{T}(\dot{B}_{p,  r}^s)} <+\infty\Big\},$$ where $$
\|f\|_{L^\rho_{T}(\dot{B}_{p,  r}^s)}\eqdefa \Big\|\big\| 2^{qs}
\|\dot{\Delta}_q f\|_{L^p}\big\|_{\ell ^{r}}\Big\|_{L^\rho_{T}}.$$
\end{Definition}
We next introduce the Besov-Chemin-Lerner space
$\widetilde{L}^q_T(\dot{B}^{s}_{p,r})$ which is initiated in
\cite{Che-Ler}.
\begin{Definition}Let $s\in \mathbb{R}$, $1\le
p,q,r\le+\infty$, $0<T\le+\infty$. The space
$\widetilde{L}^q_T(\dot{B}^s_{p,r})$ is defined by
$$\widetilde{L}^q_T(\dot{B}^s_{p,r})=\Big\{f\in \mathbb{R}_{+}\times{\cZ'}(\mathbb{R}^N):\,\|f\|_{\widetilde{L}^q_T(\dot{B}^{s}_{p,r})}<+\infty\Big\},$$
where
$$\|f\|_{\widetilde{L}^q_T(\dot{B}^{s}_{p,r})}\eqdefa \Bigl\|2^{qs}
\|\dot{\Delta}_qf(t)\|_{L^q(0,T;L^p)}\Bigr\|_{\ell^r}.$$
\end{Definition}
Obviously, $
\widetilde{L}^1_T(\dot{B}^s_{p,1})=L^1_T(\dot{B}^s_{p,1}). $ By  a
direct application  of  Minkowski's inequality, we have the
following relations between these spaces
\begin{equation*}
L^\rho_{T}(\dot{B}_{p,r}^s)\hookrightarrow\widetilde
L^\rho_{T}(\dot{B}_{p,r}^s),\,\textnormal{if}\quad  r\geq
\rho,\end{equation*}
\begin{equation*}
\widetilde L^\rho_{T}(\dot{B}_{p,r}^s)\hookrightarrow
L^\rho_{T}(\dot{B}_{p,r}^s),\, \textnormal{if}\quad \rho\geq r.
\end{equation*}
To deal with functions with different regularities for high
frequencies and low frequencies, motivated by \cite{Dan2,Dan1}, it
is more effective to work in \textit{hybrid Besov spaces}. We remark
that using hybrid Besov spaces has been crucial for proving global
well-posedness for compressible systems in critical spaces (see
\cite{CMZ2,Dan2,Dan1}).
\begin{Definition}\label{def2.6}
Let $s,t\in \mathbb{R}$. We set
$$\|f\|_{\tilde{B}^{s,t}_{2,1}}=\sum_{q\le
0}2^{qs}\|\dot{\D}_qf\|_{L^2}+\sum_{q>0}2^{qt}\|\dot{\D}_qf\|_{L^2}.$$
For $m=-\left[\f{N}{2}+1-s\right]$, we define
\begin{eqnarray}
\tilde{B}^{s,t}_{2,1}&=&\left\{f\in\mathcal{S}'(\mathbb{R}^N):
\|f\|_{\tilde{B}^{s,t}_{2,1}}<\infty\right\},\qquad \textrm{
if}\ m<0,\nonumber\\
\tilde{B}^{s,t}_{2,1}&=&\left\{f\in\mathcal{S}'(\mathbb{R}^N)/\mathcal{P}_m:
\|f\|_{\tilde{B}^{s,t}_{2,1}}<\infty\right\},\ \textrm{ if }m\ge
0.\nonumber
\end{eqnarray}
\end{Definition}
\begin{Remark}\label{2.5}
Some properties about the hybrid Besov spaces are as follows
\begin{itemize}
\item\,\, $\tilde{B}^{s,s}_{2,1}=\dot{B}^s_{2,1}$;
\item\,\, If $s\le t$, then $\tilde{B}^{s,t}_{2,1}=\dot{B}^s_{2,1}\cap \dot{B}^t_{2,1}$. Otherwise,
$\tilde{B}^{s,t}_{2,1}=\dot{B}^s_{2,1}+\dot{B}^t_{2,1}$. In
particular, $\tilde{B}^{s,\f{N}{2}}_{2,1}\hookrightarrow L^\infty$
as $s\le \f{N}{2}$;
\item\,\,Interpolation: for
$s_1, s_2, t_1, t_2\in\mathbb{R}$ and $\theta\in[0,1]$,
we have $$\|f\|_{\tilde{B}^{\theta s_1+(1-\theta)s_2,\, \theta t_1+(1-\theta)t_2}_{2,1}}\le
\|f\|^\theta_{\tilde{B}^{s_1,t_1}_{2,1}}\|f\|^{(1-\theta)}_{\tilde{B}^{s_2,t_2}_{2,1}};$$
\item\,\, If $s_1\le s_2$ and $t_1\ge t_2$, then
$\tilde{B}^{s_1,t_1}_{2,1}\hookrightarrow\tilde{B}^{s_2,t_2}_{2,1}$.
\end{itemize}
\end{Remark}
We have the following properties for the product in  Besov spaces and  hybrid Besov spaces.
\begin{Proposition}\cite{Dan5}\label{p26}
 For all $1\leq r,p, p_1, p_2\leq+\infty$,  there exists a positive universal
 constant such that
$$\|fg\|_{\dot{B}^{s}_{p,r}}\lesssim
\|f\|_{L^\infty}\|g\|_{\dot{B}^{s}_{p,r}}+\|g\|_{L^\infty}\|f\|_{\dot{B}^{s}_{p,r}},
\quad \text{if}\quad s>0;$$
$$\|fg\|_{\dot{B}^{s_1+s_2-\frac{N}{p}}_{p,r}}\lesssim
\|f\|_{\dot{B}^{s_1}_{p,r}}\|g\|_{\dot{B}^{s_2}_{p,\infty}}, \quad
\text{if}\quad s_1,s_2<\frac{N}{p},\quad \text{and}\quad
s_1+s_2>0;$$
$$\|fg\|_{\dot{B}^{s}_{p,r}}\lesssim
\|f\|_{\dot{B}^{s}_{p,r}}\|g\|_{\dot{B}^{\frac{N}{p}}_{p,\infty}\cap
L^{\infty}}, \quad \text{if}\quad |s|<\frac{N}{p};$$
$$\|fg\|_{\dot{B}^s_{2,1}}\lesssim \|f\|_{\dot{B}^{N/2}_{2,1}}\|g\|_{\dot{B}^s_{2,1}}, \quad \text{if}\quad
s\in (-N/2,N/2].$$
\end{Proposition}
\begin{Proposition}\cite{Dan3}\label{p25}
 For all $s_1, s_2>0$,  there exists a positive universal
 constant such that
$$\|fg\|_{\tilde{B}^{s_1,s_2}_{2,1}}\lesssim
\|f\|_{L^\infty}\|g\|_{\tilde{B}^{s_1,s_2}_{2,1}}+\|g\|_{L^\infty}\|f\|_{\tilde{B}^{s_1,s_2}_{2,1}}.$$
For all $s_1, s_2\le\f{N}{2}$ such that $\min\{s_1+t_1,
s_2+t_2\}>0$,  there exists a positive universal
 constant such that
$$\|fg\|_{\tilde{B}^{s_1+t_1-\f{N}{2},s_2+t_2-\f{N}{2}}_{2,1}}\lesssim
\|f\|_{\tilde{B}^{s_1,s_2}_{2,1}}\|g\|_{\tilde{B}^{t_1,t_2}_{2,1}}.$$
\end{Proposition}
For the composition of the  binary functions, we have the following estimates.
\begin{Proposition}\cite{Xu-Chi}\label{p27}
\mbox{(i) } Let $s>0$, $t\geq0$, $1\le p, q,r\le \infty$ and $(f,g)\in \big(\tilde{L}_{t}^{q}(\dot{B}^{s}_{p,r})\cap L_{t}^{\infty}(L^{\infty})\big)\times\big(\tilde{L}_{t}^{q}(\dot{B}^{s}_{p,r})\cap L_{t}^{\infty}(L^{\infty})\big)$. If $F\in W_{loc}^{[s]+2,\infty}(\mathbb{R})\times W_{loc}^{[s]+2,\infty}(\mathbb{R})$ with
$F(0,0)=0$, then $F(f,g)\in \tilde{L}_{t}^{q}(\dot{B}^{s}_{p,r})$. Moreover,  there exists a
 $C$ depending only on $s,p,N$ and $F$
such that
\begin{equation}\label{eq:2.1}
\|F(f,g)\|_{\tilde{L}_{t}^{q}(\dot{B}^{s}_{p,r})}
\leq C\big(1+\|f\|_{L_{t}^{\infty}(L^{\infty})}\|g\|_{L_{t}^{\infty}(L^{\infty})}\big)^{[s]+1}
\|(f,g)\|_{\tilde{L}_{t}^{q}(\dot{B}^{s}_{p,r})}.
\end{equation}
\mbox{ (ii)} If $(f_1,g_1)\in \tilde{L}_{t}^{\infty}(\dot{B}^{\frac{N}{p}}_{p,1})\times\tilde{L}_{t}^{\infty}(\dot{B}^{\frac{N}{p}}_{p,1})$ and $(f_2,g_2)\in \tilde{L}_{t}^{\infty}(\dot{B}^{\frac{N}{p}}_{p,1})\times\tilde{L}_{t}^{\infty}(\dot{B}^{\frac{N}{p}}_{p,1})$, and  $(f_2-f_1,g_2-g_1)$ belongs to $\tilde{L}_{t}^{q}(\dot{B}^{s}_{p,r})\times \tilde{L}_{t}^{q}(\dot{B}^{s}_{p,r})$ with $s\in (-\frac{N}{p},\frac{N}{p}]$ . If $F\in W_{loc}^{[\frac{N}{p}]+2,\infty}(\mathbb{R})\times W_{loc}^{[\frac{N}{p}]+2,\infty}(\mathbb{R})$ with
$\partial_{1}F(0,0)=0$ and $\partial_{2}F(0,0)=0$. Then, there exists a
 $C$ depending only on $s, p, N$ and $F$
such that
\begin{equation}\label{eq:2.3}\begin{split}
&\|F(f_2,g_2)-F(f_1,g_1)\|_{\tilde{L}_{t}^{q}(\dot{B}^{s}_{p,r})}
\\&\leq C\Big(1+\|(f_1,f_2)\|_{L_{t}^{\infty}(L^{\infty})}\|(g_1,g_2)\|_{L_{t}^{\infty}(L^{\infty})}\Big)^{[\frac{N}{p}]+1} \Big(\|(f_1,g_1)\|_{\tilde{L}_{t}^{\infty}(\dot{B}^{\frac{N}{p}}_{p,1})} +\|(f_2,g_2)\|_{\tilde{L}_{t}^{\infty}(\dot{B}^{\frac{N}{p}}_{p,1})}\Big)
\\&\quad\times\Big(\|f_2-f_1\|_{\tilde{L}_{t}^{q}(\dot{B}^{s}_{p,r})}+\|g_2-g_1\|_{\tilde{L}_{t}^{q}(\dot{B}^{s}_{p,r})}\Big).
\end{split}
\end{equation}
\end{Proposition}
Throughout this paper, the following estimates for the convection
terms arising in the linearized systems will be used frequently.
\begin{Proposition}\cite{Dan1}\label{p29}
    Let $F$ be an homogeneous smooth function of degree $m$. Suppose that $-N/2<s_1,t_1,s_2,t_2\leq 1+N/2$. Then,  the following two estimates hold
    \begin{align*}
    &|(F(D)\dot{\Delta}_{q}(v\cdot \nabla a),F(D)\dot{\Delta}_{q} a)|\\
    &\qquad\qquad\leq C \gamma_q2^{-q(\phi^{s_1,s_2}(q)-m)}\|v\|_{\dot{B}^{\frac{N}{2}+1}_{2,1}} \|a\|_{\tilde{B}^{s_,s_2}_{2,1}}\|F(D)\dot{\Delta}_{q} a\|_{2},\\
    &|(F(D)\dot{\Delta}_{q}(v\cdot\nabla a),\dot{\Delta}_{q} b)+(\dot{\Delta}_{q}(v\cdot\nabla b),F(D)\dot{\Delta}_{q} a)|\\
    &\qquad\qquad\leq C \gamma_q\|v\|_{\dot{B}^{\frac{N}{2}+1}_{2,1}}\times\big( 2^{-q\phi^{t_1,t_2}(q)}\|F(D)\dot{\Delta}_{q} a\|_{2}\|b\|_{\tilde{B}^{t_,t_2}_{2,1}}\\
    &\qquad\qquad\qquad+2^{-q(\phi^{s_1,s_2}(q)-m)}\|a\|_{\tilde{B}^{s_,s_2}_{2,1}}
    \|\dot{\Delta}_{q} b\|_{2}\big),
    \end{align*}
where $(\cdot,\cdot)$ denotes the $2$-inner product,  $\sum_{q\in
\mathbb{Z}} \gamma_q\leq 1$ and the operator $F(D)$ is defined by
$F(D)f:=\mathcal{F}^{-1} F(\xi)\mathcal{F} f$,
$\phi^{\alpha,\beta}(q)$ is the following characteristic function on
$\mathbb{Z}$
\begin{align*}
     \phi^{\alpha,\beta}(q)=\left\{
     \begin{array}{ll}
        \alpha,\quad &\text{if }\quad  q\leq 0,\\
        \beta, \quad & \text{if }\quad q\geq 1.
     \end{array}
     \right.
\end{align*}
\end{Proposition}
\begin{Proposition}\cite{Dan2}
\label{Pro:1}
Let $1\leq p,p_{1}\leq\infty$, and $\sigma\in\mathbb{R}$.
There exists a constant $C>0$ depending only on $\sigma$ such that for all $q\in \mathbb{Z}$, we have
$$\|[v\cdot\nabla, \partial_{\ell}\dot{\Delta}_{q}]a\|_{L^{p}}
\leq Cc_{q}2^{-q(\sigma-1)}
\|\nabla v\|_{\dot B^{\frac{N}{p_{1}}}_{p_{1},1}}
\|\nabla a\|_{\dot B^{\sigma-1}_{p,1}}, \quad \text{for}\quad -\min(\frac{N}{p_{1}},\frac{N}{p^{\prime}})<\sigma
\leq1+\min(\frac{N}{p_{1}},\frac{N}{p}),$$
$$\|[v\cdot\nabla,\dot{\Delta}_{q}]a\|_{L^{p}}
\leq Cc_{q}2^{-q\sigma}
\|\nabla v\|_{\dot B^{\frac{N}{p_{1}}}_{p_{1},\infty}\cap L^{\infty}}
\|a\|_{\dot B^{\sigma}_{p,1}}, \quad \text{for}\quad -\min(\frac{N}{p_{1}},\frac{N}{p^{\prime}})<\sigma
<1+\frac{N}{p_{1}},$$
where the commutator $[\cdot,\cdot]$ is defined by $[f,g]=fg-gf$ and $(c_{j})_{j\in \mathbb{Z}}$ denotes a
sequence such that $\sum_{q\in
\mathbb{Z}} c_{q}\leq 1$.
\end{Proposition}
\begin{Proposition}\cite{RD}
\label{Pro:3}
Assume $\mu>0$, $\sigma\in\mathbb{R}, (p,r)\in[1,\infty]^{2}$ and $1\leq \rho_{2}\leq \rho_{1}\leq\infty$. Let $u$ satisfy
\begin{align}\label{eq:heat}
\left\{
\begin{aligned}
&\partial_{t}u-\mu\Delta u=f,\\
&u\mid_{t=0}=u_{0}. \end{aligned} \right.
\end{align}
Then for all $T>0$ the following a priori estimate is fulfilled
\begin{equation}\label{eq:heat1}\mu^{\frac{1}{\rho_{1}}}\|u\|_{\wt L^{\rho_{1}}_T(\dot B^{\sigma+\frac 2\rho_{1}}_{p,r})}\lesssim\|u_{0}\|_{\dot B^{\sigma}_{p,r}}
+\mu^{\frac{1}{\rho_{2}}-1}\|f\|_{\wt L^{\rho_{2}}_T(\dot B^{\sigma-2+\frac 2\rho_{2}}_{p,r})}.
\end{equation}
\end{Proposition}
\begin{Remark}\label{2.14}
The solutions to the following Lam\'{e} system
\begin{align*}\label{eq:heat}
\left\{
\begin{aligned}
&\partial_tu-\cA u=f,\\
&u\mid_{t=0}=u_{0},
\end{aligned} \right.
\end{align*}
also satisfy  \eqref{eq:heat1}. Here, $\cA u=\mu\Delta u+(\lambda+\mu)\nabla\Dv u$.
\end{Remark}
\begin{Proposition}\cite{Dan3}\label{Prop:transport}
Let $s\in \big(-N\min(\frac1p,\frac1{p'}), 1+\frac Np\big)$, $1\le
p,r\le+\infty$, and $s=1+\frac Np$, if $r=1$. Let $v$ be a vector
field such that $\nabla v\in L^1_T(\dot{B}^{\frac{N}{p}}_{p,r}\cap
L^\infty)$. Assume that $u_0\in \dot{B}^{s}_{p,r},$ $g\in
L^1_T(\dot{B}^{s}_{p,r})$ and $u$ is the solution of the following transport equation
\begin{equation}\label{equ:transport}
\left\{
\begin{aligned}{}
&\partial_t u+v\cdot \nabla u =g,\\
&u\mid_{t=0}=u_{0}.
\end{aligned}
\right.
\end{equation}
Then there holds for $t\in[0,T]$,
\begin{equation}\label{equ:transport-estimate}
\|u\|_{\widetilde{L}^\infty_t(\dot{B}^{s}_{p,r})}\le e^{CV(t)}\Big(
\|u_0\|_{\dot{B}^{s}_{p,r}}+\int_0^t
e^{-CV(\tau)}\|g(\tau)\|_{\dot{B}^{s}_{p,r}}d\tau\Big),
\end{equation}
where
$V(t)\eqdefa\int_0^t\|\nabla
v(\tau)\|_{\dot{B}^{\frac{N}{p}}_{p,r}\cap L^\infty}d\tau.$ If
$r<+\infty$, then $f$ belongs to $\mathcal{C}([0,T];
\dot{B}^s_{p,r})$.
\end{Proposition}
We finish this subsection by listing an elementary but useful
inequality.
\begin{Lemma}\cite{MN2}\label{lemma2.13}\quad  Let $r_1,r_2>0$ satisfy $\max\{r_1,r_2\}>1$. Then
$$\int_0^t(1+t-\tau)^{-r_1}(1+\tau)^{-r_2}d\tau\leq C(r_1,r_2)(1+t)^{-\min\{r_1,r_2\}}.$$
\end{Lemma}
\par
\section{Reformulation of the  System \eqref{equ:CTFS1} and and {A priori} estimates for the linearized  system with convection terms}
\subsection{Reformulation of the  System \eqref{equ:CTFS1} }
\ \ \ \ \  To make it more convenient to study, we need some reformulations of \eqref{equ:CTFS1}. More precisely, taking a change of variables by
$$c^{\pm}=R^{\pm}-1.$$
Then, the system \eqref{equ:CTFS1} can be rewritten  as
\begin{equation}\label{3.1}
\left\{
\begin{aligned}{}
&\p_tc^{+}+\textrm{div}u^{+}=H_{1},\\
&\p_t{u}^{+}+\beta_{1}\nabla c^{+}
+\beta_{2}\nabla c^{-}-\nu_{1}^{+}\Delta u^{+}
-\nu_{2}^{+}\nabla\textrm{div}u^{+}=H_{2},
\\&\p_tc^{-}+\textrm{div}u^{-}=H_{3},
\\&\p_tu^{-}+\beta_{3}\nabla c^{+}
+\beta_{4}\nabla c^{-}-\nu_{1}^{-}\Delta u^{-}
-\nu_{2}^{-}\nabla\textrm{div}u^{-}=H_{4},
\end{aligned}
\right.
\end{equation}
with initial data
\begin{equation}\label{gama1}(c^{+},\,u^{+},\,c^{-},\,u^{-})|_{t=0}=(c^{+}_{0},\,u^{+}_{0},\,c^{-}_{0},\,u^{-}_{0}),\end{equation}
where $\beta_{1}=\frac{\mathcal{C}^{2}(1,1)\rho^{-}(1,1)}{\rho^{+}(1,1)},\quad
\beta_{2}=\beta_{3}=\mathcal{C}^{2}(1,1),\quad
\beta_{4}=\frac{\mathcal{C}^{2}(1,1)\rho^{+}(1,1)}{\rho^{-}(1,1)},\quad
\nu_{1}^{\pm}=\frac{\mu^{\pm}}{\rho^{\pm}(1,1)},\quad
\nu_{2}^{\pm}=\frac{\mu^{\pm}+\lambda^{\pm}}{\rho^{\pm}(1,1)}$
and the source terms are
\begin{align}\label{3.2}
H_{1}&=-\textrm{div}(c^{+}u^{+}),\\
\label{3.3}
H_{2}^{i}&=-g_{+}(c^{+},c^{-})\partial_{i}c^{+}
-\tilde{g}_{+}(c^{+},c^{-})\partial_{i}c^{-}-(u^{+}\cdot\nabla)u_{i}^{+}\nonumber\\
&\quad+\mu^{+}h_{+}(c^{+},c^{-})\partial_{j}c^{+}\partial_{j}u^{+}_{i}
+\mu^{+}k_{+}(c^{+},c^{-})\partial_{j}c^{-}\partial_{j}u^{+}_{i}\nonumber\\
&\quad+\mu^{+}h_{+}(c^{+},c^{-})\partial_{j}c^{+}\partial_{i}u^{+}_{j}
+\mu^{+}k_{+}(c^{+},c^{-})\partial_{j}c^{-}\partial_{i}u^{+}_{j}\\
&\quad+\lambda^{+}h_{+}(c^{+},c^{-})\partial_{i}c^{+}\partial_{j}u^{+}_{j}
+\lambda^{+}k_{+}(c^{+},c^{-})\partial_{i}c^{-}\partial_{j}u^{+}_{j}\nonumber\\
&\quad+\mu^{+}l_{+}(c^{+},c^{-})\partial_{j}^{2}u_{i}^{+}+(\mu^{+}+\lambda^{+})l_{+}
(c^{+},c^{-})\partial_{i}\partial_{j}u^{+}_{j},\quad  i,j \in \{1,2,\cdots N\},\nonumber\\
\label{3.4}
H_{3}&=-\textrm{div}(c^{-}u^{-}),\\
\label{3.5}
H_{4}^{i}&=-g_{-}(c^{+},c^{-})\partial_{i}c^{-}
-\tilde{g}_{-}(c^{+},c^{-})\partial_{i}c^{+}-(u^{-}\cdot\nabla)u_{i}^{-}\nonumber\\
&\quad+\mu^{-}h_{-}(c^{+},c^{-})\partial_{j}c^{+}\partial_{j}u^{-}_{i}
+\mu^{-}k_{-}(c^{+},c^{-})\partial_{j}c^{-}\partial_{j}u^{-}_{i}\nonumber\\
&\quad+\mu^{-}h_{-}(c^{+},c^{-})\partial_{j}c^{+}\partial_{i}u^{-}_{j}
+\mu^{-}k_{-}(c^{+},c^{-})\partial_{j}c^{-}\partial_{i}u^{-}_{j}\\
&\quad+\lambda^{-}h_{-}(c^{+},c^{-})\partial_{i}c^{+}\partial_{j}u^{-}_{j}
+\lambda^{-}k_{-}(c^{+},c^{-})\partial_{i}c^{-}\partial_{j}u^{-}_{j}\nonumber\\
&\quad+\mu^{-}l_{-}(c^{+},c^{-})\partial_{j}^{2}u_{i}^{-}+(\mu^{-}+\lambda^{-})l_{-}
(c^{+},c^{-})\partial_{i}\partial_{j}u^{-}_{j}, \quad  i,j \in \{1,2,\cdots N\},\nonumber
\end{align}
where we define the nonlinear functions of $(c^{+},c^{-})$ by
\begin{align}\label{3.6}
\left\{
\begin{aligned}
&g_{+}(c^{+},c^{-})=\frac{(\mathcal{C}^{2}\rho^{-})(c^{+}+1,c^{-}+1)}
{\rho^{+}(c^{+}+1,c^{-}+1)}
-\frac{(\mathcal{C}^{2}\rho^{-})(1,1)}
{\rho^{+}(1,1)},
\\&g_{-}(c^{+},c^{-})=\frac{(\mathcal{C}^{2}\rho^{+})(c^{+}+1,c^{-}+1)}
{\rho^{-}(c^{+}+1,c^{-}+1)}
-\frac{(\mathcal{C}^{2}\rho^{+})(1,1)}
{\rho^{-}(1,1)},
\end{aligned}
\right.
\end{align}
\begin{align}\label{3.7}
\left\{
\begin{aligned}
&h_{+}(c^{+},c^{-})=\frac{(\mathcal{C}^{2}\alpha^{-})(c^{+}+1,c^{-}+1)}
{[c^{+}+1]s_{-}^{2}(c^{+}+1,c^{-}+1)},
\\&h_{-}(c^{+},c^{-})=-\frac{\mathcal{C}^{2}(c^{+}+1,c^{-}+1)}
{(\rho^{-}s_{-}^{2})(c^{+}+1,c^{-}+1)},
\end{aligned}
\right.
\end{align}
\begin{align}\label{3.8}
\left\{
\begin{aligned}
&k_{+}(c^{+},c^{-})=-{\frac{\mathcal{C}^{2}(c^{+}+1,c^{-}+1)}
{[c^{-}+1](s_{+}^{2}\rho)(c^{+}+1,c^{-}+1)}},
\\&k_{-}(c^{+},c^{-})={\frac{(\alpha^{+}\mathcal{C}^{2})(c^{+}+1,c^{-}+1)}
{[c^{-}+1]s_{+}^{2}(c^{+}+1,c^{-}+1)}},
\end{aligned}
\right.
\end{align}
\begin{align}\label{3.9}
\tilde{g}_{+}(c^{+},c^{-})=\tilde{g}_{-}(c^{+},c^{-})=\mathcal{C}^{2}(c^{+}+1,c^{-}+1)
-\mathcal{C}^{2}(1,1),
\end{align}
\begin{align}\label{3.10}
l_{\pm}(c^{+},c^{-})=\frac{1}{\rho^{\pm}(c^{+}+1,c^{-}+1)}-\frac{1}
{\rho^{\pm}(1,1)}.
\end{align}
\subsection{{A priori} estimates for  the linearized  system with convection terms}
\ \ \ \ \  Next, we investigate some \emph{a priori} estimates for the following linearized  system with convection terms
\begin{equation}\label{3.11}
\left\{
\begin{aligned}{}
&\p_tc^{+}+{v}^{+}\cdot\nabla c^{+}+\textrm{div}u^{+}=H_{11}(c^{+}, u^{+}),\\
&\p_t{u}^{+}+{v}^{+}\cdot\nabla u^{+}+\beta_{1}\nabla c^{+}
+\beta_{2}\nabla c^{-}-\nu_{1}^{+}\Delta u^{+}
-\nu_{2}^{+}\nabla\textrm{div}u^{+}=H_{21}(c^{+}, u^{+}, c^{-}),
\\&\p_tc^{-}+{v}^{-}\cdot\nabla c^{-}+\textrm{div}u^{-}=H_{31}(c^{-}, u^{-}),
\\&\p_tu^{-}+{v}^{-}\cdot\nabla u^{-}+\beta_{3}\nabla c^{+}
+\beta_{4}\nabla c^{-}-\nu_{1}^{-}\Delta u^{-}
-\nu_{2}^{-}\nabla\textrm{div}u^{-}=H_{41}(c^{+}, u^{-}, c^{-}),
\\&(c^{+},\,u^{+},\,c^{-},\,u^{-})|_{t=0}=(c^{+}_{0},\,u^{+}_{0},\,c^{-}_{0},\,u^{-}_{0}).
\end{aligned}
\right.
\end{equation}
We will establish a uniform estimate for a mixed hyperbolic-parabolic linear system \eqref{3.11}
with  convection terms. What is crucial in this work is to exploit the smoothing effects on the velocity $u^{+}, u^{+}$ and the $L^{1}$ decay on
$c^{+},c^{-}$,  which  play a key role to control the pressure term in the proof of the Theorem \ref{th:main1}.
\begin{Lemma}\label{prop3.1}
 Denote
\begin{equation*}
V(t):=\int_0^t\|({v}^{+},{v}^{-})(\tau)\|_{\dot{B}^{\f{N}{2}+1}_{2,1}}d\tau.
\label{VV}
\end{equation*}
Let $T>0$  and $(c^{+},\,u^{+},\,c^{-},\,u^{-})$ be a solution of the system
\eqref{3.11}. Then the following estimates
hold  for  $t\in[0,T)$
\begin{equation}\label{3.12}
\begin{split}
&\bigl\|(c^{+},\,c^{-})\bigr\|_{\widetilde{L}^{\infty}([0,t];\tilde{B}^{\f{N}{2}-1,\f{N}{2}}_{2,1})}+\bigl\|(u^{+},\,u^{-})\bigr\|_{\widetilde{L}^{\infty}([0,t];
\dot{B}^{\f{N}{2}-1}_{2,1})}
\\&\qquad+\int_0^t\big\|(c^{+},\,c^{-})(\tau)\big\|_{\tilde{B}^{\f{N}{2}+1,\f{N}{2}}_{2,1}}d\tau+\int_0^t\big\|(u^{+},\,u^{-})(\tau)\big\|_{
\dot{B}^{\f{N}{2}+1}_{2,1}}d\tau
\\&\quad \lesssim e^{CV(t)}\Big( \|(c^{+}_{0},\,c^{-}_{0})\|_{\tilde{B}^{\f{N}{2}-1,\f{N}{2}}_{2,1}}+\|(u^{+}_{0},\,u^{-}_{0})\|_{
\dot{B}^{\f{N}{2}-1}_{2,1}}+\int_0^t\big\|(H_{11},\,H_{31})(\tau)\big\|_{\tilde{B}^{\f{N}{2}-1,\f{N}{2}}_{2,1}}d\tau
\\&\qquad+\int_0^t\big\|(H_{21},\,H_{41})(\tau)\big\|_{\dot{B}^{\f{N}{2}-1}_{2,1}}d\tau\Big).
\end{split}
\end{equation}
\end{Lemma}
\noindent{\bf Proof.} Applying the operator $\dot{\Delta}_q$ to the
system \eqref{3.11}, we deduce that
$(\dot{\Delta}_qc^{+},\,\dot{\Delta}_qu^{+},\,\dot{\Delta}_qc^{-},\,\dot{\Delta}_qu^{-}
)$ satisfies
\begin{equation}\label{equ:heatflow-new31}
\left\{
\begin{aligned}{}
&\p_t\dot{\Delta}_qc^{+}+\dot{ \Delta}_q(v^{+}\cdot\na c^{+})+\dv \dot{ \Delta}_qu^{+}=\dot{ \Delta}_qH_{11},\\
&\p_t \dot{\Delta}_qu^{+}+\dot{ \Delta}_q(v^{+}\cdot\nabla u^{+})-\nu_{1}^{+}\Delta \dot{ \Delta}_qu^{+}-\nu_{2}^{+}\na
\dv \dot{\Delta}_qu^{+}+\beta_1\nabla\dot{ \Delta}_qc^{+}+\beta_2\nabla\dot{ \Delta}_qc^{-}=\dot{ \Delta}_qH_{21},\\
&\p_t\dot{\Delta}_qc^{-}+\dot{ \Delta}_q(v^{-}\cdot\na c^{-})+\dv \dot{ \Delta}_qu^{-}=\dot{ \Delta}_qH_{31},\\
&\p_t \dot{\Delta}_qu^{-}+\dot{ \Delta}_q(v^{-}\cdot\nabla u^{-})-\nu_{1}^{-}\Delta \dot{ \Delta}_qu^{-}-\nu_{2}^{-}\na
\dv \dot{\Delta}_qu^{-}+\beta_3\nabla\dot{ \Delta}_qc^{+}+\beta_4\nabla\dot{ \Delta}_qc^{-}=\dot{ \Delta}_qH_{41}.
\end{aligned}
\right.
\end{equation}
Taking the $L^2$-scalar product of the
first equation of \eqref{equ:heatflow-new31} with $\dot{\Delta}_qc^{+}$ and  $(-\Delta\dot{\Delta}_qc^{+})$,
the second equation with $\dot{\Delta}_qu^{+}$, the third equation with  $\dot{\Delta}_qc^{-}$ and  $(-\Delta\dot{\Delta}_qc^{-})$ and  the fourth equation with
$\dot{\Delta}_qu^{-}$ respectively,  we obtain the
following six identities:
\begin{equation}
\label{3.14}\begin{split}
\f{1}{2}\f{d}{dt}\|\dot{\Delta}_qc^{+}\|_{L^2}^2&+(\dot{\Delta}_q(v^{+}\cdot\nabla c^{+})|\dot{\Delta}_qc^{+})+(\dv\dot{\Delta}_qu^{+}|\dot{\Delta}_qc^{+})
\\&=(\dot{\Delta}_qH_{11}|\dot{\Delta}_qc^{+}),
\end{split}
\end{equation}
\begin{equation}
\label{3.15}\begin{split}
\f{1}{2}\f{d}{dt}\|\nabla\dot{\Delta}_qc^{+}\|_{L^2}^2&+(\dot{\Delta}_q(v^{+}\cdot\nabla c^{+})|-\Delta\dot{\Delta}_qc^{+})
+(\dv\dot{\Delta}_qu^{+}|-\Delta\dot{\Delta}_qc^{+})
\\&=(\dot{\Delta}_qH_{11}|-\Delta\dot{\Delta}_qc^{+}),
\end{split}
\end{equation}
\begin{equation}
\label{3.16}\begin{split}
\f{1}{2}\f{d}{dt}\|\Dot{\Delta}_qu^{+}\|^2_{L^2}&+(\Dot{\Delta}_q(v^{+}\cdot\nabla u^{+})|\Dot{\Delta}_qu^{+})+\nu_{1}^{+}\|\nabla\Dot{\Delta}_qu^{+}\|^2_{L^2}
\\&+\nu_{2}^{+}\|\dv\Dot{\Delta}_qu^{+}\|^2_{L^2}
+\beta_{1}(\nabla\Dot{\Delta}_qc^{+}|\Dot{\Delta}_qu^{+})+\beta_{2}(\nabla\Dot{\Delta}_qc^{-}|\Dot{\Delta}_qu^{+})
\\&=(\Dot{\Delta}_qH_{21}|\Dot{\Delta}_qu^{+}),
\end{split}
\end{equation}
\begin{equation}
\label{3.17}\begin{split}
\f{1}{2}\f{d}{dt}\|\dot{\Delta}_qc^{-}\|_{L^2}^2&+(\dot{\Delta}_q(v^{-}\cdot\nabla c^{-})|\dot{\Delta}_qc^{-})+(\dv\dot{\Delta}_qu^{-}|\dot{\Delta}_qc^{-})
\\&=(\dot{\Delta}_qH_{31}|\dot{\Delta}_qc^{-}),
\end{split}
\end{equation}
\begin{equation}
\label{3.18}\begin{split}
\f{1}{2}\f{d}{dt}\|\nabla\dot{\Delta}_qc^{-}\|_{L^2}^2&+(\dot{\Delta}_q(v^{-}\cdot\nabla c^{-})|-\Delta\dot{\Delta}_qc^{-})
+(\dv\dot{\Delta}_qu^{-}|-\Delta\dot{\Delta}_qc^{-})
\\&=(\dot{\Delta}_qH_{31}|-\Delta\dot{\Delta}_qc^{-}),
\end{split}
\end{equation}
\begin{equation}
\label{3.19}\begin{split}
\f{1}{2}\f{d}{dt}\|\Dot{\Delta}_qu^{-}\|^2_{L^2}&+(\Dot{\Delta}_q(v^{-}\cdot\nabla u^{-})|\Dot{\Delta}_qu^{-})+\nu_{1}^{-}\|\nabla\Dot{\Delta}_qu^{-}\|^2_{L^2}
\\&+\nu_{2}^{-}\|\dv\Dot{\Delta}_qu^{-}\|^2_{L^2}
+\beta_{3}(\nabla\Dot{\Delta}_qc^{+}|\Dot{\Delta}_qu^{-})+\beta_{4}(\nabla\Dot{\Delta}_qc^{-}|\Dot{\Delta}_qu^{-})
\\&=(\Dot{\Delta}_qH_{41}|\Dot{\Delta}_qu^{-}),
\end{split}
\end{equation}
where $(\cdot\,|\,\cdot)$ stands for the $L^{2}$ inner product.

In order to obtain a second energy estimate, we need to derive some  identities involving
$(\Dot{\Delta}_qc^{+}|\Dot{\Delta}_qc^{-})$, $(\nabla\Dot{\Delta}_qc^{+}|\Dot{\Delta}_qu^{+})$, $(\nabla\Dot{\Delta}_qc^{-}|\Dot{\Delta}_qu^{-})$.
 Taking the $L^2$-scalar product of the
first equation of \eqref{equ:heatflow-new31} with $\dot{\Delta}_qc^{-}$ and
the third equation with $\dot{\Delta}_qc^{+}$ and then summing the results, which yields
\begin{equation} \label{3.20}\begin{split}
&\f{d}{dt}(\Dot{\Delta}_qc^{+}|\Dot{\Delta}_qc^{-})+(\dot{\Delta}_q(v^{+}\cdot\nabla c^{+})|\dot{\Delta}_qc^{-})+(\dot{\Delta}_q(v^{-}\cdot\nabla c^{-})|\dot{\Delta}_qc^{+})
\\&\qquad+(\dv\dot{\Delta}_qu^{+}|\dot{\Delta}_qc^{-})+(\dv\dot{\Delta}_qu^{-}|\dot{\Delta}_qc^{+})\\
&\quad=(\dot{\Delta}_qH_{11}|\dot{\Delta}_qc^{-})+(\dot{\Delta}_qH_{31}|\dot{\Delta}_qc^{+}).
\end{split}
\end{equation}
 On the other hand,
applying the operator $\nabla$ to the first equation in
\eqref{equ:heatflow-new3} and taking the $L^2$ scalar product with
$\Dot{\Delta}_qu^{+}$, then calculating the scalar product of the second
equation in \eqref{equ:heatflow-new3}  with $\nabla\Dot{\Delta}_qc^{+}$,
and then  summing up the results, we get
\begin{equation} \label{3.21}\begin{split}
&\f{d}{dt}(\nabla\Dot{\Delta}_qc^{+}|\Dot{\Delta}_qu^{+})-\|\dv\Dot{\Delta}_qu^{+}\|^2_{L^2}+(\nu_{1}^{+}+\nu_{2}^{+})(\dv\Dot{\Delta}_qu^{+}|\Delta\Dot{\Delta}_qc^{+})
+\beta_1\|\nabla\Dot{\Delta}_qc^{+}\|^2_{L^2}\\
&\qquad+\beta_2(\nabla\Dot{\Delta}_q c^{-}|\nabla\Dot{\Delta}_q c^{+})+(\nabla\Dot{\Delta}_q(v^{+}\cdot\nabla c^{+})|\Dot{\Delta}_qu^{+})+(\Dot{\Delta}_q(v^{+}\cdot\nabla u^{+})|\nabla\Dot{\Delta}_qc^{+})\\
&\quad=(\nabla\Dot{\Delta}_qH_{11}|\Dot{\Delta}_qu^{+})+(\Dot{\Delta}_qH_{21}|\nabla\Dot{\Delta}_qc^{+}).
\end{split}
\end{equation}
Similarly,
\begin{equation} \label{3.22}\begin{split}
&\f{d}{dt}(\nabla\Dot{\Delta}_qc^{-}|\Dot{\Delta}_qu^{-})-\|\dv\Dot{\Delta}_qu^{-}\|^2_{L^2}+(\nu_{1}^{-}+\nu_{2}^{-})(\dv\Dot{\Delta}_qu^{-}|\Delta\Dot{\Delta}_qc^{-})
+\beta_3(\nabla\Dot{\Delta}_q c^{+}|\nabla\Dot{\Delta}_q c^{-})\\
&\qquad+\beta_4\|\nabla\Dot{\Delta}_qc^{-}\|^2_{L^2}+(\nabla\Dot{\Delta}_q(v^{-}\cdot\nabla c^{-})|\Dot{\Delta}_qu^{-})+(\Dot{\Delta}_q(v^{-}\cdot\nabla u^{-})|\nabla\Dot{\Delta}_qc^{-})\\
&\quad=(\nabla\Dot{\Delta}_qH_{31}|\Dot{\Delta}_qu^{-})+(\Dot{\Delta}_qH_{41}|\nabla\Dot{\Delta}_qc^{-}).
\end{split}
\end{equation}
Define
\begin{equation*}
\begin{split}
\alpha_q^2&=\beta_1\|\Dot{\Delta}_qc^{+}\|^2_{L^2}+\beta_4\|\Dot{\Delta}_qc^{-}\|^2_{L^2}+A(\nu^{+}_{1}+\nu^{+}_{2})\|\nabla\Dot{\Delta}_qc^{+}\|^2_{L^2}
+A(\nu^{-}_{1}+\nu^{-}_{2})\|\nabla\Dot{\Delta}_qc^{-}\|^2_{L^2}
\\&\quad+\|\Dot{\Delta}_qu^{+}\|^2_{L^2}+\|\Dot{\Delta}_qu^{-}\|^2_{L^2}+2\beta_2(\Dot{\Delta}_qc^{+}|\Dot{\Delta}_qc^{-})
+2A(\nabla\Dot{\Delta}_qc^{+}|\Dot{\Delta}_qu^{+})+2A(\nabla\Dot{\Delta}_qc^{-}|\Dot{\Delta}_qu^{-}),
\end{split}
\end{equation*}
where $A=\frac{1}{4}\min\{\nu^{+}_{2},\,\nu^{-}_{2}\}>0$. Taking $M_1\in(\frac{1}{4A}, \frac{1}{A})$ and employing Young's inequality, we obtain
\begin{equation*}
\begin{split}
\big|2A(\nabla\Dot{\Delta}_qc^{+}|\Dot{\Delta}_qu^{+})\big|&\leq M_1A\|\Dot{\Delta}_qu^{+}\|^2_{L^2}+\frac{A}{M_1}\|\nabla\Dot{\Delta}_qc^{+}\|^2_{L^2}\\
&\leq \|\Dot{\Delta}_qu^{+}\|^2_{L^2}+A(\nu^{+}_{1}+\nu^{+}_{2})\|\nabla\Dot{\Delta}_qc^{+}\|^2_{L^2},
\end{split}
\end{equation*}
\begin{equation*}
\begin{split}
\big|2A(\nabla\Dot{\Delta}_qc^{-}|\Dot{\Delta}_qu^{-})\big|&\leq M_1A\|\Dot{\Delta}_qu^{-}\|^2_{L^2}+\frac{A}{M_1}\|\nabla\Dot{\Delta}_qc^{-}\|^2_{L^2}\\
&\leq \|\Dot{\Delta}_qu^{-}\|^2_{L^2}+A(\nu^{-}_{1}+\nu^{-}_{2})\|\nabla\Dot{\Delta}_qc^{-}\|^2_{L^2}.
\end{split}
\end{equation*}
Using further $\beta^{2}_2=\beta^{2}_3=\beta_{1}\beta_{4}$ and Young's inequality,  and choosing $M_2=\frac{\beta_1}{\beta_2}$,  we get
\begin{equation*}
\begin{split}
2\beta_2\big|(\Dot{\Delta}_qc^{+}|\Dot{\Delta}_qc^{-})\big|&\leq M_2\beta_2\|\Dot{\Delta}_qc^{+}\|^2_{L^2}+\frac{\beta_2}{M_2}\|\Dot{\Delta}_qc^{-}\|^2_{L^2}\\
&\leq \beta_1\|\Dot{\Delta}_qc^{+}\|^2_{L^2}+\beta_4\|\Dot{\Delta}_qc^{-}\|^2_{L^2},
\end{split}
\end{equation*}
\begin{equation*}
\begin{split}
2\beta_2\big|(\nabla\Dot{\Delta}_q c^{+}|\nabla\Dot{\Delta}_q c^{-})\big|&\leq M_2\beta_2\|\nabla\Dot{\Delta}_qc^{+}\|^2_{L^2}+\frac{\beta_2}{M_2}\|\nabla\Dot{\Delta}_qc^{-}\|^2_{L^2}\\
&\leq \beta_1\|\nabla\Dot{\Delta}_qc^{+}\|^2_{L^2}+\beta_4\|\nabla\Dot{\Delta}_qc^{-}\|^2_{L^2}.
\end{split}
\end{equation*}
Then, there exist two positive constants $c_1$ and $c_2$ such that
\begin{equation*}
\begin{split}
c_1\alpha_q^2\leq\|\Dot{\Delta}_qc^{+}\|^2_{L^2}+\|\Dot{\Delta}_qc^{-}\|^2_{L^2}+\|\nabla\Dot{\Delta}_qc^{+}\|^2_{L^2}
+\|\nabla\Dot{\Delta}_qc^{-}\|^2_{L^2}
+\|\Dot{\Delta}_qu^{+}\|^2_{L^2}+\|\Dot{\Delta}_qu^{-}\|^2_{L^2}
\leq c_2\alpha_q^2.
\end{split}
\end{equation*}
Thus, for some fix $q_0$,
\begin{equation*}\label{equ:heatflow-new3}
\alpha_q\approx\left\{
\begin{aligned}{}
&\|(\Dot{\Delta}_qc^{+},\,\Dot{\Delta}_qu^{+},\,\Dot{\Delta}_qc^{-},\,\Dot{\Delta}_qu^{-})\|_{L^2},  \hbox{~~~~~~~for~}  q\leq q_0,\\
&\|(\nabla\Dot{\Delta}_qc^{+},\,\Dot{\Delta}_qu^{+},\,\nabla\Dot{\Delta}_qc^{-},\,\Dot{\Delta}_qu^{-})\|_{L^2},   \hbox{~for~} q> q_0.
\end{aligned}
\right.
\end{equation*}
Combining with \eqref{3.14}-\eqref{3.22}, it yields, with the help of Proposition \ref{p29}, that
\begin{equation}\label{17}
\begin{split}
&\f{1}{2}\f{d}{dt}\alpha^2_q+(\nu^{+}_{2}-A)\|\dv\Dot{\Delta}_qu^{+}\|^2_{L^2}+(\nu^{-}_{2}-A)\|\dv\Dot{\Delta}_qu^{-}\|^2_{L^2}+\nu^{+}_{1}\|\nabla\Dot{\Delta}_qu^{+}\|^2_{L^2}
+\nu^{+}_{2}\|\nabla\Dot{\Delta}_qu^{-}\|^2_{L^2}\\&\qquad+A\beta_1\|\nabla\Dot{\Delta}_q c^{+}\|^2_{L^2}+A\beta_4\|\nabla\Dot{\Delta}_q c^{-}\|^2_{L^2}+2A\beta_2(\nabla\Dot{\Delta}_q c^{+}|\nabla\Dot{\Delta}_q c^{-})
\\&\quad=-\beta_1\big(\Dot{\Delta}_q(v^{+}\cdot\nabla c^{+})|\Dot{\Delta}_qc^{+}\big)+2A(\nu^{+}_1+\nu^{+}_2)\big(\Dot{\Delta}_q(v^{+}\cdot\nabla c^{+})|\Delta\Dot{\Delta}_qc^{+}\big)-\big(\Dot{\Delta}_q(v^{+}_{1}\cdot\nabla u^{+})|\Dot{\Delta}_qu^{+}\big)
\\&\qquad-\beta_4\big(\Dot{\Delta}_q(v^{-}\cdot\nabla c^{-})|\Dot{\Delta}_qc^{+}\big)+2A(\nu^{-}_1+\nu^{-}_2)\big(\Dot{\Delta}_q(v^{-}\cdot\nabla c^{-})|\Delta\Dot{\Delta}_qc^{+}\big)-\big(\Dot{\Delta}_q(v^{-}\cdot\nabla u^{-})|\Dot{\Delta}_qu^{-}\big)
\\&\qquad-\beta_2\big(\Dot{\Delta}_q(v^{+}\cdot\nabla c^{+})|\Dot{\Delta}_qc^{-}\big)-\beta_3\big(\Dot{\Delta}_q(v^{-}\cdot\nabla c^{-})|\Dot{\Delta}_qc^{+}\big)
-A\big(\nabla\Dot{\Delta}_q(v^{+}\cdot\nabla c^{+})|\Dot{\Delta}_qu^{+}\big)
\\&\qquad-A\big(\Dot{\Delta}_q(v^{+}\cdot\nabla u^{+})|\nabla\Dot{\Delta}_qc^{+}\big)
-A\big(\nabla\Dot{\Delta}_q(v^{-}\cdot\nabla c^{-})|\Dot{\Delta}_qu^{-}\big)-A\big(\Dot{\Delta}_q(v^{-}\cdot\nabla u^{-})|\nabla\Dot{\Delta}_qc^{-}\big)
\\&\qquad+\big(\Dot{\Delta}_qH_{11}|\Dot{\Delta}_qc^{+}\big)-2A(\nu^{+}_1+\nu^{+}_2)\big(\Dot{\Delta}_qH_{11}|-\Delta\Dot{\Delta}_qc^{+}\big)
+\big(\Dot{\Delta}_qH_{21}|\Dot{\Delta}_qu^{+}\big)+\big(\Dot{\Delta}_qH_{31}|\Dot{\Delta}_qc^{-}\big)\\&\qquad-2A(\nu^{-}_1+\nu^{-}_2)\big(\Dot{\Delta}_qH_{21}|-\Delta\Dot{\Delta}_qc^{-}\big)
+\big(\Dot{\Delta}_qH_{31}|\Dot{\Delta}_qu^{-}\big)+\beta_2\big(\Dot{\Delta}_qH_{11}|\Dot{\Delta}_qc^{-}\big)+\beta_3\big(\Dot{\Delta}_qH_{31}|\Dot{\Delta}_qc^{+}\big)
\\&\qquad+A\big(\nabla\Dot{\Delta}_qH_{11}|\Dot{\Delta}_qu^{+}\big)+A\big(\Dot{\Delta}_qH_{21}|\nabla\Dot{\Delta}_qc^{+}\big)
+A\big(\nabla\Dot{\Delta}_qH_{31}|\Dot{\Delta}_qu^{-}\big)+A\big(\Dot{\Delta}_qS_{41}|\nabla\Dot{\Delta}_qc^{-}\big)
\\&\quad\lesssim \alpha_q\Big(\|\Dot{\Delta}_qH_{11}\|_{L^2}+\|\nabla\Dot{\Delta}_qH_{11}\|_{L^2}+\|\Dot{\Delta}_qH_{21}\|_{L^2}+\|\Dot{\Delta}_qH_{31}\|_{L^2}
+\|\nabla\Dot{\Delta}_qH_{31}\|_{L^2}+\|\Dot{\Delta}_qH_{41}\|_{L^2}
\\&\qquad+2^{-(\frac{N}{2}-1)q}\gamma_qV'\|(c^{+},\,u^{+},\,c^{-},\,u^{-})\|_{\tilde{B}^{\f{N}{2}-1,\f{N}{2}}_{2,1}\times
\dot{B}^{\f{N}{2}-1}_{2,1}\times \tilde{B}^{\f{N}{2}-1,\f{N}{2}}_{2,1}\times
\dot{B}^{\f{N}{2}-1}_{2,1}}\Big)
\end{split}\end{equation}
with $(\gamma_q)_{q\in \mathbb{Z}}$ in  the unit sphere of $\ell^{1}(\mathbb{Z})$.

Thus, it follows
\begin{equation}\label{171}
\begin{split}
  &\frac{1}{2}\frac{d}{dt}\alpha_q^2+c_0\min(2^{2q},1)\alpha_q^2\\
  &\quad\lesssim \gamma_q2^{-(\frac{N}{2}-1)q}\Big[\|H_{11}\|_{\tilde{B}^{\f{N}{2}-1,\f{N}{2}}_{2,1}} +\|H_{21}\|_{\dot{B}^{\f{N}{2}-1}_{2,1}}+\|H_{31}\|_{\tilde{B}^{\f{N}{2}-1,\f{N}{2}}_{2,1}}+\|H_{41}\|_{\dot{B}^{\f{N}{2}-1,\f{N}{2}}_{2,1}}
  \\&\qquad+V'\|(c^{+},\,u^{+},\,c^{-},\,u^{-})\|_{\tilde{B}^{\f{N}{2}-1,\f{N}{2}}_{2,1}\times
\dot{B}^{\f{N}{2}-1}_{2,1}\times \tilde{B}^{\f{N}{2}-1,\f{N}{2}}_{2,1}\times
\dot{B}^{\f{N}{2}-1}_{2,1}}\Big)\Big]\alpha_q,
\end{split}\end{equation}
which implies that
\begin{align*}
  &2^{(\frac{N}{2}-1)q}\alpha_q+c_0\int_0^t\min(2^{2q},1)2^{(\frac{N}{2}-1)q}\alpha_q(\tau)d\tau\\
  &\quad\lesssim2^{(\frac{N}{2}-1)q}\alpha_q(0)+C\gamma_q\int_0^t\Big[\|H_{11}\|_{\tilde{B}^{\f{N}{2}-1,\f{N}{2}}_{2,1}} +\|H_{21}\|_{\dot{B}^{\f{N}{2}-1}_{2,1}}+\|H_{31}\|_{\tilde{B}^{\f{N}{2}-1,\f{N}{2}}_{2,1}}+\|H_{41}\|_{\dot{B}^{\f{N}{2}-1,\f{N}{2}}_{2,1}}
  \\&\qquad+V'\sum_q2^{(\frac{N}{2}-1)q}\alpha_q\Big].
\end{align*}
Thus, by  Gronwall's inequality, we have
\begin{equation}\label{3.24}
\begin{split}
&\bigl\|(c^{+},\,c^{-})\bigr\|_{\widetilde{L}^{\infty}([0,t];\tilde{B}^{\f{N}{2}-1,\f{N}{2}}_{2,1})}+\bigl\|(u^{+},\,u^{-})\bigr\|_{\widetilde{L}^{\infty}([0,t];
\dot{B}^{\f{N}{2}-1}_{2,1})}
\\&\qquad+\int_0^t\big\|(c^{+},\,c^{-})(\tau)\big\|_{\tilde{B}^{\f{N}{2}+1,\f{N}{2}}_{2,1}}d\tau
+\int_0^t\big\|(u^{+},\,u^{-})(\tau)\big\|_{\tilde{B}^{\f{N}{2}+1,\f{N}{2}-1}_{2,1}}d\tau
\\&\quad \lesssim e^{CV(t)}\Big( \|(c^{+}_{0},\,c^{-}_{0})\|_{\tilde{B}^{\f{N}{2}-1,\f{N}{2}}_{2,1}}+\|(u^{+}_{0},\,u^{-}_{0})\|_{
\dot{B}^{\f{N}{2}-1}_{2,1}}+\int_0^t\big\|(H_{11},\,H_{31})(\tau)\big\|_{\tilde{B}^{\f{N}{2}-1,\f{N}{2}}_{2,1}}d\tau
\\&\qquad+\int_0^t\big\|(H_{21},\,H_{41})(\tau)\big\|_{\dot{B}^{\f{N}{2}-1}_{2,1}}d\tau\Big).
\end{split}
\end{equation}
Based on the damping effects for $c^{+}$ and $c^{-}$,  we further  exploit  the smoothing
effects of $u^{+}$ and $u^{-}$ in high frequencies regime    by considering
\eqref{3.11} with $\nabla c^{+}$ and $\nabla c^{-}$ being viewed as  source
terms.  From \eqref{3.16} and \eqref{3.19}, we have
\begin{equation*}
\label{3.25}\begin{split}
&\f{1}{2}\f{d}{dt}\|\Dot{\Delta}_qu^{+}\|^2_{L^2}+C2^{2q}\|\Dot{\Delta}_qu^{+}\|^2_{L^2}\\
&\lesssim-(\Dot{\Delta}_q(v^{+}\cdot\nabla u^{+})|\Dot{\Delta}_qu^{+})
-\beta_{1}(\nabla\Dot{\Delta}_qc^{+}|\Dot{\Delta}_qu^{+})-\beta_{2}(\nabla\Dot{\Delta}_qc^{-}|\Dot{\Delta}_qu^{+})
\\&\quad+(\Dot{\Delta}_qH_{21}|\Dot{\Delta}_qu^{+}),
\end{split}
\end{equation*}
\begin{equation*}
\label{3.26}\begin{split}
&\f{1}{2}\f{d}{dt}\|\Dot{\Delta}_qu^{-}\|^2_{L^2}+C2^{2q}\|\Dot{\Delta}_qu^{+}\|^2_{L^2}\\
&\lesssim-(\Dot{\Delta}_q(v^{-}\cdot\nabla u^{-})|\Dot{\Delta}_qu^{-})
-\beta_{3}(\nabla\Dot{\Delta}_qc^{+}|\Dot{\Delta}_qu^{-})-\beta_{4}(\nabla\Dot{\Delta}_qc^{-}|\Dot{\Delta}_qu^{-})
\\&\quad+(\Dot{\Delta}_qH_{41}|\Dot{\Delta}_qu^{-}).
\end{split}
\end{equation*}
It follows that from Proposition \ref{p29}
\begin{equation*}
\label{3.25}\begin{split}
&\f{d}{dt}\sum_{q\geq q_0}2^{(\frac{N}{2}-1)q}\|\Dot{\Delta}_qu^{+}\|_{L^2}+C\sum_{q\geq q_0}2^{(\frac{N}{2}-1)q}2^{2q}\|\Dot{\Delta}_qu^{+}\|_{L^2}
\\&\quad\lesssim \sum_{q\geq q_0}2^{(\frac{N}{2}-1)q}\Big(2^{q}\|\Dot{\Delta}_qc^{+}\|_{L^{2}}+2^{q}\|\Dot{\Delta}_qc^{-}\|_{L^{2}}+\|\Dot{\Delta}_qH_{21}\|_{L^{2}}+ V'(t)\gamma_{q}2^{-(\frac{N}{2}-1)q}\|u^{+}\|_{
\dot{B}^{\frac{N}{2}-1}_{2,1}}\Big)
\\&\quad\lesssim \sum_{q\geq q_0}2^{(\frac{N}{2}-1)q}2^{q}\big(\|\Dot{\Delta}_qc^{+}\|_{L^{2}}+\|\Dot{\Delta}_qc^{-}\|_{L^{2}}\big)+\|H_{21}\|_{\dot{B}^{\frac{N}{2}-1}_{2,1}}+ V'(t)\|u^{+}\|_{
\dot{B}^{\frac{N}{2}-1}_{2,1}},
\end{split}
\end{equation*}
\begin{equation*}
\label{3.26}\begin{split}
&\f{d}{dt}\sum_{q\geq q_0}2^{(\frac{N}{2}-1)q}\|\Dot{\Delta}_qu^{-}\|_{L^2}+C\sum_{q\geq q_0}2^{(\frac{N}{2}-1)q}2^{2q}\|\Dot{\Delta}_qu^{-}\|_{L^2}
\\&\quad\lesssim \sum_{q\geq q_0}2^{(\frac{N}{2}-1)q}\Big(2^{q}\|\Dot{\Delta}_qc^{+}\|_{L^{2}}+2^{q}\|\Dot{\Delta}_qc^{-}\|_{L^{2}}+\|\Dot{\Delta}_qH_{41}\|_{L^{2}}+ V'(t)\gamma_{q}2^{-(\frac{N}{2}-1)q}\|u^{-}\|_{
\dot{B}^{\frac{N}{2}-1}_{2,1}}\Big)
\\&\quad\lesssim \sum_{q\geq q_0}2^{(\frac{N}{2}-1)q}2^{q}\big(\|\Dot{\Delta}_qc^{+}\|_{L^{2}}+\|\Dot{\Delta}_qc^{-}\|_{L^{2}}\big)+\|H_{41}\|_{\dot{B}^{\frac{N}{2}-1}_{2,1}}+ V'(t)\|u^{-}\|_{
\dot{B}^{\frac{N}{2}-1}_{2,1}},
\end{split}
\end{equation*}
which implies, with the help of \eqref{3.24}, that
\begin{equation}
\label{3.25}\begin{split}
&\int_0^tC\sum_{q\geq q_0}2^{(\frac{N}{2}-1)q}2^{2q}\Big(\|\Dot{\Delta}_qu^{+}\|_{L^2}+\|\Dot{\Delta}_qu^{-}\|_{L^2}\Big)d\tau
\\&\quad \lesssim e^{CV(t)}\Big( \|(c^{+}_{0},\,c^{-}_{0})\|_{\tilde{B}^{\f{N}{2}-1,\f{N}{2}}_{2,1}}+\|(u^{+}_{0},\,u^{-}_{0})\|_{
\dot{B}^{\f{N}{2}-1}_{2,1}}+\int_0^t\big\|(H_{11},\,H_{31})(\tau)\big\|_{\tilde{B}^{\f{N}{2}-1,\f{N}{2}}_{2,1}}d\tau
\\&\qquad+\int_0^t\big\|(H_{21},\,H_{41})(\tau)\big\|_{\dot{B}^{\f{N}{2}-1}_{2,1}}d\tau\Big).
\end{split}
\end{equation}
Combining with \eqref{3.24} and \eqref{3.25}, we finally conclude that \eqref{3.12}. Thus,  we complete the
proof of Lemma \ref{prop3.1}.
\section{Global existence for initial data near equilibrium}
\ \ \ \ \   In this section, we  show that if the initial
 data satisfy
\begin{equation*}
\|(R^{+}_{0}-1,\,R^{+}_{0}-1)\|_{\tilde{B}^{\f{N}{2}-1,\f{N}{2}}_{2,1}}+\|(u^{+}_{0},\,u^{-}_{0})\|_{
\dot{B}^{\f{N}{2}-1}_{2,1}}\leq \eta,  \label{1Hn/2}
\end{equation*} for some sufficiently small $\eta$,  then there exists a
positive constant $M$ such that
\begin{equation*}X(t)\le M\eta.\end{equation*}
This uniform estimate will enable us to
 extend the local solution $(R^{+}-1,\,u^{+},\,R^{-}-1,\,u^{-})$ obtained within an iterative scheme as in \cite{Dan1} to  a global one. To this
end, we use a contradiction argument. Define
$$T_0=\sup\big\{T\in[0,\infty): X(T)\le
M\eta\big\},$$ with $M$ to be determined later. Suppose
that $T_0<\infty$. We apply the  linear estimates in Lemma \ref{prop3.1}
to the solutions of the reformulated system
\eqref{3.1} such that for all $t\in [0,T_0]$, the
following estimates hold
\begin{equation}\label{4.1}
\begin{split}
X(T_0)& \lesssim e^{CV(t)}\Big( \|(c^{+}_{0},\,c^{-}_{0})\|_{\tilde{B}^{\f{N}{2}-1,\f{N}{2}}_{2,1}}+\|(u^{+}_{0},\,u^{-}_{0})\|_{
\dot{B}^{\f{N}{2}-1}_{2,1}}+\int_0^{T_0}\big\|(H_{11},\,H_{31})(\tau)\big\|_{\tilde{B}^{\f{N}{2}-1,\f{N}{2}}_{2,1}}d\tau
\\&\qquad+\int_0^{T_0}\big\|(H_{21},\,H_{41})(\tau)\big\|_{\dot{B}^{\f{N}{2}-1}_{2,1}}d\tau\Big),
\end{split}
\end{equation}
 where \begin{equation*}
V(T_0)=\int_0^{T_0}\|(u^{+},\,u^{-})(\tau)\|_{\dot{B}^{\f{N}{2}+1}_{2,1}}d\tau.
\label{VV}
\end{equation*}
In what follows, we derive some estimates for the nonlinear terms
$H_{11}-H_{41}$.  First, by Proposition \ref{p25},  we have
\begin{equation}\label{4.2}
\begin{split}\|(H_{11},\,H_{31})\|_{L^1_{T_0}(\tilde{B}^{\f{N}{2}-1,\f{N}{2}}_{2,1})}&\lesssim\|c^{+}\|_{L^\infty_{T_0}(\tilde{B}^{\f{N}{2}-1,\f{N}{2}}_{2,1})}\|\Dv u^{+}\|_{L^1_{T_0}(\tilde{B}^{\f{N}{2},\f{N}{2}}_{2,1})}
\\\qquad&+\|c^{-}\|_{L^\infty_{T_0}(\tilde{B}^{\f{N}{2}-1,\f{N}{2}}_{2,1})}\|\Dv u^{-}\|_{L^1_{T_0}(\tilde{B}^{\f{N}{2},\f{N}{2}}_{2,1})}\\
&\lesssim M^2\eta^2.
\end{split}
\end{equation}
Next, we bound the terms $H_{21}$ and $H_{41}$.
By the embedding $\tilde{B}^{\frac{N}{2}-1,\f
N2}_{2,1}\hookrightarrow \dot{B}^\frac N2_{2,1}\hookrightarrow
L^\infty$ and Proposition \ref{p27} (ii), we get
\begin{equation*}\label{318}
\begin{split}&\big\|l_{+}(c^{+},c^{-})\big\|_{\dot{B}^\frac N2_{2,1}}\\
\\&\leq C_0\Big(\|c^{+}\|_{L^\infty}, \|c^{-}\|_{L^\infty}\Big)\Big(\|c^{+}\|_{\dot{B}^\frac N 2_{2,1}}+\|c^{-}\|_{\dot{B}^\frac N 2_{2,1}}\Big)
\\&\lesssim \|c^{+}\|_{\tilde{B}^{\frac{N}{2}-1,\f N2}_{2,1}}+\|c^{-}\|_{\tilde{B}^{\frac{N}{2}-1,\f N2}_{2,1}},
\end{split}
\end{equation*}
\begin{equation*}\label{318}
\begin{split}&\Big\|h_{+}(c^{+},c^{-})-\frac{(\mathcal{C}^{2}\alpha^{-})(1,1)}
{s_{-}^{2}(1,1)}\Big\|_{\dot{B}^\frac N2_{2,1}}\\
\\&\leq C_0\Big(\|c^{+}\|_{L^\infty}, \|c^{-}\|_{L^\infty}\Big)\Big(\|c^{+}\|_{\dot{B}^\frac N 2_{2,1}}+\|c^{-}\|_{\dot{B}^\frac N 2_{2,1}}\Big)
\\&\lesssim \|c^{+}\|_{\tilde{B}^{\frac{N}{2}-1,\f N2}_{2,1}}+\|c^{-}\|_{\tilde{B}^{\frac{N}{2}-1,\f N2}_{2,1}},
\end{split}
\end{equation*}
thus, thanks to Proposition \ref{p26}, we easily infer
\begin{equation*}\label{318}
\begin{split}&\|l_{+}(c^{+},c^{-})\partial_{j}^{2}u_{i}^{+}\|_{L^1_{T_0}(\dot{B}^{\f{N}{2}-1}_{2,1})}\\
\\&\lesssim \Big(\|c^{+}\|_{L^\infty_{T_0}(\tilde{B}^{\f{N}{2}-1,\f{N}{2}}_{2,1})}+\|c^{-}\|_{L^\infty_{T_0}(\tilde{B}^{\f{N}{2}-1,\f{N}{2}}_{2,1})}\Big)
\|\partial_{j}^{2}u_{i}^{+}\|_{L^1_{T_0}(\dot{B}^{\f{N}{2}-1}_{2,1})}
\\&\lesssim \Big(\|c^{+}\|_{L^\infty_{T_0}(\tilde{B}^{\f{N}{2}-1,\f{N}{2}}_{2,1})}+\|c^{-}\|_{L^\infty_{T_0}(\tilde{B}^{\f{N}{2}-1,\f{N}{2}}_{2,1})}\Big)
\|u^{+}_{i}\|_{L^1_{T_0}(\dot{B}^{\f{N}{2}+1}_{2,1})}
\\&\lesssim M^2\eta^2,
\end{split}
\end{equation*}
\begin{equation*}\label{318}
\begin{split}&\|h_{+}(c^{+},c^{-})\partial_{j}c^{+}\partial_{j}u^{+}_{i}\|_{L^1_{T_0}(\dot{B}^{\f{N}{2}-1}_{2,1})}\\
\\&\lesssim \Big(1+\|c^{+}\|_{L^\infty_{T_0}(\tilde{B}^{\f{N}{2}-1,\f{N}{2}}_{2,1})}+\|c^{-}\|_{L^\infty_{T_0}(\tilde{B}^{\f{N}{2}-1,\f{N}{2}}_{2,1})}\Big)
\big\|\partial_{j}c^{+}\partial_{j}u^{+}_{i}\big\|_{L^1_{T_0}(\dot{B}^{\f{N}{2}-1}_{2,1})}
\\&\lesssim \Big(1+\|c^{+}\|_{L^\infty_{T_0}(\tilde{B}^{\f{N}{2}-1,\f{N}{2}}_{2,1})}+\|c^{-}\|_{L^\infty_{T_0}(\tilde{B}^{\f{N}{2}-1,\f{N}{2}}_{2,1})}\Big)
\|\partial_{j}c^{+}\|_{L^\infty_{T_0}(\dot{B}^{\f{N}{2}-1}_{2,1})}\|\partial_{j}u^{+}_{i}\|_{L^1_{T_0}(\dot{B}^{\f{N}{2}}_{2,1})}
\\&\lesssim \Big(1+\|c^{+}\|_{L^\infty_{T_0}(\tilde{B}^{\f{N}{2}-1,\f{N}{2}}_{2,1})}+\|c^{-}\|_{L^\infty_{T_0}(\tilde{B}^{\f{N}{2}-1,\f{N}{2}}_{2,1})}\Big)
\|c^{+}\|_{L^\infty_{T_0}(\dot{B}^{\f{N}{2}}_{2,1})}\|u^{+}_{i}\|_{L^1_{T_0}(\dot{B}^{\f{N}{2}+1}_{2,1})}
\\&\lesssim M^2\eta^2+M^3\eta^3.
\end{split}
\end{equation*}
Similarly,
\begin{equation*}\label{318}
\begin{split}\Big\|k_{+}(c^{+},c^{-})\partial_{j}c^{-}\partial_{j}u^{+}_{i}\Big\|_{L^1_{T_0}(\dot{B}^{\f{N}{2}-1}_{2,1})}
\lesssim M^2\eta^2+M^3\eta^3.
\end{split}
\end{equation*}
\begin{equation*}\label{318}
\begin{split}\Big\|u^{+}\cdot\nabla u_{i}^{+}\Big\|_{L^1_{T_0}(\dot{B}^{\f{N}{2}-1}_{2,1})}
&\lesssim \|u^{+}\|_{L^\infty_{T_0}(\dot{B}^{\f{N}{2}-1}_{2,1})}\|\nabla u^{+}\|_{L^1_{T_0}(\dot{B}^{\f{N}{2}}_{2,1})}
\\&\lesssim \|u^{+}\|_{L^\infty_{T_0}(\dot{B}^{\f{N}{2}-1}_{2,1})}\| u^{+}\|_{L^1_{T_0}(\dot{B}^{\f{N}{2}+1}_{2,1})}
\\&\lesssim M^2\eta^2.
\end{split}
\end{equation*}
According to Proposition \ref{p27}(ii) and interpolation inequality,  we have
\begin{equation*}\label{318}
\begin{split}&\Big\|g_{+}(c^{+},c^{-})\partial_{i}c^{+}
-\tilde{g}_{+}(c^{+},c^{-})\partial_{i}c^{-}\Big\|_{L^1_{T_0}(\dot{B}^{\f{N}{2}-1}_{2,1})}
\\&\lesssim\|g_{+}(c^{+},c^{-})\|_{L^2_{T_0}(\dot{B}^{\frac{N}{2}}_{2,1})}
\|\partial_{i}c^{+}\|_{L^2_{T_0}(\dot{B}^{\f{N}{2}-1}_{2,1})}+\|\tilde{g}_{+}(c^{+},c^{-})\|_{L^2_{T_0}(\dot{B}^{\frac{N}{2}}_{2,1})}
\|\partial_{i}c^{-}\|_{L^2_{T_0}(\dot{B}^{\f{N}{2}-1}_{2,1})}
\\&\lesssim\|g_{+}(c^{+},c^{-})\|_{L^2_{T_0}(\dot{B}^{\frac{N}{2}}_{2,1})}
\|c^{+}\|_{L^2_{T_0}(\dot{B}^{\f{N}{2}}_{2,1})}+\|\tilde{g}_{+}(c^{+},c^{-})\|_{L^2_{T_0}(\dot{B}^{\frac{N}{2}}_{2,1})}
\|c^{-}\|_{L^2_{T_0}(\dot{B}^{\f{N}{2}}_{2,1})}
\\&\leq C_0\Big(\|c^{+}\|_{L^\infty}, \|c^{-}\|_{L^\infty}\Big)\Big(\|c^{+}\|_{L^2_{T_0}(\dot{B}^{\f{N}{2}}_{2,1})}+\|c^{-}\|_{L^2_{T_0}(\dot{B}^{\f{N}{2}}_{2,1})}\Big)^{2}
\\&\lesssim M^2\eta^2.
\end{split}
\end{equation*}
Hence, we gather that
\begin{equation}\label{4.3}
\|H_{21}\|_{L^1_{T_0}(\dot{B}^{\f{N}{2}-1}_{2,1})}\le
C(M^2\eta^2+M^3\eta^3).
\end{equation}
Similarly,  we also have
\begin{equation}\label{4.4}
\|H_{41}\|_{L^1_{T_0}(\dot{B}^{\f{N}{2}-1}_{2,1})}\le
C(M^2\eta^2+M^3\eta^3).
\end{equation}
 Substituting
\eqref{4.2}-\eqref{4.4} into  \eqref{4.1}, we obtain that
\begin{equation*}\label{322}
\begin{split}
X(T_0)\leq
C_1e^{C_1M\eta}\Big(\eta+M^2\eta^2+M^3\eta^3\Big).
\end{split}
\end{equation*}
Choose $M=8C_1$, for sufficiently small $\eta$  such that
\begin{equation*}
e^{C_1M \eta}\leq 2, \quad
\left(1+M\eta\right)M^2\eta\leq1,
\end{equation*}
which implies that
$$X(T_0)\le \frac{M\eta}{2}. $$ This is a
contradiction with the definition of $T_0$. As a consequence, we
conclude that $T_0=\infty$. Based on the above global uniform estimates,   employing a classical Friedrich's approximation and compactness method
(cf. \cite{Dan1,Dan3,Dan2}), we can establish the  global existence
of  strong solutions of the system \eqref{equ:CTFS1}-\eqref{eq:initial data}. Here,
we omit it. This completes the proof of the existence   of a global solution to the the system \eqref{equ:CTFS1}-\eqref{eq:initial data}  in Theorem \ref{th:main1}.
\section{Uniqueness}
\ \ \ \ \ In this section, we prove the uniqueness of the solution for the
system \eqref{equ:CTFS1}-\eqref{eq:initial data}. First, let us recall the Osgood Lemma
(see \cite{fleet}), which allows us to deduce uniqueness of the
solution in the critical case.
\begin{Lemma}\cite{fleet}\label{Lem:osgood}
 Let $f\geq 0$ be a measurable function, $\gamma$
be a locally integrable function and $\mu$ be a positive, continuous
 and non decreasing function which verifies the following condition
$$\int_0^1\frac{dr}{\mu (r)}=+\infty.$$
\noindent Let also $a$ be a positive real number and let $f$
satisfy
 the inequality
$$f (t)\leq a +\int _0^t\gamma (s)\mu(f (s))ds.$$
Then,

\noindent(i)\ if $a$ is equal to zero, the function $f$ vanishes;

\noindent (ii)\ if  $a$ is not zero, then we have
$$
-\mathcal M (f (t))+\mathcal M (a)\leq \int_0^t \gamma (s)ds,
\hskip0.5cm\text{with}\hskip0.5cm \mathcal M
(x)=\int\limits_x^1\frac{dr}{\mu (r)}\cdot
$$
\end{Lemma}
Next, we need the following result of logarithmic interpolation.
\begin{Lemma}\cite{Dan5}\label{Lem:loginequ} Let $s\in \mathbb{R}$. Then for any $1\le p,\rho\le+\infty$ and $0<\epsilon\le 1$,
we have \beno \|g\|_{\widetilde{L}^\rho_T(\dot{B}^s_{p,1})}\lesssim\frac{\|g\|_{\widetilde{L}^\rho_T(\dot{B}^s_{p,\infty})}}{\epsilon}
\log\Bigl(e+\frac{\|g\|_{\widetilde{L}^\rho_T(\dot{B}^{s-\epsilon}_{p,\infty})}
+\|g\|_{\widetilde{L}^\rho_T(\dot{B}^{s+\epsilon}_{p,\infty})}}
{\|g\|_{\widetilde{L}^\rho_T(\dot{B}^s_{p,\infty})}}\Bigr). \eeno
\end{Lemma}
We assume that
$(c^{+}_{1},\,u^{+}_{1},\,c^{-}_{1},\,u^{-}_{1})$,  $(c^{+}_{2},\,u^{+}_{2},\,c^{-}_{2},\,u^{-}_{2})$
are two solutions of the system \eqref{3.1} with the
same initial data satisfying \eqref{1.6A}. 
Observe that $\partial_tc^{\pm}_i\in L^1_{loc}(\dot{B}_{2,1}^{\frac{N}{2}-1})$, hence $c^{\pm}_i\in C(\dot{B}_{2,1}^{\frac{N}{2}-1})\cap L^\infty(\dot{B}_{2,1}^{\frac{N}{2}})(i=1,2)$. This entails $c^{\pm}_i\in C([0,\infty)\times \mathbb{R}^N)$. On the other
hand, if $\eta$ is sufficiently small, we have
$$|c^{\pm}_1(t,x)|\le \frac{1}{4}
\quad\textrm{for all}\quad t\ge 0\textrm{ and }x\in\mathbb{R}^N.$$
The continuity in time for $c^{\pm}_2$ thus yields the existence of a time $T>0$ such that
$$\|c^{\pm}_i(t)\|_{L^\infty}\le \frac{1}{2} \quad\textrm{for}\quad i=1,2\textrm{  and  }t\in[0,T].$$
From the embedding theorem and \eqref{1.6A}, we have
\begin{equation}\label{5.0}\|c^{\pm}_i\|_{\widetilde{L}^\infty_t(\dot B^{1}_{N,1})}\le \eta.\end{equation}
Set $\delta c^{+}=c^{+}_{1}-c^{+}_{2}$, $\delta u^{+}=u^{+}_{1}-u^{+}_{2}$, $\delta c^{-}=c^{-}_{1}-c^{-}_{2}$ and $\delta u^{-}=u^{-}_{1}-u^{-}_{2}$.
Then $(\delta c^{+},\delta u^{+}, \delta c^{-}, \delta u^{-})$ satisfies the following system
\begin{equation}\label{5.1}
\left\{
\begin{aligned}{}
&\p_t \delta c^{+}+{u}^{+}_{2}\cdot\nabla \delta c^{+}=-\textrm{div}\delta u^{+}+\delta H_{11},\\
&\p_t\delta u^{+}-\nu_{1}^{+}\Delta \delta u^{+}
-\nu_{2}^{+}\nabla\textrm{div} \delta u^{+}=-\delta u^{+}\cdot\nabla u^{+}_{1}- u^{+}_{2}\cdot\nabla \delta u^{+}-\beta_{1}\nabla \delta c^{+}
-\beta_{2}\nabla \delta c^{-}+\delta H_{21},\\
&\p_t \delta c^{-}+{u}^{-}_{2}\cdot\nabla \delta c^{-}=-\textrm{div}\delta u^{-}+\delta H_{31},\\
&\p_t\delta u^{-}-\nu_{1}^{-}\Delta \delta u^{-}
-\nu_{2}^{-}\nabla\textrm{div} \delta u^{-}=-\delta u^{-}\cdot\nabla u^{-}_{1}- u^{-}_{2}\cdot\nabla \delta u^{-}-\beta_{3}\nabla \delta c^{+}
-\beta_{4}\nabla \delta c^{-}+\delta H_{41},\\
&(\delta c^{+},\delta u^{+},\delta c^{-},\delta u^{-})|_{t=0}=0,
\end{aligned}
\right.
\end{equation}
where \begin{align*}\delta H_{11}&=H_{11}(c^{+}_{1}, u^{+}_{1})-H_{11}(c^{+}_{2}, u^{+}_{2}), & \qquad \delta H_{21}&=H_{11}(c^{+}_{1}, u^{+}_{1}, c^{-}_{1})-H_{11}(c^{+}_{2}, u^{+}_{2}, c^{-}_{2}),\\
\delta H_{31}&=H_{31}(c^{-}_{1}, u^{-}_{1})-H_{31}(c^{-}_{2}, u^{-}_{2}),& \qquad \delta H_{41}&=H_{41}(c^{+}_{1}, u^{-}_{1}, c^{-}_{1})-H_{41}(c^{+}_{2}, u^{-}_{2}, c^{-}_{2}).\end{align*}
In what follows, we set $U^{\pm}_i(t)=\int_0^t\|u^{\pm}_i(\tau)\|_{\dot{B}^{
\frac{N}{2}+1}_{2,1}}d\tau$ for $i=1,2$, and denote by $A_T$ a
constant depending  on $\|c^{\pm}_1\|_{\widetilde{L}^\infty_T(\dot B^{
\frac{N}{2}}_{2,1})}$ and $\|c^{\pm}_2\|_{\widetilde{L}^\infty_T(\dot B^{
\frac{N}{2}}_{2,1})}$. To begin with, we shall prove uniqueness on the time interval $[0, T]$ by
estimating $(\delta c^{+}, \delta u^{+},\delta c^{-}, \delta u^{-})$ in the following functional space:
\begin{equation*}
\begin{split}F_{T}&=L^{\infty}\big([0,T]; \dot{B}_{N,\infty}^{0}\big)\times \big(L^{\infty}([0,T]; \dot{B}_{N,\infty}^{-1})\cap \tilde{L}^{1}([0,T]; \dot{B}_{N,\infty}^{1})\big)^{N}
\\ &\quad\times L^{\infty}\big([0,T]; \dot{B}_{N,\infty}^{0}\big)\times \big(L^{\infty}([0,T]; \dot{B}_{N,\infty}^{-1})\cap \tilde{L}^{1}([0,T]; \dot{B}_{N,\infty}^{1})\big)^{N}.
\end{split}
\end{equation*}
We apply Proposition \ref{Prop:transport} to get for any $t\in
[0,T]$,
\begin{equation}
\begin{split}\label{eq:dens-diffe2} \|\delta
c^{+}(t)\|_{\dot{B}^{0}_{N,\infty}}&\lesssim e^{CU^{+}_2(t)}\int_0^t\big(\|\delta
H_{11}(\tau)\|_{\dot{B}^{0}_{N,\infty}}+ \|\textrm{div}\delta u^{+}(\tau)\|_{\dot{B}^{0}_{N,\infty}}\big)d\tau,\\
\|\delta
c^{-}(t)\|_{\dot{B}^{0}_{N,\infty}}&\lesssim e^{CU^{-}_2(t)}\int_0^t\big(\|\delta
H_{31}(\tau)\|_{\dot{B}^{0}_{N,\infty}}+ \|\textrm{div}\delta u^{-}(\tau)\|_{\dot{B}^{0}_{N,\infty}}\big)d\tau.
\end{split}
\end{equation}
From Proposition \ref{p26}, we have
 \beno &&\|\delta
H_{11}(\tau)\|_{\dot{B}^{0}_{N,\infty}}+ \|\textrm{div}\delta u^{+}(\tau)\|_{\dot{B}^{0}_{N,\infty}}\\
&&\quad\lesssim \|u^{+}_2(\tau)\|_{\dot{B}^{2}_{N,1}}\|\delta
c^{+}(\tau)\|_{\dot{B}^{0}_{N,\infty}}+\big(1+\|c^{+}_1(\tau)\|_{\dot{B}^{1}_{N,1}}\big)\|\delta
u^{+}(\tau)\|_{\dot{B}^{1}_{N,1}}, \eeno
\beno &&\|\delta
H_{31}(\tau)\|_{\dot{B}^{0}_{N,\infty}}+ \|\textrm{div}\delta u^{-}(\tau)\|_{\dot{B}^{0}_{N,\infty}}\\
&&\quad\lesssim\|u^{-}_2(\tau)\|_{\dot{B}^{2}_{N,1}}\|\delta
c^{-}(\tau)\|_{\dot{B}^{0}_{N,\infty}}+\big(1+\|c^{-}_1(\tau)\|_{\dot{B}^{1}_{N,1}}\big)\|\delta
u^{-}(\tau)\|_{\dot{B}^{1}_{N,1}}. \eeno
Plugging the above two  inequalities
into (\ref{eq:dens-diffe2}), we get by Gronwall's inequality that
\begin{equation}
\begin{split}\label{eq:dens-diffe3} \|(\delta c^{+},\delta c^{-})(t)\|_{\dot{B}^{0}_{N,\infty}}\lesssim
e^{CU_2(t)}\int_0^t\big(1+\|(c^{+}_1,c^{-}_1) (\tau)\|_{\dot{B}^{1}_{N,1}}\big)\|(\delta
u^{+}, \delta
u^{-})(\tau)\|_{\dot{B}^{1}_{N,1}}d\tau,
\end{split}
\end{equation}
where $ U_2(t)=\max\big\{U^{+}_2(t),U^{-}_2(t)\big\}.$

Applying Remark \ref{2.14} to the second equation and  the fourth equation of \eqref{5.1},  we have
\begin{equation}
\begin{split}\label{eq:dens-diffe4}
 &\|\delta u^{+}\|_{\widetilde{L}^1_t(\dot B^{1}_{N,\infty})}+\|\delta
u^{+}\|_{\widetilde{L}^2_t(\dot B^{0}_{N,\infty})}\\&\lesssim  \|\delta
H_{21}\|_{\widetilde{L}^1_t(\dot B^{-1}_{N,\infty})}+ \|\delta u^{+}\cdot\nabla u^{+}_{1}\|_{\widetilde{L}^1_t(\dot B^{-1}_{N,\infty})}
+ \|u^{+}_{2}\cdot\nabla \delta u^{+}\|_{\widetilde{L}^1_t(\dot B^{0}_{N,\infty})}
\\&\quad+\|\nabla \delta c^{+}\|_{\widetilde{L}^1_t(\dot B^{-1}_{N,\infty})}
+\|\nabla \delta c^{-}\|_{\widetilde{L}^1_t(\dot B^{-1}_{N,\infty})},
\end{split}
\end{equation}
\begin{equation}
\begin{split}\label{eq:dens-diffe5}
 &\|\delta
u^{-}\|_{\widetilde{L}^1_t(\dot B^{1}_{N,\infty})}+\|\delta
u^{-}\|_{\widetilde{L}^2_t(\dot B^{0}_{N,\infty})} \\ &\lesssim \|\delta
H_{41}\|_{\widetilde{L}^1_t(\dot B^{-1}_{N,\infty})}+ \|\delta u^{-}\cdot\nabla u^{-}_{1}\|_{\widetilde{L}^1_t(\dot B^{-1}_{N,\infty})}
+ \|u^{-}_{2}\cdot\nabla \delta u^{-}\|_{\widetilde{L}^1_t(\dot B^{0}_{N,\infty})}
\\&\quad+\|\nabla \delta c^{+}\|_{\widetilde{L}^1_t(\dot B^{-1}_{N,\infty})}
+\|\nabla \delta c^{-}\|_{\widetilde{L}^1_t(\dot B^{-1}_{N,\infty})}.
\end{split}
\end{equation}
Employing  Proposition \ref{p26} and Proposition \ref{p27},  we have
\begin{equation}
\begin{split}\label{eq:dens-diffe6}&\|\delta
H_{21}\|_{\widetilde{L}^1_t(\dot B^{-1}_{N,\infty})}+ \|\delta u^{+}\cdot\nabla u^{+}_{1}\|_{\widetilde{L}^1_t(\dot B^{-1}_{N,\infty})}
+ \|u^{+}_{2}\cdot\nabla \delta u^{+}\|_{\widetilde{L}^1_t(\dot B^{0}_{N,\infty})}
\\&\quad+\|\nabla \delta c^{+}\|_{\widetilde{L}^1_t(\dot B^{-1}_{N,\infty})}
+\|\nabla \delta c^{-}\|_{\widetilde{L}^1_t(\dot B^{-1}_{N,\infty})}
\\&\lesssim A_T\int_{0}^{T}\Big(1+\|u^{+}_2(\tau)\|_{\dot B^{2}_{N,1}}\Big) \|(\delta
c^{+},\delta c^{-})(\tau)\|_{\dot
B^{0}_{N,\infty}}d\tau+\|(u^{+}_1,u^{+}_2)\|_{\widetilde{L}^2_t(\dot B^{1}_{N,1})}\|\delta
u^{+}\|_{\widetilde{L}^2_t(\dot B^{0}_{N,\infty})}
\\&\qquad+ A_T\|(c^{+}_1,c^{-}_1)\|_{\widetilde{L}^\infty_t(\dot B^{1}_{N,1})}
\|\delta u^{+}\|_{\widetilde{L}^1_t(\dot B^{1}_{N,\infty})},
\end{split}
\end{equation}
\begin{equation}
\begin{split}\label{eq:dens-diffe7}&\|\delta
H_{41}\|_{\widetilde{L}^1_t(\dot B^{-1}_{N,\infty})}+ \|\delta u^{-}\cdot\nabla u^{-}_{1}\|_{\widetilde{L}^1_t(\dot B^{-1}_{N,\infty})}
+ \|u^{-}_{2}\cdot\nabla \delta u^{-}\|_{\widetilde{L}^1_t(\dot B^{0}_{N,\infty})}
\\&\quad+\|\nabla \delta c^{+}\|_{\widetilde{L}^1_t(\dot B^{-1}_{N,\infty})}
+\|\nabla \delta c^{-}\|_{\widetilde{L}^1_t(\dot B^{-1}_{N,\infty})}
\\&\lesssim A_T\int_{0}^{T}\Big(1+\|u^{-}_2(\tau)\|_{\dot B^{2}_{N,1}}\Big) \|(\delta
c^{+},\delta c^{-})(\tau)\|_{\dot
B^{0}_{N,\infty}}d\tau+\|(u^{-}_1,u^{-}_2)\|_{\widetilde{L}^2_t(\dot B^{1}_{N,1})}\|\delta
u^{-}\|_{\widetilde{L}^2_t(\dot B^{0}_{N,\infty})}
\\&\qquad+ A_T\|(c^{+}_1,c^{-}_1)\|_{\widetilde{L}^\infty_t(\dot B^{1}_{N,1})}
\|\delta u^{-}\|_{\widetilde{L}^1_t(\dot B^{1}_{N,\infty})}.
\end{split}
\end{equation}
 Employing   \eqref{5.0}  and taking  $T$ small enough such that
\begin{equation}\label{eq:dens-diffe8}A_T\|(c^{+}_1,c^{-}_1)\|_{\widetilde{L}^\infty_t(\dot B^{1}_{N,1})}+ \|(u^{-}_1,u^{-}_2)\|_{\widetilde{L}^2_t(\dot B^{1}_{N,1})}\ll 1.\end{equation}
Combining with \eqref{eq:dens-diffe4}-\eqref{eq:dens-diffe8}, we have
\begin{equation}
\begin{split}\label{eq:dens-diffe9}
 \|(\delta u^{+},\delta u^{-})\|_{\widetilde{L}^1_t(\dot B^{1}_{N,\infty})}
 \lesssim A_T\int_{0}^{T}\Big(1+\|(u^{+}_2,u^{-}_2)(\tau)\|_{\dot B^{2}_{N,1}}\Big) \|(\delta
c^{+},\delta c^{-})(\tau)\|_{\dot
B^{0}_{N,\infty}}d\tau.
\end{split}
\end{equation}
From Lemma \ref{Lem:loginequ}, it follows that
\begin{equation*}
\begin{split} \|(\delta
u^{+},\delta u^{-})\|_{L^1_t(\dot B^{1}_{N,1})}&\lesssim\|(\delta
u^{+},\delta u^{-})\|_{\widetilde{L}^1_t(\dot{B}^1_{N,\infty})}\\
&\times\log\Big(e+\frac{\|(\delta
u^{+},\delta u^{-})\|_{\widetilde{L}^1_t(\dot{B}^{0}_{N,\infty})} +\|(\delta
u^{+},\delta u^{-})\|_{\widetilde{L}^1_t(\dot{B}^{2}_{N,\infty})}} {\|(\delta
u^{+},\delta u^{-})\|_{\widetilde{L}^1_t(\dot{B}^1_{N,\infty})}}\Big),
\end{split}
\end{equation*} which
together with  \eqref{eq:dens-diffe3} yields that for any $t\in
[0,T]$,
\begin{equation}\label{eq:dens-diffe10}
\begin{split}\|(\delta c^{+},\delta c^{-})(t)\|_{\dot{B}^{0}_{N,\infty}}
&\lesssim\|(\delta
u^{+},\delta u^{-})\|_{\widetilde{L}^1_t(\dot{B}^1_{N,\infty})}\\
&\times\log\Big(e+\frac{\|(\delta
u^{+},\delta u^{-})\|_{\widetilde{L}^1_t(\dot{B}^{0}_{N,\infty})} +\|(\delta
u^{+},\delta u^{-})\|_{\widetilde{L}^1_t(\dot{B}^{2}_{N,\infty})}} {\|(\delta
u^{+},\delta u^{-})\|_{\widetilde{L}^1_t(\dot{B}^1_{N,\infty})}}\Big).
\end{split}
\end{equation}
Combining
the above inequality  with  \eqref {eq:dens-diffe9}, we have
\begin{equation}\label{eq:dens-diffe10}
\begin{split}\|(\delta u^{+},\delta u^{-})\|_{\widetilde{L}^1_t(\dot B^{1}_{N,\infty})}
\lesssim\displaystyle\int_0^t\bigl(1+\|(u_1,u_2)(\tau)\|_{\dot
B^{2}_{N,1}}\bigr)V(\tau)
\log\Bigl(e+\frac{C_T} {V(\tau)}\Bigr) d\tau,
\end{split}
\end{equation}
 where $V(t)=\|(\delta
u^{+},\delta u^{-})\|_{\widetilde{L}^1_t(\dot{B}^1_{N,\infty})},\quad C_T=\|(\delta
u^{+},\delta u^{-})\|_{\widetilde{L}^1_t(\dot{B}^{0}_{N,\infty})} +\|(\delta
u^{+},\delta u^{-})\|_{\widetilde{L}^1_t(\dot{B}^{2}_{N,\infty})}$. Notice that
$1+\|(u^{+}_2,u^{-}_2)(t)\|_{\dot B^{2}_{N,1}}$ is integrable on $[0,T]$, and \beno \int_0^1\f {dr}
{r\log(e+C_Tr^{-1})}dr=+\infty. \eeno From the  Osgood lemma \ref
{Lem:osgood}, we conclude that $(\delta u^{+},\delta u^{-})=0$ on
$[0,T]$. This gives by inequality
(\ref{eq:dens-diffe3}), that $(\delta c^{+},\delta c^{-})=0$.
A standard continuity argument ensures that $(c^{+}_{1},\,u^{+}_{1},\,c^{-}_{1},\,u^{-}_{1})=(c^{+}_{2},\,u^{+}_{2},\,c^{-}_{2},\,u^{-}_{2})$ on $[0,+\infty)$.

\section{Time decay estimates}
\ \ \ \ \ In this section, we will establish the time decay rates of the  global  strong solutions constructed in Theorem \ref{th:main1}.  We divide the proof into several steps.

\subsubsection*{Step 1: Low frequencies}
We first exhibit  the
smoothing properties of the system \eqref{3.1} in the low frequencies regime.  The
key to these remarkable properties is given by the following lemma.
\begin{Lemma}\label{lem6.1} Let $(c^{+},\,u^{+},\,c^{-},\,u^{-})$ be a solution of the system
\eqref{3.1}. Then, there exist two positive constants $c_0$ and $C$ depending only on $\beta_i(i=1,2,3,4)$ and $\nu^{\pm}_{i}(i=1,2)$
 respectively, such that the following inequality holds for all $t\geq0$,
\begin{equation}\label{6.1}
\begin{split}
\bigl\|(\dot{\Delta}_qc^{+},\,\dot{\Delta}_qu^{+},\,\dot{\Delta}_qc^{-},\,\dot{\Delta}_qu^{-}
)\bigr\|^\ell_{L_{2}}
& \leq C \Big( e^{-c_{0}2^{2q}t}\bigl\|(\dot{\Delta}_qc^{+}_{0},\,\dot{\Delta}_qu^{+}_{0},\,\dot{\Delta}_qc^{-}_{0},\,\dot{\Delta}_qu^{-}_{0}
)\bigr\|^\ell_{L_{2}}\\&\quad+\int_0^te^{-c2^{2q}(t-\tau)}\bigl\|(\dot{\Delta}_qH_{1},\,\dot{\Delta}_qH_{2},\,\dot{\Delta}_qH_{3},\,\dot{\Delta}_qH_{4}
)(\tau)\bigr\|^\ell_{L_{2}}d\tau\Big).
\end{split}
\end{equation}
\end{Lemma}
\noindent{\bf Proof.} By the same  derivation process of \eqref{171}, in the case  $H_{1}\equiv H_{2}\equiv H_{3}\equiv H_{4}\equiv0$,  we have
\begin{equation}\label{6.2}
\begin{split}
  \frac{1}{2}\frac{d}{dt}\alpha_q+c_02^{2q}\alpha_q\leq 0,
\end{split}\end{equation}
where
\begin{equation*}
\begin{split}
\alpha_q&\approx\|\Dot{\Delta}_qc^{+}\|_{L^2}+\|\Dot{\Delta}_qc^{-}\|_{L^2}+\|\nabla\Dot{\Delta}_qc^{+}\|_{L^2}
+\|\nabla\Dot{\Delta}_qc^{-}\|_{L^2}
+\|\Dot{\Delta}_qu^{+}\|_{L^2}+\|\Dot{\Delta}_qu^{-}\|_{L^2}\\
&\approx\bigl\|(\dot{\Delta}_qc^{+},\,\dot{\Delta}_qu^{+},\,\dot{\Delta}_qc^{-},\,\dot{\Delta}_qu^{-}
)\bigr\|^\ell_{L_{2}}.
\end{split}
\end{equation*}
Thus, \begin{equation}\label{6.3}
\begin{split}
\bigl\|(\dot{\Delta}_qc^{+},\,\dot{\Delta}_qu^{+},\,\dot{\Delta}_qc^{-},\,\dot{\Delta}_qu^{-}
)\bigr\|^\ell_{L_{2}} \leq Ce^{-c_{0}2^{2q}t}\bigl\|(\dot{\Delta}_qc^{+}_{0},\,\dot{\Delta}_qu^{+}_{0},\,\dot{\Delta}_qc^{-}_{0},\,\dot{\Delta}_qu^{-}_{0}
)\bigr\|^\ell_{L_{2}}.
\end{split}
\end{equation}
Furthermore, taking advantage of the Duhamel formula, we can readily deduce \eqref{6.1}.

Denoting by $A(D)$ the semi-group associated to the system \eqref{3.1}, we have for all $q\in\mathbb{Z},$
\begin{equation}
\label{low.1}
\left(\begin{array}{c}\ddq c^{+}(t)\\\ddq u^{+}(t)\\\ddq c^{-}(t)\\\ddq u^{-}(t)\end{array}\right)
=e^{tA(D)}\left(\begin{array}{c}\ddq c^{+}_0\\\ddq u^{+}_0\\\ddq c^{-}_0\\\ddq u^{-}_0\end{array}\right)
+\int_0^te^{(t-\tau)A(D)}\left(\begin{array}{c}\ddq H_1(\tau)\\\ddq H_2(\tau)\\\ddq H_3(\tau)\\\ddq H_4(\tau)\end{array}\right)d\tau.
\end{equation}
Based on  \eqref{6.3} and \eqref{low.1},  we get for all $q\leq q_0,$
\begin{align*}\|e^{tA(D)}\ddq U\|_{L^2}&\lesssim  e^{-c_{0}2^{2q}t}\|\ddq U\|_{L^2}.
\end{align*}
Hence, multiplying by $t^{\frac{N}{4}+\frac{s}{2}}2^{qs}$ and summing up on $q\leq q_{0}$, we get
\begin{equation}
\begin{split}
\label{low.2}
t^{\frac{N}{4}+\frac{s}{2}}\sum_{q\leq q_0}2^{qs}\|e^{tA(D)}\ddq U\|_{L^2}
&\lesssim
\sum_{q\leq q_0}2^{qs}e^{-c_{0}2^{2q}t}\|\ddq U\|_{L^2}t^{\frac{N}{4}+\frac{s}{2}}\\
&\lesssim
\sum_{q\leq q_0}2^{q(s+\frac{N}{2})}e^{-c_{0}2^{2q}t}\|\ddq U\|_{L^2}2^{q(-\frac{N}{2})}t^{\frac{N}{4}+\frac{s}{2}}\\
&\lesssim
\|U\|_{\dot B^{-\frac{N}{2}}_{2,\infty}}^\ell\sum_{q\leq q_0}2^{q(s+\frac{N}{2})}e^{-c_{0}2^{2q}t}t^{\frac{N}{4}+\frac{s}{2}}\\
&\lesssim
\|U\|_{\dot B^{-\frac{N}{2}}_{2,\infty}}^\ell\sum_{q\leq q_0}2^{q(s+\frac{N}{2})}e^{-c_{0}2^{2q}t}
t^{\frac{1}{2}(s+\frac{N}{2})}.
\end{split}
\end{equation}
As for any $\sigma>0$ there  exists a constant $C_\sigma$ so that
\begin{equation}\label{low.3}
\sup_{t\geq0}\sum_{q\in\mathbb{Z}}t^{\frac\sigma2}2^{q\sigma}e^{-c_{0}2^{2q}t}\leq C_\sigma.
\end{equation}
We get from \eqref{low.2} and \eqref{low.3} that for $s>-N/2,$
$$
\sup_{t\geq0}\, t^{\frac N4+\frac s2}\|e^{tA(D)}U\|_{\dot B^s_{2,1}}^\ell
\lesssim\|U\|_{\dot B^{-\frac{N}{2}}_{2,\infty}}^\ell.
$$
Furthermore, it is  obvious that  for $s>-N/2,$
$$
\|e^{tA(D)}U\|_{\dot B^s_{2,1}}^\ell
\lesssim \|U\|_{\dot B^{-\frac{N}{2}}_{2,\infty}}^\ell\sum_{q\leq q_0}2^{q(s+\frac{N}{2})}\lesssim\|U\|_{\dot B^{-\frac{N}{2}}_{2,\infty}}^\ell .
$$
Hence, setting $\langle t\rangle\eqdefa\sqrt{1+t^{2}}$, we get
\begin{equation}
\label{U}
\sup_{t\geq0}\, \langle t\rangle^{\frac N4+\frac s2}\|e^{tA(D)}U\|_{\dot B^s_{2,1}}^\ell
\lesssim\|U\|_{\dot B^{-\frac{N}{2}}_{2,\infty}}^\ell.
\end{equation}
Thus, from \eqref{6.1} and \eqref{U}, we have
\begin{equation}
\begin{split}
&\|(c^{+},\,u^{+},\,c^{-},\,u^{-}
)\|_{\dot B^s_{2,1}}^\ell\\
&\lesssim \sup_{t\geq0}\, \langle t\rangle^{-(\frac N4+\frac s2)}\|(c^{+}_{0},\,u^{+}_{0},\,c^{-}_{0},\,u^{-}_{0}
)\|_{\dot B^{-\frac{N}{2}}_{2,\infty}}^\ell+\int_{0}^{t}\langle t-\tau\rangle^{-(\frac N4+\frac s2)}\|(H_{1},\,H_{2},\,H_{3},\,H_{4})(\tau)
\|_{\dot B^{-\frac{N}{2}}_{2,\infty}}^\ell d\tau.
\end{split}
\end{equation}
We claim that for all $s\in(-N/2,2]$ and $t\geq0$,  then
\begin{equation}\label{s1234low1}
\int_{0}^{t}\langle t-\tau\rangle^{-(\frac N4+\frac s2)}\|(H_{1},\,H_{2},\,H_{3},\,H_{4})(\tau)
\|_{\dot B^{-\frac{N}{2}}_{2,\infty}}^\ell d\tau
\lesssim\langle t\rangle^{-(\frac N4+\frac s2)}
\Big(X^2(t)+D^2(t)+D^3(t)+D^4(t)\Big),
\end{equation}
where $X(t)$ and $D(t)$ have been defined in \eqref{1.6} and \eqref{1.9}, respectively.\medbreak
Owing to the embedding $L^{1}\hookrightarrow\dot B^{-\frac{N}{2}}_{2,\infty}$,  it suffices to prove \eqref{s1234low1}  with $\|(H_{1},H_{2},H_{3},H_{4})(\tau)\|_{L^1}^{\ell}$ instead of $\|(H_{1},H_{2},H_{3},H_{4})(\tau)\|_{\dot  B^{-\frac{N}{2}}_{2,\infty}}^{\ell}$.\medbreak
To bound the term with $H_{1}$, we use the following decomposition:
$$ H_{1}=u^{+}\cdot\nabla c^{+}+ c^{+}\, \div (u^{+})^\ell + c^{+}\,\div (u^{+})^h.$$
Now, from H\"{o}lder's inequality,  the embedding $\dot B^{0}_{2,1}\hookrightarrow L^{2}$, the definitions of $D(t), \alpha$ and Lemma \ref{lemma2.13}, one may write for all $s\in(\varepsilon-\frac N2,2]$,
\begin{equation}
\begin{split}\label{s11}
&\int_0^t\langle t-\tau\rangle^{-(\frac N4+\frac s2)}\|(u^{+}\cdot\nabla c^{+})(\tau)\|_{L^1}\,d\tau\\
&\quad\lesssim\int_0^t\langle t-\tau\rangle^{-(\frac N4+\frac s2)}\|u^{+}\|_{L^2}\|\nabla c^{+}\|_{L^2}\,d\tau\\
&\quad\lesssim\int_0^t\langle t-\tau\rangle^{-(\frac N4+\frac s2)}\|u^{+}\|_{\dot  B^{0}_{2,1}}\|\nabla c^{+}\|_{\dot  B^{0}_{2,1}}\,d\tau\\
&\quad\lesssim\int_0^t\langle t-\tau\rangle^{-(\frac N4+\frac s2)}\big(\|u^{+}\|_{\dot  B^{0}_{2,1}}^{\ell}+\|u^{+}\|_{\dot  B^{0}_{2,1}}^h\big)\big(\|\nabla c^{+}\|_{\dot  B^{0}_{2,1}}^{\ell}+\|\nabla c^{+}\|_{\dot  B^{0}_{2,1}}^{h}\big)\,d\tau\\
&\quad\lesssim\big(\sup_{0\leq\tau\leq t}\langle \tau\rangle^{\frac N4}\|u^{+}(\tau)\|_{\dot  B^{0}_{2,1}}^{\ell}\big)
\big(\sup_{0\leq\tau\leq t}\langle \tau\rangle^{\frac N4+\frac 12}\|\nabla c^{+}(\tau)\|_{\dot  B^{0}_{2,1}}^{\ell}\big)
\int_0^t\langle t-\tau\rangle^{-(\frac N4+\frac s2)}\langle \tau\rangle^{-(\frac N2+\frac 12)}\,d\tau\\
&\qquad+\big(\sup_{0\leq\tau\leq t}\langle \tau\rangle^{\frac N4}\|u^{+}(\tau)\|_{\dot  B^{0}_{2,1}}^{\ell}\big)
\big(\sup_{0\leq\tau\leq t}\langle \tau\rangle^\alpha\|\nabla c^{+}(\tau)\|_{\dot  B^{\frac{N}{2}-1}_{2,1}}^h\big)
\int_0^t\langle t-\tau\rangle^{-(\frac N4+\frac s2)}\langle \tau\rangle^{-(\frac N4+\alpha)}\,d\tau\\
&\qquad+\big(\sup_{0\leq\tau\leq t}\langle \tau\rangle^\alpha\|u^{+}(\tau)\|_{\dot  B^{\frac{N}{2}-1}_{2,1}}^h\big)
\big(\sup_{0\leq\tau\leq t}\langle \tau\rangle^{\frac N4+\frac 12}\|\nabla c^{+}(\tau)\|_{\dot  B^{0}_{2,1}}^{\ell}\big)
\int_0^t\langle t-\tau\rangle^{-(\frac N4+\frac s2)}\langle \tau\rangle^{-(\alpha+\frac N4+\frac 12)}\,d\tau\\
&\qquad+\big(\sup_{0\leq\tau\leq t}\langle \tau\rangle^\alpha\|u^{+}(\tau)\|_{\dot  B^{\frac{N}{2}-1}_{2,1}}^h\big)
\big(\sup_{0\leq\tau\leq t}\langle \tau\rangle^{\alpha}\|\nabla c^{+}(\tau)\|_{\dot  B^{\frac{N}{2}-1}_{2,1}}^h\big)
\int_0^t\langle t-\tau\rangle^{-(\frac N4+\frac s2)}\langle \tau\rangle^{-2\alpha}\,d\tau\\
&\quad\lesssim D^{2}(t)\int_0^t\langle t-\tau\rangle^{-(\frac N4+\frac s2)}\langle \tau\rangle^{-\min({\frac N2+\frac 12,\alpha+\frac{N}{4},2\alpha})}\,d\tau\\
&\quad\lesssim \langle t\rangle^{-(\frac N4+\frac s2)} D^{2}(t).
\end{split}
\end{equation}
The term  $c^{+}\,\div (u^{+})^\ell$ may be treated along the same lines,  and we have
\begin{equation}
\begin{split}
\label{s12}
&\int_0^t\langle t-\tau\rangle^{-(\frac N4+\frac s2)}\|c^{+}\div (u^{+})^{\ell} \|_{L^1}\,d\tau\\
&\quad\lesssim\int_0^t\langle t-\tau\rangle^{-(\frac N4+\frac s2)}\big(\|c^{+}\|_{\dot  B^{0}_{2,1}}^{\ell}+\|c^{+}\|_{\dot  B^{0}_{2,1}}^h\big)\|\nabla u^{+}\|_{\dot  B^{0}_{2,1}}^{\ell}\,d\tau\\
&\quad\lesssim\big(\sup_{0\leq\tau\leq t}\langle \tau\rangle^{\frac N4}\|c^{+}(\tau)\|_{\dot  B^{0}_{2,1}}^{\ell}\big)
\big(\sup_{0\leq\tau\leq t}\langle \tau\rangle^{\frac{N}{4}+\frac 12}\| u^{+}(\tau)\|_{\dot  B^{1}_{2,1}}^\ell\big)
\int_0^t\langle t-\tau\rangle^{-(\frac N4+\frac s2)}\langle \tau\rangle^{-(\frac{N}{2}+\frac 12)}\,d\tau\\
&\qquad+\big(\sup_{0\leq\tau\leq t}\langle \tau\rangle^\alpha\|c^{+}(\tau)\|_{\dot  B^{\frac{N}{2}}_{2,1}}^h\big)
\big(\sup_{0\leq\tau\leq t}\langle \tau\rangle^{\frac N4+\frac 12}\| u^{+}(\tau)\|_{\dot  B^{1}_{2,1}}^{\ell}\big)
\int_0^t\langle t-\tau\rangle^{-(\frac N4+\frac s2)}\langle \tau\rangle^{-(\alpha+\frac N4+\frac 12)}\,d\tau\\
&\quad\lesssim D^{2}(t)\int_0^t\langle t-\tau\rangle^{-(\frac N4+\frac s2)}\langle \tau\rangle^{-\min({\frac N2+\frac 12,\alpha+\frac{N}{4}+\frac 12})}\,d\tau\\
&\quad\lesssim \langle t\rangle^{-(\frac N4+\frac s2)} D^{2}(t).
\end{split}
\end{equation}
Regarding the term with $c^{+}\,\div (u^{+})^h,$ we get for all $t\geq2$ that
\begin{equation*}
\begin{split}
&\int_0^t\langle t-\tau\rangle^{-(\frac N4+\frac s2)}\|c^{+}\div (u^{+})^h(\tau)\|_{L^1}\,d\tau\\
&\quad\lesssim\int_0^t\langle t-\tau\rangle^{-(\frac N4+\frac s2)} \|c^{+}(\tau)\|_{\dot  B^{0}_{2,1}}\|\div u^{+}(\tau)\|^h_{\dot  B^{0}_{2,1}}\,d\tau\\
&\quad\lesssim\int_0^1\langle t-\tau\rangle^{-(\frac N4+\frac s2)} \|c^{+}(\tau)\|_{\dot  B^{0}_{2,1}}\|\div u^{+}(\tau)\|^h_{\dot  B^{0}_{2,1}}\,d\tau\\
&\qquad+\int_1^t\langle t-\tau\rangle^{-(\frac N4+\frac s2)} \|c^{+}(\tau)\|_{\dot  B^{0}_{2,1}}\|\div u^{+}(\tau)\|^h_{\dot  B^{0}_{2,1}}\,d\tau\\
&\quad\eqdefa I_{1}+I_{2}.
\end{split}
\end{equation*}
From the definitions of $X(t)$ and $D(t)$, we  obtain
\begin{align*}
I_{1}&\lesssim\langle t\rangle^{-(\frac N4+\frac s2)}
\sup _{0\leq\tau\leq1}\|c^{+}  (\tau)\|_{\dot  B^{0}_{2,1}}\int_0^1\|\div u^{+}(\tau)\|^h_{\dot  B^{0}_{2,1}}\,d\tau\\
&\lesssim\langle t\rangle^{-(\frac N4+\frac s2)}
\sup _{0\leq\tau\leq1}\|c^{+}(\tau)\|_{\dot  B^{0}_{2,1}}\int_0^1\| u^{+}(\tau)\|^h_{\dot  B^{\frac{N}{2}+1}_{2,1}}\,d\tau\\
&\lesssim\langle t\rangle^{-(\frac N4+\frac s2)} D(1)X(1),
\end{align*}
and,  using the fact that $\langle \tau\rangle\approx\tau$ when $\tau\geq1$, we get
\begin{align*}
I_{2} &\lesssim\int_1^t\langle t-\tau\rangle^{-(\frac N4+\frac s2)}
\big(\|c^{+}(\tau)\|_{\dot  B^{0}_{2,1}}^{\ell}+\|c^{+}(\tau)\|_{\dot  B^{0}_{2,1}}^{h}\big)\|\div u^{+}(\tau)\|^h_{\dot  B^{0}_{2,1}}\,d\tau\\
&\lesssim\big(\sup_{1\leq\tau\leq t}\langle \tau\rangle^{\frac N4}\|c^{+}(\tau)\|_{\dot  B^{0}_{2,1}}^{\ell}\big)
\big(\sup_{1\leq\tau\leq t}\| \tau\nabla u^{+}(\tau)\|^{h}_{\dot  B^{\frac{N}{2}}_{2,1}}\big)
\int_1^t\langle t-\tau\rangle^{-(\frac N4+\frac s2)}\langle \tau\rangle^{-(\frac N4+1)}\,d\tau\\
&\quad+\big(\sup_{1\leq\tau\leq t}\langle \tau\rangle^{\alpha}\|c^{+}(\tau)\|_{\dot  B^{\frac{N}{2}}_{2,1}}^h\big)
\big(\sup_{1\leq\tau\leq t}\| \tau\nabla u^{+}(\tau)\|^{h}_{\dot  B^{\frac{N}{2}}_{2,1}}\big)
\int_1^t\langle t-\tau\rangle^{-(\frac N4+\frac s2)}\langle \tau\rangle^{-(\alpha+1)}\,d\tau\\
&\lesssim D^{2}(t)\int_1^t\langle t-\tau\rangle^{-(\frac N4+\frac s2)}\langle \tau\rangle^{-\min({\alpha+1,\frac{N}{4}+1})}\,d\tau\\
&\lesssim \langle t\rangle^{-(\frac 34+\frac s2)} D^{2}(t).
\end{align*}
Thus, for $t\geq2,$ we conclude that
\begin{align}
\label{s112}
&\int_0^t\langle t-\tau\rangle^{-(\frac N4+\frac s2)}\|c^{+}\div (u^{+})^h(\tau)\|_{L^1}\,d\tau\nonumber\\
&\quad\lesssim\langle t\rangle^{-(\frac N4+\frac s2)} \big(D^{2}(t)+X^{2}(t)\big).
\end{align}
The case $t\leq2$ is obvious as $\langle t\rangle\approx1$ and
$\langle t-\tau\rangle\approx1$ for $0\leq\tau\leq t\leq 2$, and
\begin{equation}
\begin{split}
\label{s14}
&\int_0^t\|c^{+}\,\div (u^{+})^h\|_{L^1}\,d\tau\\
&\quad\lesssim \|c^{+}\|_{L^\infty_t(L^2)}\|\div(u^{+})\|^h_{L_t^1(L^2)}\\
&\quad\lesssim \|c^{+}\|_{L^\infty_t(\dot  B^{0}_{2,1})}
\|\div (u^{+})\|^h_{L_t^1(\dot  B^{0}_{2,1})}\\
&\quad\lesssim \|c^{+}\|_{L^\infty_t(\dot  B^{0}_{2,1})}
\| u^{+}\|^h_{L_t^1(\dot  B^{\frac{N}{2}+1}_{2,1})}\\
&\quad\lesssim X(t)D(t).
\end{split}
\end{equation}
From \eqref{s11}-\eqref{s14}, we get
\begin{equation*}\label{s1}
\int_0^t\langle t-\tau\rangle^{-(\frac{N}{4}+\frac{s}{2})} \big\|H_{1}(\tau)\big\|_{\dot B^{-\frac{N}{2}}_{2,\infty}}^\ell d\tau
\lesssim\langle t\rangle^{-(\frac N4+\frac s2)}
\big(X^2(t)+D^2(t)\big).
\end{equation*}
The term  $H_{3}$ may be treated along the same lines,  and we obtain
\begin{equation*}\label{s3}
\int_0^t\langle t-\tau\rangle^{-(\frac{N}{4}+\frac{s}{2})} \big\|H_{3}(\tau)\big\|_{\dot B^{-\frac{N}{2}}_{2,\infty}}^\ell d\tau
\lesssim\langle t\rangle^{-(\frac N4+\frac s2)}
\big(X^2(t)+D^2(t)\big).
\end{equation*}
Next, to bound the first term of $H_{2}^{i}$, we write that
\begin{equation}
\begin{split}\label{s21}
&\int_0^t\langle t-\tau\rangle^{-(\frac N4+\frac s2)}\|g_{+}(c^{+},c^{-})\partial_{i}c^{+}(\tau)\|_{L^1}\,d\tau\\
&\quad\lesssim\int_0^t\langle t-\tau\rangle^{-(\frac N4+\frac s2)}
\|g_{+}(c^{+},c^{-})\|_{L^2}\|\nabla c^{+}\|_{L^2}\,d\tau\\
&\quad\lesssim\int_0^t\langle t-\tau\rangle^{-(\frac N4+\frac s2)}\|g_{+}(c^{+},c^{-})\|_{\dot  B^{0}_{2,1}}\|\nabla c^{+}\|_{\dot  B^{0}_{2,1}}\,d\tau\\
&\quad\lesssim\int_0^t\langle t-\tau\rangle^{-(\frac N4+\frac s2)}\|(c^{+},c^{-})\|_{\dot  B^{0}_{2,1}}\|\nabla c^{+}\|_{\dot  B^{0}_{2,1}}\,d\tau\\
&\quad\lesssim\int_0^t\langle t-\tau\rangle^{-(\frac N4+\frac s2)}
\big(\|(c^{+},c^{-})\|_{\dot B^{0}_{2,1}}^{\ell}
+\|(c^{+},c^{-})\|_{\dot  B^{0}_{2,1}}^h\big)
\big(\|\nabla c^{+}\|_{\dot  B^{0}_{2,1}}^{\ell}+\|\nabla c^{+}\|_{\dot  B^{0}_{2,1}}^{h}\big)\,d\tau\\
&\quad\lesssim\big(\sup_{0\leq\tau\leq t}\langle \tau\rangle^{\frac N4}\|(c^{+},c^{-})(\tau)\|_{\dot  B^{0}_{2,1}}^{\ell}\big)
\big(\sup_{0\leq\tau\leq t}\langle \tau\rangle^{\frac N4+\frac 12}\|\nabla c^{+}(\tau)\|_{\dot  B^{0}_{2,1}}^{\ell}\big)
\int_0^t\langle t-\tau\rangle^{-(\frac N4+\frac s2)}\langle \tau\rangle^{-(\frac N2+\frac 12)}\,d\tau\\
&\qquad+\big(\sup_{0\leq\tau\leq t}\langle \tau\rangle^{\frac N4}\|(c^{+},c^{-})(\tau)\|_{\dot  B^{0}_{2,1}}^{\ell}\big)
\big(\sup_{0\leq\tau\leq t}\langle \tau\rangle^\alpha\|\nabla c^{+}(\tau)\|_{\dot  B^{\frac{N}{2}-1}_{2,1}}^h\big)
\int_0^t\langle t-\tau\rangle^{-(\frac N4+\frac s2)}\langle \tau\rangle^{-(\alpha+\frac N4)}\,d\tau\\
&\qquad+\big(\sup_{0\leq\tau\leq t}\langle \tau\rangle^\alpha\|(c^{+},c^{-})(\tau)\|_{\dot  B^{\frac{N}{2}}_{2,1}}^h\big)
\big(\sup_{0\leq\tau\leq t}\langle \tau\rangle^{\frac N4+\frac 12}\|\nabla c^{+}(\tau)\|_{\dot  B^{0}_{2,1}}^{\ell}\big)
\int_0^t\langle t-\tau\rangle^{-(\frac N4+\frac s2)}\langle \tau\rangle^{-(\alpha+\frac N4+\frac 12)}\,d\tau\\
&\qquad+\big(\sup_{0\leq\tau\leq t}\langle \tau\rangle^\alpha\|(c^{+},c^{-})(\tau)\|_{\dot  B^{\frac{N}{2}}_{2,1}}^h\big)
\big(\sup_{0\leq\tau\leq t}\langle \tau\rangle^{\alpha}\|\nabla c^{+}(\tau)\|_{\dot  B^{\frac{N}{2}-1}_{2,1}}^h\big)
\int_0^t\langle t-\tau\rangle^{-(\frac N4+\frac s2)}\langle \tau\rangle^{-2\alpha}\,d\tau\\
&\quad\lesssim D^{2}(t)\int_0^t\langle t-\tau\rangle^{-(\frac N4+\frac s2)}\langle \tau\rangle^{-\min({\frac N2+\frac 12,\alpha+\frac{N}{4},2\alpha})}\,d\tau\\
&\quad\lesssim \langle t\rangle^{-(\frac N4+\frac s2)} D^{2}(t),
\end{split}
\end{equation}
where $g_{+}$ stands for some smooth function vanishing at $0$.\\
Similar to  \eqref{s21}, we have
\begin{equation}
\begin{split}\label{s22}
&\int_0^t\langle t-\tau\rangle^{-(\frac N4+\frac s2)}\|\tilde{g}_{+}(c^{+},c^{-})\partial_{i}c^{-}\|_{L^1}\,d\tau\\
&\quad\lesssim D^{2}(t)\int_0^t\langle t-\tau\rangle^{-(\frac N4+\frac s2)}\langle \tau\rangle^{-\min({\frac N2+\frac 12,\alpha+\frac{N}{4},2\alpha})}\,d\tau\\
&\quad\lesssim \langle t\rangle^{-(\frac N4+\frac s2)} D^{2}(t).
\end{split}
\end{equation}
To bound the term with $(u^{+}\cdot\nabla)u_{i}^{+}$, we employ the following decomposition:
$$(u^{+}\cdot\nabla)u_{i}^{+}=(u^{+}\cdot\nabla)(u_{i}^{+})^{\ell}
+(u^{+}\cdot\nabla)(u_{i}^{+})^{h}.
$$
For the term  $(u^{+}\cdot\nabla)(u_{i}^{+})^{\ell}$, we have
\begin{equation*}
\begin{split}
&\int_0^t\langle t-\tau\rangle^{-(\frac N4+\frac s2)}
\|(u^{+}\cdot\nabla)(u_{i}^{+})^{\ell}\|_{L^1}\,d\tau\\
&\quad\lesssim\int_0^t\langle t-\tau\rangle^{-(\frac N4+\frac s2)}\big(\|u^{+}\|_{\dot  B^{0}_{2,1}}^{\ell}+\|u^{+}\|_{\dot  B^{0}_{2,1}}^h\big)\|\nabla u^{+}\|_{\dot  B^{0}_{2,1}}^{\ell}\,d\tau\\
&\quad\lesssim\big(\sup_{0\leq\tau\leq t}\langle \tau\rangle^{\frac N4}\|u^{+}(\tau)\|_{\dot  B^{0}_{2,1}}^{\ell}\big)
\big(\sup_{0\leq\tau\leq t}\langle \tau\rangle^{\frac{N}{4}+\frac 12}\| u^{+}(\tau)\|_{\dot  B^{1}_{2,1}}^\ell\big)
\\&\qquad\times\int_0^t\langle t-\tau\rangle^{-(\frac N4+\frac s2)}\langle \tau\rangle^{-(\frac{N}{2}+\frac 12)}\,d\tau\\
&\qquad+\big(\sup_{0\leq\tau\leq t}\langle \tau\rangle^\alpha\|u^{+}(\tau)\|_{\dot  B^{\frac{N}{2}-1}_{2,1}}^h\big)
\big(\sup_{0\leq\tau\leq t}\langle \tau\rangle^{\frac N4+\frac 12}\| u^{+}(\tau)\|_{\dot  B^{1}_{2,1}}^{\ell}\big)
\\&\qquad\times\int_0^t\langle t-\tau\rangle^{-(\frac N4+\frac s2)}\langle \tau\rangle^{-(\alpha+\frac N4+\frac 12)}\,d\tau\\
&\quad\lesssim D^{2}(t)\int_0^t\langle t-\tau\rangle^{-(\frac N4+\frac s2)}\langle \tau\rangle^{-\min({\frac N2+\frac 12,\alpha+\frac{N}{4}+\frac 12})}\,d\tau\\
&\quad\lesssim \langle t\rangle^{-(\frac N4+\frac s2)} D^{2}(t).
\end{split}
\end{equation*}
Regarding the term with $(u^{+}\cdot\nabla)(u_{i}^{+})^{h},$ we have for all $t\geq2$ that
\begin{equation*}
\begin{split}
&\int_0^t\langle t-\tau\rangle^{-(\frac N4+\frac s2)}
\|(u^{+}\cdot\nabla)(u_{i}^{+})^{h}(\tau)\|_{L^1}\,d\tau\\
&\quad\lesssim\int_0^t\langle t-\tau\rangle^{-(\frac N4+\frac s2)} \|u^{+}(\tau)\|_{L^2}\|\nabla u^{+}(\tau)\|^h_{L^2}\,d\tau\\
&\quad\lesssim\int_0^t\langle t-\tau\rangle^{-(\frac N4+\frac s2)} \|u^{+}(\tau)\|_{\dot  B^{0}_{2,1}}\|\nabla u^{+}(\tau)\|^h_{\dot  B^{0}_{2,1}}\,d\tau\\
&\quad\lesssim\int_0^1\langle t-\tau\rangle^{-(\frac N4+\frac s2)} \|u^{+}(\tau)\|_{\dot  B^{0}_{2,1}}\|\nabla u(\tau)\|^h_{\dot  B^{0}_{2,1}}\,d\tau\\
&\qquad+\int_1^t\langle t-\tau\rangle^{-(\frac N4+\frac s2)} \|u^{+}(\tau)\|_{\dot  B^{0}_{2,1}}\|\nabla u^{+}(\tau)\|^h_{\dot  B^{0}_{2,1}}\,d\tau\\
&\quad\eqdefa K_{1}+K_{2}.
\end{split}
\end{equation*}
From the definitions of $X(t)$ and $D(t)$, we  get
\begin{align*}
K_{1}\lesssim\langle t\rangle^{-(\frac N4+\frac s2)}
\sup _{0\leq\tau\leq1}\|u^{+}(\tau)\|_{\dot  B^{0}_{2,1}}\int_0^1\|\nabla u^{+}(\tau)\|^h_{\dot  B^{0}_{2,1}}\,d\tau\lesssim\langle t\rangle^{-(\frac N4+\frac s2)} D(1)X(1),
\end{align*}
and,  using the fact that $\langle \tau\rangle\approx\tau$ when $\tau\geq1$,
\begin{align*}
K_{2}&\lesssim\int_1^t\langle t-\tau\rangle^{-(\frac N4+\frac s2)}
\big(\|u^{+}(\tau)\|_{\dot  B^{0}_{2,1}}^{\ell}+\|u^{+}(\tau)\|_{\dot  B^{0}_{2,1}}^{h}\big)\|\nabla u^{+}(\tau)\|^h_{\dot  B^{0}_{2,1}}\,d\tau\\
&\lesssim\big(\sup_{1\leq\tau\leq t}\langle \tau\rangle^{\frac N4}\|u^{+}(\tau)\|_{\dot  B^{0}_{2,1}}^{\ell}\big)
\big(\sup_{1\leq\tau\leq t}\| \tau\nabla u^{+}(\tau)\|^{h}_{\dot  B^{\frac{N}{2}}_{2,1}}\big)
\int_1^t\langle t-\tau\rangle^{-(\frac N4+\frac s2)}\langle \tau\rangle^{-(\frac N4+1)}\,d\tau\\
&\quad+\big(\sup_{1\leq\tau\leq t}\langle \tau\rangle^{\alpha}\|u^{+}(\tau)\|_{\dot  B^{\frac{N}{2}-1}_{2,1}}^h\big)
\big(\sup_{1\leq\tau\leq t}\| \tau\nabla u^{+}(\tau)\|^{h}_{\dot  B^{\frac{N}{2}}_{2,1}}\big)
\int_1^t\langle t-\tau\rangle^{-(\frac N4+\frac s2)}\langle \tau\rangle^{-(\alpha+1)}\,d\tau\\
&\lesssim D^{2}(t)\int_1^t\langle t-\tau\rangle^{-(\frac N4+\frac s2)}\langle \tau\rangle^{-\min({\alpha+1,\frac{N}{4}+1})}\,d\tau\\
&\lesssim \langle t\rangle^{-(\frac N4+\frac s2)} D^{2}(t).
\end{align*}
Therefore, for $t\geq2,$ we deduce that
\begin{align}\label{S23h2}
&\int_0^t\langle t-\tau\rangle^{-(\frac N4+\frac s2)}
\|(u^{+}\cdot\nabla)(u_{i}^{+})^{h}(\tau)\|_{L^1}\,d\tau\nonumber\\
&\quad\lesssim\langle t\rangle^{-(\frac N4+\frac s2)} \big(D^{2}(t)+X^{2}(t)\big).
\end{align}
The case $t\leq2$ is obvious as $\langle t\rangle\approx1$ and
$\langle t-\tau\rangle\approx1$ for $0\leq\tau\leq t\leq 2$, and
\begin{equation}
\begin{split}
\label{s23h1}
&\int_0^t\|(u^{+}\cdot\nabla)(u_{i}^{+})^{h}\|_{L^1}\,d\tau\\
&\quad\lesssim \|u^{+}\|_{L^\infty_t(L^2)}\|\nabla u^{+}\|^h_{L_t^1(L^2)}\\
&\quad\lesssim \|u^{+}\|_{L^\infty_t(\dot  B^{0}_{2,1})}\|\nabla u^{+}\|^h_{L_t^1(\dot  B^{0}_{2,1})}\\
&\quad\lesssim \|u^{+}\|_{L^\infty_t(\dot  B^{0}_{2,1})}\|u^{+}\|^h_{L_t^1(\dot  B^{\frac{N}{2}+1}_{2,1})}\\
&\quad\lesssim X(t)D(t).
\end{split}
\end{equation}
To deal with the term  $\mu^{+}h_{+}(c^{+},c^{-})\partial_{j}c^{+}\partial_{j}u^{+}_{i}$, we take the following decomposition:
$$\mu^{+}h_{+}(c^{+},c^{-})\partial_{j}c^{+}\partial_{j}u^{+}_{i}
=\mu^{+}h_{+}(c^{+},c^{-})\partial_{j}c^{+}\partial_{j}(u^{+}_{i})^{\ell}
+\mu^{+}h_{+}(c^{+},c^{-})\partial_{j}c^{+}\partial_{j}(u^{+}_{i})^{h}.
$$
For the term  $\mu^{+}h_{+}(c^{+},c^{-})\partial_{j}c^{+}\partial_{j}(u^{+}_{i})^{\ell}$, we have
\begin{equation*}
\begin{split}
\label{s24l}
&\int_0^t\langle t-\tau\rangle^{-(\frac N4+\frac s2)}
\|\mu^{+}h_{+}(c^{+},c^{-})\partial_{j}c^{+}\partial_{j}(u^{+}_{i})^{\ell}\|_{L^1}\,d\tau\\
&\quad\lesssim\int_0^t\langle t-\tau\rangle^{-(\frac N4+\frac s2)}
\|h_{+}(c^{+},c^{-})\|_{L^2}
\|\nabla c^{+}\nabla (u^{+})^{\ell}\|_{L^2}\,d\tau\\
&\quad\lesssim\int_0^t\langle t-\tau\rangle^{-(\frac N4+\frac s2)}
\|h_{+}(c^{+},c^{-})\|_{\dot  B^{0}_{2,1}}
\|\nabla c^{+}\nabla (u^{+})^{\ell}\|_{\dot  B^{0}_{2,1}}\,d\tau\\
&\quad\lesssim\int_0^t\langle t-\tau\rangle^{-(\frac N4+\frac s2)}
\big(1+\|(c^{+},c^{-})\|_{\dot  B^{0}_{2,1}}\big)
\|\nabla c^{+}\|_{\dot  B^{0}_{2,1}}\|\nabla(u^{+})^{\ell}\|_{\dot  B^{\frac N2}_{2,1}}\,d\tau\\
&\quad\lesssim\int_0^t\langle t-\tau\rangle^{-(\frac N4+\frac s2)}
\|\nabla c^{+}\|_{\dot  B^{0}_{2,1}}\|\nabla(u^{+})^{\ell}\|_{\dot  B^{\frac N2}_{2,1}}\,d\tau\\
&\qquad+\int_0^t\langle t-\tau\rangle^{-(\frac N4+\frac s2)}
\|(c^{+},c^{-})\|_{\dot  B^{0}_{2,1}}
\|\nabla c^{+}\|_{\dot  B^{0}_{2,1}}\|\nabla(u^{+})^{\ell}\|_{\dot  B^{\frac N2}_{2,1}}\,d\tau\\
&\quad\eqdefa L_{1}+L_{2}.
\end{split}
\end{equation*}
We bound the two terms $L_{1}$ and $L_{2}$ as follows respectively,
\begin{equation*}
\begin{split}
&L_{1}\lesssim\int_0^t\langle t-\tau\rangle^{-(\frac N4+\frac s2)}
\big(\|\nabla c^{+}\|_{\dot  B^{0}_{2,1}}^{\ell}+\|\nabla c^{+}\|_{\dot  B^{0}_{2,1}}^h\big)
\|\nabla u^{+}\|_{\dot  B^{0}_{2,1}}^{\ell}\,d\tau\\
&\quad\lesssim
\big(\sup_{0\leq\tau\leq t}\langle \tau\rangle^{\frac N4+\frac 12}\|\nabla c^{+}(\tau)\|_{\dot  B^{0}_{2,1}}^{\ell}\big)
\big(\sup_{0\leq\tau\leq t}\langle \tau\rangle^{\frac{N}{2}+\frac 12}\| \nabla u^{+}(\tau)\|_{\dot  B^{\frac N2}_{2,1}}^\ell\big)
\\&\qquad\times\int_0^t\langle t-\tau\rangle^{-(\frac N4+\frac s2)}\langle \tau\rangle^{-(\frac{3N}{4}+1)}\,d\tau\\
&\qquad+
\big(\sup_{0\leq\tau\leq t}\langle \tau\rangle^{\alpha}\|\nabla c^{+}(\tau)\|_{\dot  B^{\frac N2 -1}_{2,1}}^{h}\big)
\big(\sup_{0\leq\tau\leq t}\langle \tau\rangle^{\frac{N}{2}+\frac 12}\| \nabla u^{+}(\tau)\|_{\dot  B^{\frac N2}_{2,1}}^\ell\big)
\\&\qquad\times\int_0^t\langle t-\tau\rangle^{-(\frac N4+\frac s2)}\langle\tau\rangle^{-(\frac{N}{2}+\frac 12+\alpha)}\,d\tau\\
&\quad\lesssim D^{2}(t)\int_0^t\langle t-\tau\rangle^{-(\frac N4+\frac s2)}\langle \tau\rangle^{-\min(\frac{3N}{4}+1,\frac{N}{2}+\frac 12+\alpha)}\,d\tau\\
&\quad\lesssim \langle t\rangle^{-(\frac N4+\frac s2)} D^{2}(t),
\end{split}
\end{equation*}
and
\begin{equation*}
\begin{split}
&L_{2}\lesssim\int_0^t\langle t-\tau\rangle^{-(\frac N4+\frac s2)}
\big(\|(c^{+},c^{-})\|_{\dot  B^{0}_{2,1}}^{\ell}+\|(c^{+},c^{-})\|_{\dot  B^{0}_{2,1}}^h\big)
\|\nabla c^{+}\|_{\dot  B^{0}_{2,1}}
\|\nabla u^{+}\|_{\dot  B^{0}_{2,1}}^{\ell}\,d\tau\\
&\quad\lesssim\big(\sup_{0\leq\tau\leq t}\langle \tau\rangle^{\frac N4}\|(c^{+},c^{-})(\tau)\|_{\dot  B^{0}_{2,1}}^{\ell}\big)
\int_0^t\langle t-\tau\rangle^{-(\frac N4+\frac s2)}\langle \tau\rangle^{-\frac N4}
\|\nabla c^{+}\|_{\dot  B^{0}_{2,1}}
\|\nabla(u^{+})^{\ell}\|_{\dot  B^{\frac N2}_{2,1}}\,d\tau\\
&\qquad+\big(\sup_{0\leq\tau\leq t}\langle \tau\rangle^{\alpha}\|(c^{+},c^{-})(\tau)\|_{\dot  B^{\frac{N}{2}}_{2,1}}^{h}\big)
\int_0^t\langle t-\tau\rangle^{-(\frac N4+\frac s2)}\langle\tau\rangle^{-\alpha}
\|\nabla c^{+}\|_{\dot  B^{0}_{2,1}}
\|\nabla(u^{+})^{\ell}\|_{\dot  B^{\frac N2}_{2,1}}\,d\tau\\
&\quad\lesssim D^{3}(t)\int_0^t\langle t-\tau\rangle^{-(\frac N4+\frac s2)}\langle \tau\rangle^{-\min(N+1,\frac{3N}{4}+\frac 12+\alpha,2\alpha+\frac N2+\frac 12)}\,d\tau\\
&\quad\lesssim \langle t\rangle^{-(\frac N4+\frac s2)} D^{3}(t).
\end{split}
\end{equation*}
Thus,
\begin{equation}
\begin{split}
\label{s24l}
\int_0^t\langle t-\tau\rangle^{-(\frac N4+\frac s2)}
\|\mu^{+}h_{+}(c^{+},c^{-})\partial_{j}c^{+}\partial_{j}(u^{+}_{i})^{\ell}\|_{L^1}\,d\tau\lesssim \langle t\rangle^{-(\frac N4+\frac s2)}\big( D^{2}(t)+D^{3}(t)\big).
\end{split}
\end{equation}
Regarding the term  $\mu^{+}h_{+}(c^{+},c^{-})\partial_{j}c^{+}\partial_{j}(u^{+}_{i})^{h},$ we also get,
\begin{equation*}
\begin{split}
&\int_0^t\langle t-\tau\rangle^{-(\frac N4+\frac s2)}
\|\mu^{+}h_{+}(c^{+},c^{-})\partial_{j}c^{+}\partial_{j}(u^{+}_{i})^{h}
(\tau)\|_{L^1}\,d\tau\\
&\quad\lesssim\int_0^t\langle t-\tau\rangle^{-(\frac N4+\frac s2)} \|h_{+}(c^{+},c^{-})(\tau)\|_{L^2}
\|\partial_{j}c^{+}\partial_{j}(u^{+}_{i})^{h}(\tau)\|_{L^2}\,d\tau\\
&\quad\lesssim\int_0^t\langle t-\tau\rangle^{-(\frac N4+\frac s2)} \big(1+\|(c^{+},c^{-})(\tau)\|_{\dot  B^{0}_{2,1}}\big)
\|\nabla c^{+}(\tau)\|_{\dot  B^{0}_{2,1}}
\|\nabla u^{+}(\tau)\|^h_{\dot  B^{\frac N2}_{2,1}}\,d\tau\\
&\quad\lesssim\int_0^t\langle t-\tau\rangle^{-(\frac N4+\frac s2)}
\|\nabla c^{+}(\tau)\|_{\dot  B^{0}_{2,1}}
\|\nabla u^{+}(\tau)\|^h_{\dot  B^{\frac N2}_{2,1}}\,d\tau\\
&\qquad+
\int_0^t\langle t-\tau\rangle^{-(\frac N4+\frac s2)} \|(c^{+},c^{-})(\tau)\|_{\dot  B^{0}_{2,1}}
\|\nabla c^{+}(\tau)\|_{\dot  B^{0}_{2,1}}
\|\nabla u^{+}(\tau)\|^h_{\dot  B^{\frac N2}_{2,1}}\,d\tau\\
&\quad\eqdefa M_{1}+M_{2}.
\end{split}
\end{equation*}
We deal with the two terms $M_{1}$ and $M_{2}$ in the following,
if $t\geq2,$
\begin{equation*}
\begin{split}
&M_{1}\lesssim\int_0^1\langle t-\tau\rangle^{-(\frac N4+\frac s2)}
\|\nabla c^{+}(\tau)\|_{\dot  B^{0}_{2,1}}
\|\nabla u^{+}(\tau)\|^h_{\dot  B^{\frac N2}_{2,1}}\,d\tau\\
&\qquad+\int_1^t\langle t-\tau\rangle^{-(\frac N4+\frac s2)}
\|\nabla c^{+}(\tau)\|_{\dot  B^{0}_{2,1}}
\|\nabla u^{+}(\tau)\|^h_{\dot  B^{\frac N2}_{2,1}}\,d\tau\\
&\quad\eqdefa M_{11}+M_{12}.
\end{split}
\end{equation*}
Using the definitions of $X(t)$ and $D(t)$, we  obtain
\begin{align*}
M_{11}\lesssim\langle t\rangle^{-(\frac N4+\frac s2)}
\sup _{0\leq\tau\leq1}\|\nabla c^{+}(\tau)\|_{\dot  B^{0}_{2,1}}
\int_0^1\|\nabla u^{+}(\tau)\|^h_{\dot B^{\frac{N}{2}}_{2,1}}\,d\tau\lesssim\langle t\rangle^{-(\frac N4+\frac s2)} D(1)X(1),
\end{align*}
and,  employing the fact that $\langle \tau\rangle\approx\tau$ when $\tau\geq1$, we have
\begin{align*}
M_{12}&\lesssim
\big(\sup_{1\leq\tau\leq t}\langle \tau\rangle^{\frac N4+\frac 12}
\|\nabla c^{+}(\tau)\|_{\dot  B^{0}_{2,1}}^{\ell}\big)
\big(\sup_{1\leq\tau\leq t}\| \tau\nabla u^{+}(\tau)\|^{h}_{\dot  B^{\frac{N}{2}}_{2,1}}\big)
\\&\qquad\times\int_1^t\langle t-\tau\rangle^{-(\frac N4+\frac s2)}\langle \tau\rangle^{-(\frac N4+\frac 32)}\,d\tau\\
&\qquad+
\big(\sup_{1\leq\tau\leq t}\langle \tau\rangle^{\alpha}
\|\nabla c^{+}(\tau)\|_{\dot  B^{\frac{N}{2}-1}_{2,1}}^{h}\big)
\big(\sup_{1\leq\tau\leq t}\| \tau\nabla u^{+}(\tau)\|^{h}_{\dot  B^{\frac{N}{2}}_{2,1}}\big)
\\&\qquad\times\int_1^t\langle t-\tau\rangle^{-(\frac N4+\frac s2)}\langle \tau\rangle^{-(\alpha+1)}\,d\tau\\
&\lesssim \langle t\rangle^{-(\frac N4+\frac s2)} D^{2}(t).
\end{align*}
For the term  $M_{2}$, we have
\begin{equation*}
\begin{split}
M_{2}&\lesssim\int_0^1\langle t-\tau\rangle^{-(\frac N4+\frac s2)} \|(c^{+},c^{-})(\tau)\|_{\dot  B^{0}_{2,1}}
\|\nabla c^{+}(\tau)\|_{\dot  B^{0}_{2,1}}
\|\nabla u^{+}(\tau)\|^h_{\dot  B^{\frac N2}_{2,1}}\,d\tau\\
&\qquad+\int_1^t\langle t-\tau\rangle^{-(\frac N4+\frac s2)} \|(c^{+},c^{-})(\tau)\|_{\dot  B^{0}_{2,1}}
\|\nabla c^{+}(\tau)\|_{\dot  B^{0}_{2,1}}
\|\nabla u^{+}(\tau)\|^h_{\dot  B^{\frac N2}_{2,1}}\,d\tau\\
&\eqdefa M_{21}+M_{22}.
\end{split}
\end{equation*}
Remembering the definitions of $X(t)$ and $D(t)$, we  obtain
\begin{align*}
M_{21}&\lesssim\langle t\rangle^{-(\frac N4+\frac s2)}
\sup _{0\leq\tau\leq1}\|(c^{+},c^{-})(\tau)\|_{\dot  B^{0}_{2,1}}
\sup _{0\leq\tau\leq1}\|\nabla c^{+}(\tau)\|_{\dot  B^{0}_{2,1}}
\int_0^1\|\nabla u^{+}(\tau)\|^h_{\dot B^{\frac{N}{2}}_{2,1}}\,d\tau\\
&\lesssim\langle t\rangle^{-(\frac N4+\frac s2)} D^{2}(1)X(1),
\end{align*}
and
\begin{align*}
M_{22}&\lesssim\big(\sup_{1\leq\tau\leq t}\langle \tau\rangle^{\frac N4}
\|(c^{+},c^{-})(\tau)\|_{\dot  B^{0}_{2,1}}^{\ell}\big)
\int_1^t\langle t-\tau\rangle^{-(\frac N4+\frac s2)}
\langle \tau\rangle^{-\frac N4}
\|\nabla c^{+}(\tau)\|_{\dot  B^{0}_{2,1}}
\|\nabla u^{+}(\tau)\|^h_{\dot  B^{\frac N2}_{2,1}}\,d\tau\\
&\quad+\big(\sup_{1\leq\tau\leq t}\langle \tau\rangle^{\alpha}\|(c^{+},c^{-})(\tau)\|_{\dot  B^{\frac{N}{2}}_{2,1}}^h\big)
\int_1^t\langle t-\tau\rangle^{-(\frac N4+\frac s2)}
\langle\tau\rangle^{-\alpha}
\|\nabla c^{+}(\tau)\|_{\dot  B^{0}_{2,1}}
\|\nabla u^{+}(\tau)\|^h_{\dot  B^{\frac N2}_{2,1}}\,d\tau\\\
&\lesssim D^{3}(t)\int_1^t\langle t-\tau\rangle^{-(\frac N4+\frac s2)}\langle \tau\rangle^{-\min({2\alpha+1,\alpha+\frac{N}{4}+1,\frac{N}{2}+\frac{3}{2}})}\,d\tau\\
&\lesssim \langle t\rangle^{-(\frac N4+\frac s2)} D^{3}(t).
\end{align*}
Therefore, for $t\geq2,$  we obatin
\begin{align}
\label{s24h2}
&\int_0^t\langle t-\tau\rangle^{-(\frac N4+\frac s2)}
\|\mu^{+}h_{+}(c^{+},c^{-})\partial_{j}c^{+}\partial_{j}(u^{+}_{i})^{h}(\tau)\|_{L^1}\,d\tau\nonumber\\
&\quad\lesssim\langle t\rangle^{-(\frac N4+\frac s2)} \big(X^{2}(t)+D^{2}(t)+D^{3}(t)+D^{4}(t)\big).
\end{align}
The case $t\leq2$ is obvious as $\langle t\rangle\approx1$ and
$\langle t-\tau\rangle\approx1$ for $0\leq\tau\leq t\leq 2$, and
\begin{equation}
\begin{split}
\label{s24h1}
&\int_0^t\|\mu^{+}h_{+}(c^{+},c^{-})\partial_{j}c^{+}\partial_{j}(u^{+}_{i})^{h}
(\tau)\|_{L^1}\,d\tau\\
&\quad\lesssim \|h_{+}(c^{+},c^{-})\|_{L^\infty_t(L^2)}
\|\nabla c^{+}\nabla (u^{+})^h\|_{L_t^1(L^2)}\\
&\quad\lesssim\big( \|h_{+}(c^{+},c^{-})-\frac{(\mathcal{C}^{2}\alpha^{-})(1,1)}
{s_{-}^{2}(1,1)}\|_{\dot  B^{0}_{2,1}}+1\big)
\|\nabla c^{+}\|_{L^\infty_t(\dot  B^{0}_{2,1})}
\|\nabla u^{+}\|^h_{L_t^1(\dot  B^{\frac{N}{2}}_{2,1})}\\
&\quad\lesssim \big(1+\|(c^{+},c^{-})\|_{L^\infty_t(\dot  B^{0}_{2,1})}\big)
\|\nabla c^{+}\|_{L^\infty_t(\dot  B^{0}_{2,1})}
\|u^{+}\|^h_{L_t^1(\dot  B^{\frac{N}{2}+1}_{2,1})}\\
&\quad\lesssim \big(1+D(t)\big)D(t)X(t)\\
&\quad\lesssim X^{2}(t)+D^{2}(t)+D^{4}(t).
\end{split}
\end{equation}
From \eqref{s24l}-\eqref{s24h1}, we finally conclude that
\begin{equation}
\begin{split}
\label{s24}
\int_0^t\|\mu^{+}h_{+}(c^{+},c^{-})\partial_{j}c^{+}\partial_{j}u^{+}_{i}
(\tau)\|_{L^1}\,d\tau\lesssim \langle t\rangle^{-(\frac N4+\frac s2)} \big(X^{2}(t)+D^{2}(t)+D^{3}(t)+D^{4}(t)\big).
\end{split}
\end{equation}
Similarly,  we also obtain  the 	corresponding estimates of other terms
$\mu^{+}k_{+}(c^{+},c^{-})\partial_{j}c^{-}\partial_{j}u^{+}_{i},$\\
$\mu^{+}h_{+}(c^{+},c^{-})\partial_{j}c^{+}\partial_{i}u^{+}_{j},$
$\mu^{+}k_{+}(c^{+},c^{-})\partial_{j}c^{-}\partial_{i}u^{+}_{j},
\lambda^{+}h_{+}(c^{+},c^{-})\partial_{i}c^{+}\partial_{j}u^{+}_{j}~\text{and}~
\lambda^{+}k_{+}(c^{+},c^{-})\partial_{i}c^{-}\partial_{j}u^{+}_{j}$. Here, we omit the details.

From the low-high frequency  decomposition for $\mu^{+}l_{+}(c^{+},c^{-})\partial_{j}^{2}u_{i}^{+},$  we have
\begin{equation*}
\mu^{+}l_{+}(c^{+},c^{-})\partial_{j}^{2}u_{i}^{+}
=\mu^{+}l_{+}(c^{+},c^{-})\partial_{j}^{2}(u_{i}^{+})^\ell
+\mu^{+}l_{+}(c^{+},c^{-})\partial_{j}^{2}(u_{i}^{+})^h,
\end{equation*}
where $l_{+}$ stands for some smooth function vanishing at $0$.
Thus,
 \begin{equation}
 \begin{split}
 \label{s10l}
 &\int_0^t\langle t-\tau\rangle^{-(\frac N4+\frac s2)}\|\mu^{+}l_{+}(c^{+},c^{-})\partial_{j}^{2}(u_{i}^{+})^\ell\|_{L^1}\,d\tau\\
 &\quad\lesssim\int_0^t\langle t-\tau\rangle^{-(\frac N4+\frac s2)}\|l_{+}(c^{+},c^{-})\|_{\dot  B^{0}_{2,1}}
 \|\nabla^{2}u^{+}\|_{\dot  B^{0}_{2,1}}^\ell\,d\tau\\
 &\quad\lesssim\int_0^t\langle t-\tau\rangle^{-(\frac N4+\frac s2)}\|(c^{+},c^{-})\|_{\dot  B^{0}_{2,1}}
 \|\nabla^{2}u^{+}\|_{\dot  B^{0}_{2,1}}^\ell\,d\tau\\
 &\quad\lesssim \big(\sup_{\tau\in[0,t]} \langle\tau\rangle^{\frac N4}\|(c^{+},c^{-})(\tau)\|_{\dot  B^{0}_{2,1}}^{\ell}\big)
 \big(\sup_{\tau\in[0,t]} \langle\tau\rangle^{\frac N4+1}
 \|\nabla^{2}u^{+}\|_{\dot  B^{0}_{2,1}}^\ell)\big)
\\&\qquad\times \int_0^t \langle t-\tau\rangle^{-(\frac N4+\frac s2)}
 \langle\tau\rangle^{-(\frac N2+1)}\,d\tau\\
 &\qquad+ \big(\sup_{0\leq\tau\leq t}\langle \tau\rangle^{\alpha}\|(c^{+},c^{-})(\tau)\|_{\dot  B^{\frac{N}{2}}_{2,1}}^h\big)
 \big(\sup_{\tau\in[0,t]} \langle\tau\rangle^{\frac N4+1}
 \|\nabla^{2}u^{+}\|_{\dot  B^{0}_{2,1}}^\ell)\big)
 \\&\qquad\times\int_0^t \langle t-\tau\rangle^{-(\frac N4+\frac s2)}\langle\tau\rangle^{-(\frac N4+\alpha+1)}\,d\tau\\
 &\quad\lesssim D^{2}(t)\int_0^t\langle t-\tau\rangle^{-(\frac N4+\frac s2)}\langle \tau\rangle^{-\min(\frac N2+1,\frac N4+\alpha+1)}\,d\tau\\
 &\quad\lesssim \langle t\rangle^{-(\frac N4+\frac s2)} D^{2}(t).
 \end{split}
 \end{equation}
To handle the term $\mu^{+}l_{+}(c^{+},c^{-})\partial_{j}^{2}(u_{i}^{+})^h$, we consider the cases $t\geq2$ and $t\leq2$ respectively. When $t\geq2$, then we have
 \begin{align*}
 &\int_0^t\langle t-\tau\rangle^{-(\frac N4+\frac s2)}\|\mu^{+}l_{+}(c^{+},c^{-})\partial_{j}^{2}(u_{i}^{+})^h\|_{L^1}\,d\tau
 \\&\quad\lesssim\int_0^t\langle t-\tau\rangle^{-(\frac N4+\frac s2)}\|l_{+}(c^{+},c^{-})\|_{\dot  B^{0}_{2,1}}
 \|\nabla^{2} u^{+}\|^h_{\dot  B^{0}_{2,1}}\,d\tau\\
 &\quad\lesssim\int_0^1\langle t-\tau\rangle^{-(\frac N4+\frac s2)}\|(c^{+},c^{-})\|_{\dot  B^{0}_{2,1}}
 \|\nabla^{2} u^{+}\|^h_{\dot  B^{0}_{2,1}}\,d\tau
\\&\qquad +\int_1^t\langle t-\tau\rangle^{-(\frac N4+\frac s2)}\|(c^{+},c^{-})\|_{\dot  B^{0}_{2,1}}\|\nabla^{2} u^{+}\|^h_{\dot  B^{0}_{2,1}}\,d\tau\\
 &\quad\eqdefa N_{1}+N_{2}.
 \end{align*}
From the definitions of $X(t)$ and $D(t)$, we  obtain
\begin{align*}
 N_{1}&=\int_0^1\langle t-\tau\rangle^{-(\frac N4+\frac s2)}
 \|(c^{+},c^{-})\|_{\dot  B^{0}_{2,1}}
 \|\nabla^{2} u^{+}\|^h_{\dot  B^{0}_{2,1}}\,d\tau\,d\tau\\
 &\lesssim\langle t\rangle^{-(\frac N4+\frac s2)} \big(\sup_{\tau\in[0,1]}\|(c^{+},c^{-})\|_{\dot  B^{0}_{2,1}}\big)
\int_0^1 \|u^{+}\|^h_{\dot  B^{\frac{N}{2}+1}_{2,1}}\,d\tau\\
&\lesssim\langle t\rangle^{-(\frac N4+\frac s2)} D(1)X(1),
\end{align*}
and, using the fact that $\langle \tau\rangle\approx\tau$ when $\tau\geq1$, we have
\begin{align*}
 N_{2}&=\int_1^t\langle t-\tau\rangle^{-(\frac N4+\frac s2)}
 \|(c^{+},c^{-})\|_{\dot  B^{0}_{2,1}}\|\nabla^{2} u^{+}\|^h_{\dot  B^{0}_{2,1}}\,d\tau\\
&\lesssim\int_1^t\langle t-\tau\rangle^{-(\frac N4+\frac s2)}
\big(\|(c^{+},c^{-})\|_{\dot  B^{0}_{2,1}}^{\ell}+\|(c^{+},c^{-})\|_{\dot  B^{0}_{2,1}}^{h}\big)\|\nabla^{2} u^{+}\|^h_{\dot  B^{0}_{2,1}}\,d\tau\\
&\lesssim\int_1^t\langle t-\tau\rangle^{-(\frac N4+\frac s2)}
\|(c^{+},c^{-})\|_{\dot  B^{0}_{2,1}}^{\ell}\|\nabla^{2} u^{+}\|^h_{\dot B^{0}_{2,1}}\,d\tau
+\int_1^t\langle t-\tau\rangle^{-(\frac N4+\frac s2)}
\|(c^{+},c^{-})\|_{\dot  B^{0}_{2,1}}^{h})\|\nabla^{2} u^{+}\|^h_{\dot  B^{0}_{2,1}}\,d\tau\\
&\lesssim\big(\sup_{0\leq\tau\leq t}\langle \tau\rangle^{\frac N4}
\|(c^{+},c^{-})(\tau)\|_{\dot  B^{0}_{2,1}}^{\ell}\big)
\big(\sup_{0\leq\tau\leq t} \|\tau\nabla u(\tau)\|_{\dot  B^{\frac{N}{2}}_{2,1}}^{h}\big)
\int_1^t\langle t-\tau\rangle^{-(\frac N4+\frac s2)}\langle \tau\rangle^{-(\frac N4+1)}\,d\tau\\
&\quad+\big(\sup_{0\leq\tau\leq t}\langle \tau\rangle^{\alpha}\|\nabla (c^{+},c^{-})(\tau)\|_{\dot  B^{\frac{N}{2}-1}_{2,1}}^h\big)
\big(\sup_{0\leq\tau\leq t} \|\tau\nabla u(\tau)\|_{\dot  B^{\frac{N}{2}}_{2,1}}^h\big)
\int_1^t\langle t-\tau\rangle^{-(\frac N4+\frac s2)}\langle \tau\rangle^{-(\alpha+1)}d\tau\\
&\lesssim D^{2}(t)
\int_1^t\langle t-\tau\rangle^{-(\frac N4+\frac s2)}\langle \tau\rangle^{-\min(\alpha+1,\frac N4+1)}d\tau\\
&\lesssim \langle t\rangle^{-(\frac N4+\frac s2)} D^{2}(t).
\end{align*}
Thus, for $t\geq2,$ we arrive at
\begin{equation}
\label{s10h2}
\begin{split}
\int_0^t\langle t-\tau\rangle^{-(\frac N4+\frac s2)}\|\mu^{+}l_{+}(c^{+},c^{-})\partial_{j}^{2}(u_{i}^{+})^h\|_{L^1}\,d\tau
\lesssim\langle t\rangle^{-(\frac N4+\frac s2)}(X^{2}(t)+D^{2}(t)).
\end{split}
\end{equation}
The case $t\leq2$ is obvious as $\langle t\rangle\approx1$ and $\langle t-\tau\rangle\approx1$ for $0\leq\tau\leq t\leq2,$
\begin{equation}
\label{s10h1}
\begin{split}
&\int_0^t\|\mu^{+}l_{+}(c^{+},c^{-})\partial_{j}^{2}(u_{i}^{+})^h\|_{L^1}\,d\tau\\
&\quad\lesssim\int_0^t\|l_{+}(c^{+},c^{-})\|_{\dot  B^{0}_{2,1}}
\|\nabla^{2}u^+\|_{\dot  B^{0}_{2,1}}^h\,d\tau\\
&\quad\lesssim\big(\sup_{\tau\in[0,1]}\|(c^{+},c^{-})(\tau)\|_{\dot  B^{0}_{2,1}}\big)
\int_0^1 \|u\|_{\dot  B^{\frac{N}{2}+1}_{2,1}}^h\,d\tau\\
&\quad\lesssim D(t)X(t).
\end{split}
\end{equation}
From \eqref{s10l}-\eqref{s10h1}, we get
\begin{equation}
\label{s10}
\begin{split}
\int_0^t\langle t-\tau\rangle^{-(\frac N4+\frac s2)}\|\mu^{+}l_{+}(c^{+},c^{-})\partial_{j}^{2}u_{i}^{+}\|_{L^1}\,d\tau
\lesssim \langle t\rangle^{-(\frac N4+\frac s2)} \big(D^{2}(t)+X^{2}(t)\big).
\end{split}
\end{equation}
Similarly,
\begin{equation}
\label{s1last}
\begin{split}
\int_0^t\langle t-\tau\rangle^{-(\frac N4+\frac s2)}\|(\mu^{+}+\lambda^{+})l_{+}(c^{+},c^{-})\partial_{i}\partial_{j}
u^{+}_{j}\|_{L^1}\,d\tau\lesssim \langle t\rangle^{-(\frac N4+\frac s2)} \big(D^{2}(t)+X^{2}(t)\big).
\end{split}
\end{equation}
Thus,
\begin{equation*}\label{s2}
\int_0^t\langle t-\tau\rangle^{-(\frac{N}{4}+\frac{s}{2})} \big\|H_{2}(\tau)\big\|_{\dot B^{-\frac{N}{2}}_{2,\infty}}^\ell d\tau
\lesssim\langle t\rangle^{-(\frac N4+\frac s2)}
\big(X^2(t)+D^2(t)+D^3(t)+D^4(t)\big).
\end{equation*}
The term  $H_{4}$ may be treated along the same lines,  and we have
\begin{equation*}\label{s4}
\int_0^t\langle t-\tau\rangle^{-(\frac{N}{4}+\frac{s}{2})} \big\|H_{4}(\tau)\big\|_{\dot B^{-\frac{N}{2}}_{2,\infty}}^\ell d\tau
\lesssim\langle t\rangle^{-(\frac N4+\frac s2)}
\big(X^2(t)+D^2(t)+D^3(t)+D^4(t)\big).
\end{equation*}
Thus, we complete the proof of \eqref{s1234low1}. Combining  with \eqref{U} and \eqref{s1234low1}, we conclude that for all $t\geq0$ and $s\in(\varepsilon-\frac{N}{2},2],$
\begin{equation}
\label{low}
\langle t\rangle^{\frac N4+\frac s2}
 \|(c^{+},\,u^{+},\,c^{-},\,u^{-}\|^{\ell}_{\dot  B^{s}_{2,1}}\lesssim D_{0}+X^{2}(t)+D^{2}(t)+D^{3}(t)+D^{4}(t).
\end{equation}
\subsubsection*{Step 2: High frequencies}
Now,  the starting point is Inequality \eqref{17} which implies that for $q\geq q_0$ and for some
$c_0=c(q_0)>0,$ we have
\begin{align*}
\frac12\frac d{dt}\alpha_q^2+c_0\alpha_q^2
&\leq\Bigl(\|(\dot{\Delta}_qH_{11},\,\dot{\Delta}_qH_{21},
\,\dot{\Delta}_qH_{31},\,\dot{\Delta}_qH_{41},\,\nabla\dot{\Delta}_qH_{11},
\,\nabla\dot{\Delta}_qH_{31},
)\|_{L^2}+\|R_q(u^{+},c^{+})\|_{L^2}\\
&\quad+\|R_q(u^{+},u^{+})\|_{L^2}+\|R_q(u^{-},c^{-})\|_{L^2}
+\|R_q(u^{-},u^{-})\|_{L^2}\\
&\quad+\|\wt R_q(u^{+},c^{+})\|_{L^2}+\|\wt R_q(u^{-},c^{-})\|_{L^2}+\|\nabla (u^{+},u^{-})\|_{L^\infty}\alpha_q\Bigr)\alpha_q,
\end{align*}
in which
$$R_q(u,b)\eqdefa [u^{\pm}\cdot\nabla,\ddq]b=u^{\pm}\cdot\nabla\ddq b-\ddq(u^{\pm}\cdot\nabla b) \quad\quad for \quad b\in\{c^{\pm},u^{\pm}\},$$
$$\wt R_q^i(u,b)\eqdefa [u^{\pm}\cdot\nabla,\partial_i\ddq]c^{\pm}=u^{\pm}\cdot\nabla\partial_i\ddq c^{\pm}-\partial_i\ddq(u^{\pm}\cdot\nabla c^{\pm}).$$
After time integration, we have
\begin{align*}
e^{c_0t}\alpha_q(t)
&\leq\alpha_q(0)+\int_0^te^{c_0\tau}\Bigl(\|(\dot{\Delta}_qH_{11},\,\dot{\Delta}_qH_{21},
\,\dot{\Delta}_qH_{31},\,\dot{\Delta}_qH_{41},\,\nabla\dot{\Delta}_qH_{11},
\,\nabla\dot{\Delta}_qH_{31},
)\|_{L^2}\\
&\quad+\|R_q(u^{+},c^{+})\|_{L^2}+\|R_q(u^{+},u^{+})\|_{L^2}+\|R_q(u^{-},c^{-})\|_{L^2}
+\|R_q(u^{-},u^{-})\|_{L^2}\\
&\quad+\|\wt R_q(u^{+},c^{+})\|_{L^2}+\|\wt R_q(u^{-},c^{-})\|_{L^2}+\|\nabla (u^{+},u^{-})\|_{L^\infty}\alpha_q\Bigr) d\tau.
\end{align*}
For  $q\geq q_0$, we have  $\alpha_q\approx\|(\nabla\ddq c^{+},\ddq u^{+},\nabla\ddq c^{-},\ddq u^{-})\|_{L^2}$. Then,
\begin{align*}
&\langle t\rangle^\alpha\|(\nabla\ddq c^{+},\ddq u^{+},\nabla\ddq c^{-},\ddq u^{-})(t)\|_{L^2}\lesssim \langle t\rangle^\alpha e^{-c_0t}\|(\nabla\ddq c^{+},\ddq u^{+},\nabla\ddq c^{-},\ddq u^{-})(0)\|_{L^2}\\
&\quad+\langle t\rangle^\alpha\int_0^te^{c_0(\tau-t)}\Bigl(\|(\dot{\Delta}_qH_{11},\,\dot{\Delta}_qH_{21},
\,\dot{\Delta}_qH_{31},\,\dot{\Delta}_qH_{41},\,\nabla\dot{\Delta}_qH_{11},
\,\nabla\dot{\Delta}_qH_{31},
)\|_{L^2}\\
&\quad+\|R_q(u^{+},c^{+})\|_{L^2}+\|R_q(u^{+},u^{+})\|_{L^2}+\|R_q(u^{-},c^{-})\|_{L^2}
+\|R_q(u^{-},u^{-})\|_{L^2}\\
&\quad+\|\wt R_q(u^{+},c^{+})\|_{L^2}+\|\wt R_q(u^{-},c^{-})\|_{L^2}+\|\nabla (u^{+},u^{-})\|_{L^\infty}\alpha_q\Bigr) d\tau,
\end{align*}
and thus, by  multiplying both sides by $2^{q(\frac N2-1)},$ taking the supremum on $[0,T],$
and  then summing up over $q\geq q_0,$
\begin{equation}
\begin{split}\label{S7}
&\|\langle t\rangle^\alpha(\nabla c^{+}, u^{+},\nabla c^{-}, u^{-})\|_{\wt L^\infty_T(\dot B^{\frac N2-1}_{2,1})}^{h}
\\&\qquad\lesssim\|(\nabla c^{+}_{0}, u^{+}_{0},\nabla c^{-}_{0}, u^{-}_{0})\|_{\dot B^{\frac N2-1}_{2,1}}^h
\\&\qquad+\sum_{q\geq q_0}\sup_{0\leq t\leq T}\Big(\langle t\rangle^\alpha\!\int_0^t\!e^{c_0(\tau-t)}2^{q(\frac N2-1)}S_q\,d\tau\Big)
\end{split}\end{equation}
with $S_q\eqdefa\sum_{i=1}^8 S_q^i$ and
$$\displaylines{
S_q^1\eqdefa \|(\dot{\Delta}_qH_{11},\,\dot{\Delta}_qH_{21},
\,\dot{\Delta}_qH_{31},\,\dot{\Delta}_qH_{41},\,\nabla\dot{\Delta}_qH_{11},
\,\nabla\dot{\Delta}_qH_{31},
)\|_{L^2},\quad
S_q^2\eqdefa\|R_q(u^{+},c^{+})\|_{L^2},\cr
S_q^3\eqdefa\|R_q(u^{+},u^{+})\|_{L^2},\quad
S_q^4\eqdefa\|R_q(u^{-},c^{-})\|_{L^2}, \quad
S_q^5\eqdefa\|R_q(u^{-},u^{-})\|_{L^2}, \quad
S_q^6\eqdefa\|\wt R_q(u^{+},c^{+})\|_{L^2}, \cr
S_q^7\eqdefa\|\wt R_q(u^{-},c^{-})\|_{L^2}, \quad
S_q^8\eqdefa\|\nabla (u^{+},u^{-})\|_{L^\infty}\|(\nabla\ddq c^{+},\ddq u^{+},\nabla\ddq c^{-},\ddq u^{-})\|_{L^2}.}$$
Bounding the sum, for $0\leq t\leq 2$, and taking advantage of Proposition \ref{Pro:1}, we end up with
\begin{equation}
\begin{split}
\label{s7t02}
\sum_{q\geq q_0}\sup_{0\leq t\leq 2}&\Big(\langle t\rangle^\alpha\!\int_0^t\!e^{c_0(\tau-t)}2^{q(\frac N2-1)}S_q(\tau)\,d\tau\Big)
\lesssim\int_0^2 \sum_{q\geq q_0}2^{q(\frac N2-1)}S_q(\tau)\,d\tau\\
&\quad\lesssim\int_0^2\Big(\|(H_{11},H_{21},H_{31},H_{41},\nabla H_{11},\nabla H_{31})\|_{\dot B^{\frac N2-1}_{2,1}}^{h}\\
&\qquad+\|\nabla (u^{+},u^{-})\|_{\dot B^{\frac N2}_{2,1}}\|(c^{+},u^{+},c^{-},u^{-},\nabla c^{+},\nabla c^{-})\|_{\dot B^{\frac N2-1}_{2,1}}\Big)d\tau\\
&\quad\lesssim\int_0^2\Big(\|(\nabla H_{11},H_{21},\nabla H_{31},H_{41})\|_{\dot B^{\frac N2-1}_{2,1}}^{h}\\
&\qquad+\|\nabla (u^{+},u^{-})\|_{\dot B^{\frac N2}_{2,1}}\|(c^{+},u^{+},c^{-},u^{-},\nabla c^{+},\nabla c^{-})\|_{\dot B^{\frac N2-1}_{2,1}}\Big)d\tau\\
&\quad\eqdefa Q_{1}+Q_{2}.
\end{split}
\end{equation}
From Propositions \ref{p26}-\ref{p27},   we bound  the terms  $Q_{1}$ and $Q_{2}$  as follows
\begin{align*}
\int_0^2\|\nabla H_{11}\|_{\dot B^{\frac N2-1}_{2,1}}^{h}d\tau
&\lesssim \int_0^2\| H_{11}\|_{\dot B^{\frac N2}_{2,1}}^{h}d\tau\\
&\lesssim \int_0^2\|c^{+}\div u^{+}\|_{\dot B^{\frac N2}_{2,1}}^{h}d\tau\\
&\lesssim \int_0^2\|c^{+}\|_{\dot B^{\frac N2}_{2,1}} \|\div u^{+}\|_{\dot B^{\frac N2}_{2,1}}d\tau\\
&\lesssim \| c^{+}\|_{L_{t}^{\infty}(\dot B^{\frac N2-1\frac N2}_{2,1})} \| u^{+}\|_{L_{t}^{1}(\dot B^{\frac N2+1}_{2,1})}\\
&\lesssim X^{2}(2).
\end{align*}
For the term  $\nabla H_{31}$, along the same lines, we have
\begin{align*}
\int_0^2\|\nabla H_{31}\|_{\dot B^{\frac N2-1}_{2,1}}^{h}d\tau
&\lesssim X^{2}(2).
\end{align*}
Combining interpolation inequality and H\"{o}lder's inequality, we deduce that
\begin{equation}
\begin{split}\label{s21step2}
&\int_0^t\|g_{+}(c^{+},c^{-})\partial_{i}c^{+}(\tau)\|_{\dot B^{\frac N2-1}_{2,1}}\,d\tau\\
&\quad\lesssim\|(c^{+},c^{-})\|_{L^{2}_{t}{(\dot B^{\frac N2}_{2,1})}}
\|\nabla c^{+}\|_{L^{2}_{t}{(\dot B^{\frac N2-1}_{2,1})}}\,\\
&\quad\lesssim\bigl(\|(c^{+},c^{-})\|_{L^{1}_{t}{(\dot B^{\frac N2+1 \frac N2}_{2,1})}}\bigr)^{\frac{1}{2}}
\big(\|(c^{+},c^{-})\|_{L^{\infty}_{t}{(\dot B^{\frac N2-1 \frac N2}_{2,1})}}\big)^{\frac{1}{2}}
\\&\qquad\times\bigl(\|c^{+}\|_{L^{1}_{t}{(\dot B^{\frac N2+1 \frac N2}_{2,1})}}\bigr)^{\frac{1}{2}}
\big(\|c^{+}\|_{L^{\infty}_{t}{(\dot B^{\frac N2-1 \frac N2}_{2,1})}}\big)^{\frac{1}{2}}\\
&\quad\lesssim X^{2}(2),
\end{split}
\end{equation}
where $g_{+}$ stands for some smooth function vanishing at $0$.

Similarly,
\begin{equation}
\begin{split}\label{s22step2}
&\int_0^t\|\tilde{g}_{+}(c^{+},c^{-})\partial_{i}c^{-}\|_{\dot B^{\frac N2-1}_{2,1}}\,d\tau\\
&\quad\lesssim X^{2}(2).
\end{split}
\end{equation}
For the term with $\mu^{+}h_{+}(c^{+},c^{-})\partial_{j}c^{+}\partial_{j}u^{+}_{i}$, we obtain
\begin{equation}
\begin{split}\label{s23step2}
&\int_0^t\|\mu^{+}h_{+}(c^{+},c^{-})\partial_{j}c^{+}\partial_{j}u^{+}_{i}(\tau)
\|_{\dot B^{\frac N2-1}_{2,1}}\,d\tau\\
&\quad\lesssim\int_0^t\|h_{+}(c^{+},c^{-})\|_{\dot B^{\frac N2}_{2,1}}
\|\nabla c^{+}\nabla u^{+}(\tau)\|_{\dot B^{\frac N2-1}_{2,1}}\,d\tau\\
&\quad\lesssim\int_0^t\big(\|h_{+}(c^{+},c^{-})-\frac{(\mathcal{C}^{2}\alpha^{-})(1,1)}
{s_{-}^{2}(1,1)}\|_{\dot B^{\frac N2}_{2,1}}+1\big)
\|\nabla c^{+}\|_{\dot B^{\frac N2-1}_{2,1}}
\|\nabla u^{+}(\tau)\|_{\dot B^{\frac N2}_{2,1}}\,d\tau\\
&\quad\lesssim\big(1+\|(c^{+},c^{-})\|_{L^{\infty}_{t}{(\dot B^{\frac N2-1 \frac N2}_{2,1})}}\big)
\|c^{+}\|_{L^{\infty}_{t}{(\dot B^{\frac N2-1 \frac N2}_{2,1})}}
\| u^{+}\|_{L^{1}_{t}{(\dot B^{\frac N2+1}_{2,1})}}\,d\tau\\
&\quad\lesssim X^{2}(2)+X^{3}(2).
\end{split}
\end{equation}
The estimates of other terms such as
$\mu^{+}k_{+}(c^{+},c^{-})\partial_{j}c^{-}\partial_{j}u^{+}_{i},~
\mu^{+}h_{+}(c^{+},c^{-})\partial_{j}c^{+}\partial_{i}u^{+}_{j},~
\mu^{+}k_{+}(c^{+},c^{-})\partial_{j}c^{-}\partial_{i}u^{+}_{j},\\
\lambda^{+}h_{+}(c^{+},c^{-})\partial_{i}c^{+}\partial_{j}u^{+}_{j}~\text{and}~
\lambda^{+}k_{+}(c^{+},c^{-})\partial_{i}c^{-}\partial_{j}u^{+}_{j}$
are similar to \eqref{s23step2}. Here, we omit the details.

For $\mu^{+}l_{+}(c^{+},c^{-})\partial_{j}^{2}u_{i}^{+},$ thus
\begin{equation}
\begin{split}\label{s2step2last2}
&\int_0^t\|\mu^{+}l_{+}(c^{+},c^{-})\partial_{j}^{2}u_{i}^{+}(\tau)\|_{\dot B^{\frac N2-1}_{2,1}}\,d\tau\\
&\quad\lesssim\|(c^{+},c^{-})\|_{L^{\infty}_{t}{(\dot B^{\frac N2}_{2,1})}}
\|\nabla^{2} u^{+}\|_{L^{1}_{t}{(\dot B^{\frac N2-1}_{2,1})}}\,d\tau\\
&\quad\lesssim\|(c^{+},c^{-})\|_{L^{\infty}_{t}{(\dot B^{\frac N2-1 \frac N2}_{2,1})}}
\| u^{+}\|_{L^{1}_{t}{(\dot B^{\frac N2+1}_{2,1})}}\,d\tau\\
&\quad\lesssim X^{2}(2),
\end{split}
\end{equation}
where $l_{+}$ stands for some smooth function vanishing at $0$.

Similarly,
\begin{align*}
\int_0^2\| H_{21}\|_{\dot B^{\frac N2-1}_{2,1}}^{h}d\tau
&\lesssim X^{2}(2)+X^{3}(2),
\end{align*}
\begin{align*}
\int_0^2\| H_{41}\|_{\dot B^{\frac N2-1}_{2,1}}^{h}d\tau
&\lesssim X^{2}(2)+X^{3}(2),
\end{align*}
and
\begin{align*}
&\int_0^2\|\nabla (u^{+},u^{-})\|_{\dot B^{\frac N2}_{2,1}}\|(c^{+},u^{+},c^{-},u^{-},\nabla c^{+},\nabla c^{-})\|_{\dot B^{\frac N2-1}_{2,1}}d\tau\\
&\quad\lesssim\|\nabla (u^{+},u^{-})\|_{L_{T}^{1}(\dot B^{\frac N2}_{2,1})}
\|(c^{+},u^{+},c^{-},u^{-},\nabla c^{+},\nabla c^{-})\|_{L_{T}^{\infty}(\dot B^{\frac N2-1}_{2,1})}\\
&\quad\lesssim\|(u^{+},u^{-})\|_{L_{T}^{1}(\dot B^{\frac N2+1}_{2,1})}
\big(\|(c^{+},c^{-})\|_{L_{T}^{\infty}(\dot B^{\frac N2-1 \frac N2}_{2,1})}
+\|(u^{+},u^{-})\|_{L_{T}^{\infty}(\dot B^{\frac N2-1}_{2,1})}\big)\\
&\quad\lesssim X^{2}(2).
\end{align*}
Therefore, for the case $t\leq2$,
\begin{equation}\label{l2}
\sum_{q\geq q_0}\sup_{0\leq t\leq 2}\Big(
\langle t\rangle^\alpha\!\int_0^t\!e^{c_0(\tau-t)}2^{q(\frac{N}{2}-1)}S_q\,d\tau\Big)
\lesssim X^2(2)+X^3(2).
\end{equation}
To bound the supremum on $[2,T],$ we split the integral on $[0,t]$ into  integrals
on  $[0,1]$ and $[1,t],$ respectively.
The $[0,1]$ part of the integral is easy to handle, and we have
$$\begin{aligned}
&\sum_{q\geq q_0}\sup_{2\leq t\leq T}\Big(
\langle t\rangle^\alpha
\!\int_0^1\!e^{c_0(\tau-t)}2^{q(\frac{N}{2}-1)}S_q(\tau)\,d\tau\Big)\\
&\quad\leq  \sum_{q\geq q_0}\sup_{2\leq t\leq T} \Big(\langle t\rangle^\alpha e^{-\frac{c_0}2t} \int_0^1 2^{q(\frac{N}{2}-1)}S_q\,d\tau\Big)\\
&\quad\lesssim \int_0^1  \sum_{q\geq q_0}  2^{q(\frac{N}{2}-1)}S_q\,d\tau.
\end{aligned}
$$
Hence
\begin{equation}\label{l1}
\sum_{q\geq q_0}\sup_{2\leq t\leq T}
\Big(\langle t\rangle^\alpha\!\int_0^1\!e^{c_0(\tau-t)}2^{q(\frac{N}{2}-1)}S_q(\tau)\,d\tau\Big)\lesssim X^2(1)+X^3(1).
\end{equation}
Let us finally consider the $[1,t]$ part of the integral for $2\leq t\leq T.$  We shall use repeatedly the following inequalities
\begin{equation}\label{eq:135}
\|\tau\nabla u^{\pm}\|_{\wt L^\infty_t(\dot{B}^{\frac N2}_{2,1})}\lesssim D(t),
\end{equation}
which is straightforward as regards to  the high frequencies of $u^{\pm}$ and  stem from
$$
\|\tau\nabla u^{\pm}\|_{\wt L^\infty_t(\dot B^{\frac N2}_{2,1})}^\ell\lesssim
\|\langle\tau\rangle^{\frac N4+\frac 12}\nabla u^{\pm}\|_{\wt L^\infty_t(\dot B^{\frac N2}_{2,1})}^\ell
\lesssim\|\langle\tau\rangle^{\frac N4+\frac 12} u^{\pm}\|_{L^\infty_t(\dot B^{1}_{2,1})}^\ell\leq D(t),
$$
for the low frequencies of  $u^{\pm}$.

Regarding the contribution of $S_q^1,$  by Lemma \ref{lemma2.13} we first notice that
\begin{align*}\label{eq:decay6}
&\sum_{q\geq q_0}\sup_{2\leq t\leq T}\Big(\langle t\rangle^\alpha\!\int_1^t\!e^{c_0(\tau-t)}2^{q(\frac{N}{2}-1)}S_q^1(\tau)\,d\tau\Big)\\
&=\sum_{q\geq q_0}\sup_{2\leq t\leq T}\Big(
\langle t\rangle^\alpha\!\int_1^t\!e^{c_0(\tau-t)}2^{q(\frac{N}{2}-1)}
\|(\dot{\Delta}_qH_{11},\,\dot{\Delta}_qH_{21},
\,\dot{\Delta}_qH_{31},\,\dot{\Delta}_qH_{41},\,\nabla\dot{\Delta}_qH_{11},
\,\nabla\dot{\Delta}_qH_{31},
)\|_{L^2}(\tau)\,d\tau\Big)\\
&\lesssim \|\tau^\alpha(H_{11},H_{21},
H_{31},H_{41},\,\nabla H_{11},
\,\nabla H_{31},)\|_{\wt L^\infty_T(\dot B^{\frac N2-1}_{2,1})}^h\\
&\lesssim \|\tau^\alpha(\nabla H_{11},H_{21},\nabla H_{31}, H_{41})\|_{\wt L^\infty_T(\dot B^{\frac N2-1}_{2,1})}^h.
\end{align*}
Now, the product laws in tilde spaces ensures that
$$
\|\tau^\alpha\nabla  H_{11}\|_{\wt L^\infty_T(\dot B^{\frac N2-1}_{2,1})}^h
\lesssim\|\tau^{\alpha-1} c^{+}\|_{\wt L^\infty_T(\dot B^{\frac N2}_{2,1})}
\|\tau\div u^{+}\|_{\wt L^\infty_T(\dot B^{\frac N2}_{2,1})}.
$$
The high frequencies of the first term is obviously bounded by $D(T)$. That is,
\begin{equation}
\label{ch32}
\|\tau^{\alpha-1}c^{+}\|^h_{\wt L_T^\infty(\dot B^{\frac N2}_{2,1})}
\lesssim \|\tau^{\alpha}c^{+}\|^h_{\wt L_T^\infty(\dot B^{\frac N2}_{2,1})}
\lesssim D(T),
\end{equation}
as for the low frequencies, we notice that if $N\leq4$ for all small enough $\varepsilon>0,$
\begin{equation}\label{cl32one}
\begin{split}
\|\tau^{\alpha-1}c^{+}\|^\ell_{\wt L_T^\infty(\dot B^{\frac N2}_{2,1})}
&\lesssim
\|\tau^{\alpha-1}c^{+}\|^\ell_{L_T^\infty(\dot B^{\frac N2-2\varepsilon}_{2,1})}\\
&\lesssim
\|\tau^{\alpha-1-\frac N2+\varepsilon}\tau^{\frac N2-\varepsilon}c^{+}\|^\ell_{L_T^\infty(\dot B^{\frac N2-2\varepsilon}_{2,1})}\\
&\lesssim
\|\tau^{\frac N2-\varepsilon}c^{+}\|^\ell_{L_T^\infty(\dot B^{\frac N2-2\varepsilon}_{2,1})}\\
&\lesssim
 D(T),
\end{split}
\end{equation}
and if $N\geq5$,
\begin{equation}\label{cl32two}
\begin{split}
\|\tau^{\alpha-1}c^{+}\|^\ell_{\wt L_T^\infty(\dot B^{\frac N2}_{2,1})}
&\lesssim
\|\tau^{\alpha-1}c^{+}\|^\ell_{L_T^\infty(\dot B^{2}_{2,1})}\\
&\lesssim
\|\tau^{\alpha-2-\frac N4}\tau^{\frac N4+1}c^{+}\|^\ell_{L_T^\infty(\dot B^{2}_{2,1})}\\
&\lesssim
\|\tau^{\frac N4+1}c^{+}\|^\ell_{L_T^\infty(\dot B^{2}_{2,1})}\\
&\lesssim D(T).
\end{split}
\end{equation}
Combining with \eqref{ch32}, \eqref{cl32one} and \eqref{cl32two}, we obtain
\begin{align}
\label{c32}
\|\tau^{\alpha-1}c^{+}\|_{\wt L_T^\infty(\dot B^{\frac N2}_{2,1})}
&\lesssim D(T).
\end{align}
Therefore, using \eqref{eq:135} and  \eqref{c32} we get
$$
\|\tau^\alpha\nabla H_{11}\|_{\wt L^\infty_T(\dot B^{\frac N2-1}_{2,1})}^h\lesssim D^2(T).
$$
Similarly,
$$
\|\tau^\alpha\nabla H_{31}\|_{\wt L^\infty_T(\dot B^{\frac N2-1}_{2,1})}^h\lesssim D^2(T).
$$
Next, we shall use repeatedly the following inequality
\begin{equation}
\begin{split}
\label{notc32}
\|(c^{+},c^{-})\|_{\wt L^\infty_T(\dot B^{\frac N2}_{2,1})}
&\lesssim\|(c^{+},c^{-})\|_{\wt L^\infty_T(\dot B^{\frac N2-1,\frac N2}_{2,1})}\\
&\lesssim X(T).
\end{split}
\end{equation}
For the first term of $H_{21}$, employing \eqref{notc32} and  the definition of $D(t),$ we have
\begin{align*}
\|\tau^\alpha g_{+}(c^{+},c^{-})\partial_{i}(c^{+})^{h}\|_{\wt L^\infty_T(\dot B^{\frac N2-1}_{2,1})}
&\lesssim \|(c^{+},c^{-})\|_{\wt L^\infty_T(\dot B^{\frac N2}_{2,1})}
\|\tau^\alpha \nabla c^{+}\|_{\wt L^\infty_T(\dot B^{\frac N2-1}_{2,1})}^h\\
&\lesssim X(T)D(T).
\end{align*}
According to   \eqref{cl32one}, \eqref{cl32two} and  the definition of $D(t),$  we have
\begin{align*}
\|\tau^\alpha g_{+}(c^{+},c^{-})\partial_{i}(c^{+})^\ell\|_{\wt L^\infty_T(\dot B^{\frac N2-1}_{2,1})}
&\lesssim \|\tau^{1-\varepsilon} (c^{+},c^{-})\|_{\wt L^\infty_T(\dot B^{\frac N2}_{2,1})}\|\tau^{\alpha-1+\varepsilon} \nabla c^{+}\|^\ell_{\wt L^\infty_T(\dot B^{\frac N2-1}_{2,1})}\\
&\lesssim \| \tau^{1-\varepsilon} (c^{+},c^{-})\|_{\wt L^\infty_T(\dot B^{\frac N2}_{2,1})}
D(T)\\
&\lesssim \Big(\|\tau^{1-\varepsilon} (c^{+},c^{-})\|_{ L^\infty_T(\dot B^{\frac N2-2\varepsilon}_{2,1})}^{\ell}
+\|\tau (c^{+},c^{-})\|_{\wt L^\infty_T(\dot B^{\frac N2}_{2,1})}^{h}\Big)
D(T)\\
&\lesssim \Big(\|\tau^{\frac{N}{2}-\varepsilon} (c^{+},c^{-})\|_{ L^\infty_T(\dot B^{\frac N2-2\varepsilon}_{2,1})}^{\ell}
+\|\langle \tau\rangle^{\alpha} (c^{+},c^{-})\|_{\wt L^\infty_T(\dot B^{\frac N2}_{2,1})}^{h}\Big)
D(T)\\
&\lesssim  D^2(T).
\end{align*}
Thus,
\begin{equation}
\label{s21step223}
\begin{split}
\|\tau^\alpha g_{+}(c^{+},c^{-})\partial_{i}c^{+} \|^{h}_{\wt L^\infty_T(\dot B^{\frac N2-1}_{2,1})}
\lesssim X^{2}(T)+ D^2(T).
\end{split}
\end{equation}
Similarly,
\begin{equation*}
\begin{split}
\|\tau^\alpha \tilde{g}_{+}(c^{+},c^{-})\partial_{i}c^{-} \|^{h}_{\wt L^\infty_T(\dot B^{\frac N2-1}_{2,1})}
\lesssim X^{2}(T)+ D^2(T).
\end{split}
\end{equation*}
Employing Propositions \ref{p26}-\ref{p27} ,  \eqref{eq:135} , \eqref{c32} and \eqref{notc32}, for the term $\mu^{+}h_{+}(c^{+},c^{-})\partial_{j}c^{+}\partial_{j}u^{+}_{i}$, we have
\begin{equation}
\begin{split}
\label{L31}
&\Big\| \tau^\alpha \mu^{+}h_{+}(c^{+},c^{-})\partial_{j}c^{+}\partial_{j}u^{+}_{i}\big)\Big\|_{\wt L^\infty_T(\dot B^{\frac N2-1}_{2,1})}^{h}\\
&\quad\lesssim\|h_{+}(c^{+},c^{-})\|_{\wt L^\infty_T(\dot B^{\frac N2}_{2,1})}
\| \tau^\alpha \nabla c^{+} \nabla u^{+}\|_{\wt L^\infty_T(\dot B^{\frac N2-1}_{2,1})}\\
&\quad\lesssim\Big(1+\|h_{+}(c^{+},c^{-})-\frac{(\mathcal{C}^{2}\alpha^{-})(1,1)}
{s_{-}^{2}(1,1)}\|_{\wt L^\infty_T(\dot B^{\frac N2}_{2,1})}
\Big)
\| \tau^{\alpha-1} \nabla c^{+}\|_{\wt L^\infty_T(\dot B^{\frac N2-1}_{2,1})}
\| \tau\nabla u^{+}\|_{\wt L^\infty_T(\dot B^{\frac N2}_{2,1})}\\
&\quad\lesssim
\Big(1+\|(c^{+},c^{-})\|_{\wt L^\infty_T(\dot B^{\frac N2}_{2,1})}\Big)
\| \tau^{\alpha-1}  c^{+}\|_{\wt L^\infty_T(\dot B^{\frac N2-1 \frac N2}_{2,1})}
\| \tau\nabla u^{+}\|_{\wt L^\infty_T(\dot B^{\frac N2}_{2,1})}\\
&\quad\lesssim  X^{2}(T)+D^{2}(T)+D^{4}(T).
\end{split}
\end{equation}
 The terms
$\mu^{+}k_{+}(c^{+},c^{-})\partial_{j}c^{-}\partial_{j}u^{+}_{i},~
\mu^{+}h_{+}(c^{+},c^{-})\partial_{j}c^{+}\partial_{i}u^{+}_{j},~
\mu^{+}k_{+}(c^{+},c^{-})\partial_{j}c^{-}\partial_{i}u^{+}_{j},\\
\lambda^{+}h_{+}(c^{+},c^{-})\partial_{i}c^{+}\partial_{j}u^{+}_{j}~\text{and}~
\lambda^{+}k_{+}(c^{+},c^{-})\partial_{i}c^{-}\partial_{j}u^{+}_{j}$
may be treated along the same lines.

From \eqref{eq:135} and \eqref{c32}, we also see that
\begin{align*}
\|\tau^\alpha \mu^{+}l_{+}(c^{+},c^{-})\partial_{j}^{2}u_{i}^{+}\|^{h}_{\wt L^\infty_T(\dot B^{\frac N2-1}_{2,1})}
&\lesssim\|\tau\nabla^2u^{+}\|_{\wt L^\infty_T(\dot B^{\frac N2-1}_{2,1})}
\|\tau^{\alpha-1} (c^{+},c^{-})\|_{\wt L^\infty_T(\dot B^{\frac N2}_{2,1})}\\
&\lesssim D^2(T),
\end{align*}
\begin{align*}
\|\tau^\alpha (\mu^{+}+\lambda^{+})l_{+}
(c^{+},c^{-})\partial_{i}\partial_{j}u^{+}_{j}\|^{h}_{\wt L^\infty_T(\dot B^{\frac N2-1}_{2,1})}
&\lesssim\|\tau\nabla^2u^{+}\|_{\wt L^\infty_T(\dot B^{\frac N2-1}_{2,1})}
\|\tau^{\alpha-1} (c^{+},c^{-})\|_{\wt L^\infty_T(\dot B^{\frac N2}_{2,1})}\\
&\lesssim D^2(T).
\end{align*}
Thus,
$$
\|\tau^\alpha H_{21}\|_{\wt L^\infty_T(\dot B^{\frac N2-1}_{2,1})}^h\lesssim X^{2}(T)+D^{2}(T)+D^{4}(T),
$$
$$
\|\tau^\alpha H_{41}\|_{\wt L^\infty_T(\dot B^{\frac N2-1}_{2,1})}^h\lesssim X^{2}(T)+D^{2}(T)+D^{4}(T).
$$
To bound the term with $S_q^2,$ we use the following fact that
$$
\int_1^te^{c_0(\tau-t)}\|R_q(u^{+},c^{+})\|_{L^2}\,d\tau\lesssim
\|R_q(\tau u^{+},\tau^{\alpha-1}c^{+})\|_{L^\infty_t(L^2)}
\int_1^te^{c_0(\tau-t)} \tau^{-\alpha}\,d\tau.
$$
Hence, thanks to Lemma \ref{lemma2.13} and Proposition \ref{Pro:1}, we have
\begin{equation*}
\begin{split}
&\sum_{q\geq q_0}\sup_{2\leq t\leq T}
\Big( \tau^\alpha\!\int_1^t\!e^{c_0(\tau-t)}2^{q(\frac N2-1)}S_q^2(\tau)\,d\tau\Big)\\
&\quad\lesssim\sum_{q\geq q_0}\sup_{2\leq t\leq T}
\Big( \tau^\alpha\!\int_1^t\!e^{c_0(\tau-t)}2^{q(\frac N2-1)}\|R_q(u^{+},c^{+})\|_{L^2}\,d\tau\Big)\\
&\quad\lesssim \|\tau\nabla u^{+}\|_{\wt L^\infty_T(\dot B^{\frac N2}_{2,1})} \|\tau^{\alpha-1}c^{+}\|_{\wt L^\infty_T(\dot B^{\frac N2-1}_{2,1})}.
\end{split}
\end{equation*}
The first term on the right-side of the above inequality may be bounded from \eqref{eq:135},
and the high frequencies of the last one on the right-side are obviously bounded by $D(T).$
To bound the term   $\|\tau^{\alpha-1}c^{+}\|_{\wt L^\infty_T(\dot B^{\frac N2-1}_{2,1})}^\ell,$
if $N\leq6$ we have the following inequality
\begin{equation*}
\begin{split}
\label{c111one}
\|\tau^{\alpha-1}c^{+}\|_{\wt L^\infty_T(\dot B^{\frac N2-1}_{2,1})}^\ell
&\lesssim\|\tau^{\alpha-1}c^{+}\|_{L^\infty_T(\dot B^{\frac N2-1-2\varepsilon}_{2,1})}^\ell\\
&\lesssim \|\tau^{\alpha-\frac N2-\frac 12+\varepsilon}\tau^{\frac N2-\frac 12-\varepsilon}c^{+}\|_{L^\infty_T(\dot B^{\frac N2-1-2\varepsilon}_{2,1})}^\ell\\
&\lesssim D(T).
\end{split}
\end{equation*}
If $N\geq7$,
\begin{equation*}
\begin{split}
\label{c111two}
\|\tau^{\alpha-1}c^{+}\|_{\wt L^\infty_T(\dot B^{\frac N2-1}_{2,1})}^\ell
&\lesssim\|\tau^{\alpha-1}c^{+}\|_{L^\infty_T(\dot B^{2}_{2,1})}^\ell\\
&\lesssim \|\tau^{\alpha-\frac N4-2}\tau^{\frac N4+1}c^{+}\|_{L^\infty_T(\dot B^{2}_{2,1})}^\ell\\
&\lesssim \|\tau^{\frac N4+1}c^{+}\|_{L^\infty_T(\dot B^{2}_{2,1})}^\ell\\
&\lesssim D(T).
\end{split}
\end{equation*}
We eventually get
\begin{equation}
\begin{split}
\label{6.27}
\sum_{q\geq q_0}\sup_{2\leq t\leq T}
\Bigl(\tau^\alpha\!\int_1^t\!e^{c_0(\tau-t)}2^{q(\frac N2-1)}S_q^2(\tau)\,d\tau\Big)
\lesssim D^2(T).
\end{split}
\end{equation}
Similarly,
\begin{equation}
\label{6.28}
\begin{split}
\sum_{q\geq q_0}\sup_{2\leq t\leq T}
\Big(\tau^\alpha\!\int_1^t\!e^{c_0(\tau-t)}2^{q(\frac N2-1)}\big(S_q^3(\tau)+S_q^4+S_q^5(\tau)\big)\,d\tau\Big)\lesssim D^2(T).
\end{split}
\end{equation}
Finally, using product laws, \eqref{eq:135}, \eqref{c32} and Lemma \ref{lemma2.13}, we obtain
\begin{equation}\label{6.29}
\begin{split}
&\sum_{q\geq q_0}\sup_{2\leq t\leq T}
\Big( t^\alpha\!\int_1^t\!e^{c_0(\tau-t)}2^{q(\frac N2-1)}S_q^6(\tau)\,d\tau\Big)\\
&\quad\lesssim\sum_{q\geq q_0}\sup_{2\leq t\leq T}
\Big( t^\alpha\!\int_1^t\!e^{c_0(\tau-t)}2^{q(\frac N2-1)}
\|\wt R_q(u^{+},c^{+})\|_{L^2}\,d\tau\Big)\\
&\quad\lesssim \|\tau\nabla u^{+}\|_{\wt L^\infty_T(\dot B^{\frac N2}_{2,1})}\|\tau^{\alpha-1}\nabla c^{+} \|_{\wt L^\infty_T(\dot B^{\frac N2-1}_{2,1})}
\sup_{2\leq t\leq T}\Big(  t^\alpha \int_1^te^{c_0(\tau-t)}\tau^{-\alpha}\,d\tau\Big)\\
&\quad\lesssim D^2(T),
\end{split}
\end{equation}
\begin{equation}\label{6.30}
\begin{split}
\sum_{q\geq q_0}\sup_{2\leq t\leq T}
\Big( t^\alpha\!\int_1^t\!e^{c_0(\tau-t)}2^{q(\frac N2-1)}S_q^7(\tau)\,d\tau\Big)\lesssim D^2(T),
\end{split}
\end{equation}
and
\begin{equation}\label{6.31}
\begin{split}
&\sum_{q\geq q_0}\sup_{2\leq t\leq T}\Big( t^\alpha\!\int_1^t\!e^{c_0(\tau-t)}2^{q(\frac N2-1)}S_k^8(\tau)\,d\tau\Big)\\
&\lesssim\sum_{q\geq q_0}\sup_{2\leq t\leq T}\Big( t^\alpha\!\int_1^t\!e^{c_0(\tau-t)}2^{q(\frac N2-1)}
\|\nabla (u^{+},u^{-})\|_{L^\infty}\|(\nabla\ddq c^{+},\ddq u^{+},\nabla\ddq c^{-},\ddq u^{-})\|_{L^2}(\tau)\,d\tau\Big)\\
&\lesssim \|\tau\nabla (u^{+},u^{-})\|_{\wt L^\infty_T(\dot B^{\frac N2}_{2,1})}\|\tau^{\alpha-1}(\nabla c^{+},u^{+},\nabla c^{-},u^{-})\|_{\wt L^\infty_T(\dot B^{\frac N2-1}_{2,1})}
\sup_{2\leq t\leq T}\Big(  t^\alpha \int_1^te^{c_0(\tau-t)}\tau^{-\alpha}\,d\tau\Big)\\
&\lesssim D^2(T).
\end{split}
\end{equation}
Putting all the above inequalities together, we conclude that
\begin{equation}\label{6.32}
\sum_{q\geq q_0}\sup_{2\leq t\leq T}\Big(t^\alpha\!\int_1^t\!e^{c_0(\tau-t)}2^{q(\frac N2-1)}S_q(\tau)\,d\tau\Big)
\lesssim X^{2}(T)+X^{3}(T)+D^{2}(T)+D^{4}(T).
\end{equation}
Then plugging \eqref{l2}, \eqref{l1} and \eqref{6.32} into  \eqref{S7} yields
\begin{equation}
\begin{split}\label{high}
&\|\langle \tau\rangle^\alpha(\nabla c^{+}, u^{+},\nabla c^{-}, u^{-})\|_{\wt L^\infty_T(\dot B^{\frac N2-1}_{2,1})}^{h}
\\&\qquad\lesssim\|(\nabla c^{+}_{0}, u^{+}_{0},\nabla c^{-}_{0}, u^{-}_{0})\|_{\dot B^{\frac N2-1}_{2,1}}^h
+X^2(T)+X^3(T)+D^2(T)+D^4(T).
\end{split}
\end{equation}
\subsubsection*{Step 3: Decay estimates with gain of regularity  for the high frequencies of $\nabla u^{+}, \nabla u^{-}$.}
This step is devoted to bounding the last  two terms of $D(t)$. We first deal with the term $\|\tau\nabla  u^{+}\|_{\wt L^\infty_t(\dot B^{\frac{N}{2}}_{2,1})}^h$  and shall use the fact that the velocity $u^{+}$ satisfies
\begin{equation}
\begin{split}
\p_t{u}^{+}-\cA u^{+}=F:=-\beta_{1}\nabla c^{+}-\beta_{2}\nabla c^{-}+H_{2},
\end{split}
\end{equation}
where $\cA =
\nu_{1}^{+}\Delta
-\nu_{2}^{+}\nabla\textrm{div}.$\\
So,
$$
\partial_t(t\cA u^{+})-\cA(t\cA u^{+})=\cA u^{+}+t\cA F.
$$
We deduce from Remark \ref{2.14}  that
$$
\|\tau\cA u^{+}\|_{\wt L_t^\infty(\dot B^{\frac N2-1}_{2,1})}^h
\lesssim \|\cA u^{+}\|_{L_t^1(\dot B^{\frac N2-1}_{2,1})}^h +\|\tau\cA F\|_{\wt L^\infty_t(\dot B^{\frac N2-3}_{2,1})}^h,
$$
whence, using the bounds given by Theorem \ref{th:main1}, we have
\begin{equation}\label{step3}
\|\tau\nabla u^{+}\|_{\wt L_t^\infty(\dot B^{\frac N2}_{2,1})}^h
\lesssim  X(0)+\|\tau F\|_{\wt L^\infty_t(\dot B^{\frac N2-1}_{2,1})}^h.
\end{equation}
In order to bound the first two terms of $F$, we notice that, from $\alpha\geq1$ and \eqref{high}, we have
\begin{equation*}
\begin{split}
\|\tau(\nabla c^{+},\nabla c^{-})\|_{\wt L^\infty_t(\dot B^{\frac N2-1}_{2,1})}^h&\lesssim
\|\langle\tau\rangle^\alpha \nabla c^{+},\nabla c^{-})\|_{\wt L^\infty_t(\dot B^{\frac N2-1}_{2,1})}^h
\\&\lesssim X(0)+X^2(t)+X^3(t)+D^2(t)+D^4(t).\end{split}
\end{equation*}
Next, the product and composition estimates adapted to tilde spaces give
$$
\|\tau\, g_{+}(c^{+},c^{-})\partial_{i}c^{+}\|_{\wt L^\infty_t(\dot B^{\frac N2-1}_{2,1})}^h
\lesssim \|\tau^{\frac 12} (c^{+},c^{-})\|_{\wt L^\infty_t(\dot B^{\frac N2}_{2,1})}\|\tau^{\frac 12} \nabla c^{+}\|_{\wt L^\infty_t(\dot B^{\frac N2-1}_{2,1})}
\lesssim D^2(t).
$$
Similarly,  we have
$$
\|\tau\, \tilde{g}_{+}(c^{+},c^{-})\partial_{i}c^{-}\|_{\wt L^\infty_t(\dot B^{\frac N2-1}_{2,1})}^h
\lesssim D^2(t).
$$
Furthermore, from \eqref{eq:135} and   the definition of $X(t)$, we have
\begin{align*}
\|\tau \, (u^{+}\cdot\nabla)u_{i}^{+}\|_{\wt L^\infty_t(\dot B^{\frac N2-1}_{2,1})}^h
\lesssim \|u^{+}\|_{\wt L^\infty_t(\dot B^{\frac N2-1}_{2,1})}
\|\tau \nabla u^{+}\|_{\wt L^\infty_t(\dot B^{\frac N2}_{2,1})}\lesssim X(t)D(t).
\end{align*}
Employing \eqref{eq:135} and \eqref{notc32}, we get
\begin{equation}
\label{s24step3}
\begin{split}
&\|\tau\, \mu^{+}h_{+}(c^{+},c^{-})\partial_{j}c^{+}\partial_{j}u^{+}_{i}\|_{\wt L^\infty_t(\dot B^{\frac N2-1}_{2,1})}^h\\
&\quad\lesssim \Big(1+\| (c^{+},c^{-})\|_{\wt L^\infty_t(\dot B^{\frac N2}_{2,1})}\Big)
\|\nabla c^{+}\|_{\wt L^\infty_t(\dot B^{\frac N2-1}_{2,1})}
\|\tau \nabla u^{+}\|_{\wt L^\infty_t(\dot B^{\frac N2}_{2,1})}\\
&\quad\lesssim X^{2}(t)+D^{2}(t)+X^{4}(t).
\end{split}
\end{equation}
The terms
$\mu^{+}k_{+}(c^{+},c^{-})\partial_{j}c^{-}\partial_{j}u^{+}_{i},~
\mu^{+}h_{+}(c^{+},c^{-})\partial_{j}c^{+}\partial_{i}u^{+}_{j},~
\mu^{+}k_{+}(c^{+},c^{-})\partial_{j}c^{-}\partial_{i}u^{+}_{j},\\
\lambda^{+}h_{+}(c^{+},c^{-})\partial_{i}c^{+}\partial_{j}u^{+}_{j}~\text{and}~
\lambda^{+}k_{+}(c^{+},c^{-})\partial_{i}c^{-}\partial_{j}u^{+}_{j}$
may be treated along the same lines.

From \eqref{eq:135} and   \eqref{notc32}, we have
\begin{align*}
&\|\tau\, \mu^{+}l_{+}(c^{+},c^{-})\partial_{j}^{2}u_{i}^{+}\|_{\wt L^\infty_t(\dot B^{\frac N2-1}_{2,1})}^h\\
&\lesssim \| (c^{+},c^{-})\|_{\wt L^\infty_t(\dot B^{\frac N2}_{2,1})}
 \|\tau\nabla u^{+}\|_{\wt L^\infty_t(\dot B^{\frac N2}_{2,1})}\\
&\lesssim X(t)D(t),
\end{align*}
\begin{align*}
&\|\tau\, (\mu^{+}+\lambda^{+})l_{+}
(c^{+},c^{-})\partial_{i}\partial_{j}u^{+}_{j}\|_{\wt L^\infty_t(\dot B^{\frac N2-1}_{2,1})}^h\\
&\lesssim X(t)D(t).
\end{align*}
Therefore,
\begin{equation}\label{last1}
\|\tau\nabla u^{+}\|_{\wt L_t^\infty(\dot B^{\frac N2}_{2,1})}^h\\
\lesssim X(0)+X^2(t)+X^{4}(t)+D^2(t).
\end{equation}
Similarly,
\begin{equation}\label{last2}
\|\tau\nabla u^{-}\|_{\wt L_t^\infty(\dot B^{\frac N2}_{2,1})}^h\\
\lesssim X(0)+X^2(t)+X^{4}(t)+D^2(t).
\end{equation}
Finally, adding up these obtained inequalities \eqref{last1} and  \eqref{last2} to \eqref{low} and \eqref{high} yields for all $t\geq0,$
\begin{align*}D(t)&\lesssim  X(0)+D_{0}+\|(\nabla c^{+}_{0}, u^{+}_{0},\nabla c^{-}_{0}, u^{-}_{0})\|_{\dot B^{\frac N2-1}_{2,1}}^h\\&\quad+X^{2}(t)+X^{3}(t)+X^{4}(t)+D^{2}(t)+D^{3}(t)+D^{4}(t)\\
&\lesssim  D_{0}+\|(\nabla c^{+}_{0}, u^{+}_{0},\nabla c^{-}_{0}, u^{-}_{0})\|_{\dot B^{\frac N2-1}_{2,1}}^h\\&\quad+X^{2}(t)+X^{3}(t)+X^{4}(t)+D^{2}(t)+D^{3}(t)+D^{4}(t),
\end{align*}
where we have used $X(0)^{\ell}=\|( c_{0}^{+}, u_{0}^{+}, c_{0}^{-}, u_{0}^{-})\|^{\ell}_{\dot B^{\frac N2-1}_{2,1}}\lesssim \|( c_{0}^{+}, u_{0}^{+}, c_{0}^{-}, u_{0}^{-})\|^{\ell}_{\dot B^{-\frac N2}_{2,\infty}}$.
As Theorem \ref{th:main1} ensures that $X(t)$ is small, one can conclude that \eqref{1.8} is fulfilled for all time if $D_{0}$ and $\|(\nabla R^{+}_{0},\,u^{+}_{0},\,\nabla R^{-}_{0},\,u^{-}_{0})\|^h_{\dot B^{\frac N2-1}_{2,1}}$ are small enough. This completes the proof of Theorem \ref{th:decay}.
\begin{center}

\end{center}
\end{document}